\documentclass[12pt,a4paper,twoside]{article}

\newcommand{\pageformat}[6]{\setlength{\hoffset}{-1in}
                  \setlength{\voffset}{-1in}
                  \addtolength{\hoffset}{#5}
                            \addtolength{\voffset}{#6}
                            \setlength{\oddsidemargin}{#1}
                            \setlength{\evensidemargin}{#2}
                            \setlength{\textwidth}{\paperwidth}
                  \addtolength{\textwidth}{-\oddsidemargin}
                  \addtolength{\textwidth}{-\evensidemargin}
                  \addtolength{\textwidth}{-\marginparsep}
                  \addtolength{\textwidth}{-\marginparwidth}
                            \setlength{\topmargin}{#3}
                            \setlength{\textheight}{\paperheight}
                  \addtolength{\textheight}{-\topmargin}
                  \addtolength{\textheight}{-\headheight}
                  \addtolength{\textheight}{-\headsep}
                  \addtolength{\textheight}{-\footskip}
                  \addtolength{\textheight}{-#4}}
\pageformat{2cm}{3cm}{25mm}{25mm}{1pt}{0pt}

\usepackage{ifthen}
\newboolean{article}
    \setboolean{article}{true}
\newboolean{report}
\newboolean{book}
\newboolean{letter}
\newboolean{german}
\newboolean{italian}
\newboolean{nobaselinestretch}
\newboolean{nosectionappendix}
\newboolean{oldtoc}
\newboolean{nosectionequn}
\newboolean{notheorem}

\ifthenelse{\boolean{german}}{
    \usepackage{german}}{}

\usepackage[latin1]{inputenc}
\usepackage{hyperref}
\usepackage{arttoc}

\usepackage{amsmath}
\usepackage{amssymb}
\usepackage[mathscr]{eucal}

\ifthenelse{\boolean{notheorem}}{}{
    \usepackage{theorem}}

\ifthenelse{\boolean{nobaselinestretch}}{}{
    \renewcommand{\baselinestretch}{1.25}}

\newenvironment{env}[2]{\begin{#1}#2\end{#1}}{}
    \newcommand{\beq}[1]{\begin{env}{equation}{#1}}
    \newcommand{\beqn}[1]{\begin{env}{equation*}{#1}}
    \newcommand{\bal}[1]{\begin{env}{align}{#1}}
    \newcommand{\baln}[1]{\begin{env}{align*}{#1}}
    \newcommand{\bga}[1]{\begin{env}{gather}{#1}}
    \newcommand{\bgan}[1]{\begin{env}{gather*}{#1}}
    \newcommand{\bflal}[1]{\begin{env}{flalign}{#1}}
    \newcommand{\bflaln}[1]{\begin{env}{flalign*}{#1}}
    \newcommand{\bmu}[1]{\begin{env}{multline}{#1}}
    \newcommand{\bmun}[1]{\begin{env}{multline*}{#1}}
    \newcommand{\bsp}[1]{\begin{env}{split}{#1}}

    \newcommand{\eeq}{\end{env}}
    \newcommand{\eeqn}{\end{env}}
    \newcommand{\eal}{\end{env}}
    \newcommand{\ealn}{\end{env}}
    \newcommand{\ega}{\end{env}}
    \newcommand{\egan}{\end{env}}
    \newcommand{\eflal}{\end{env}}
    \newcommand{\eflaln}{\end{env}}
    \newcommand{\emu}{\end{env}}
    \newcommand{\emun}{\end{env}}
    \newcommand{\esp}{\end{env}}

\newcommand{\lf}{\vspace{2ex}}

\renewcommand{\bf}[1]{\textbf{#1}}
\renewcommand{\it}[1]{\textit{#1}}
\renewcommand{\sc}[1]{\textsc{#1}}
\renewcommand{\sf}[1]{\textsf{#1}}

\renewcommand{\tt}[1]{\texttt{#1}}
\newcommand{\hl}[1]{\bf{\it{#1}}}
\newcommand{\mrm}[1]{\mathrm{#1}}
\newcommand{\mbf}[1]{\mathbf{#1}}
\newcommand{\msf}[1]{\text{\small $\sf{#1}$}}

\newcommand{\cmc}[1]{\mathcal{#1}}
\newcommand{\eus}[1]{\mathscr{#1}}
\newcommand{\euf}[1]{\mathfrak{#1}}
\newcommand{\bb}[1]{\mathbb{#1}}

\newcommand{\mfootnotesize}[1]{{\setlength{\arraycolsep}{.5ex}\text{\footnotesize$#1$}}}
\newcommand{\mscriptsize}[1]{{\setlength{\arraycolsep}{.3ex}\text{\scriptsize$#1$}}}
\newcommand{\mtiny}[1]{{\setlength{\arraycolsep}{.3ex}\text{\tiny$#1$}}}
\newcommand{\nbd}[1]{$#1$\nobreakdash--}
\newcommand{\ol}[1]{\overline{#1}}
\newcommand{\ul}[1]{\underline{#1}}
\newcommand{\wt}[1]{\widetilde{#1}}
\newcommand{\wh}[1]{\widehat{#1}}
\newcommand{\ve}{\varepsilon}
\newcommand{\vt}{\vartheta}
\newcommand{\vk}{\varkappa}
\newcommand{\vp}{\varphi}
\newcommand{\om}{\omega}

\newcommand{\abs}[1]{\left\lvert#1\right\rvert}
\newcommand{\norm}[1]{\left\lVert#1\right\rVert}
\newcommand{\babs}[1]{\bigl\lvert#1\bigr\rvert}
\newcommand{\bnorm}[1]{\bigl\lVert#1\bigr\rVert}
\newcommand{\Babs}[1]{\Bigl\lvert#1\Bigr\rvert}

\newcommand{\snorm}[1]{\norm{\smash{#1}}}
\newcommand{\sabs}[1]{\abs{\smash{#1}}}
\newcommand{\family}[1]{\left(#1\right)}

\newcommand{\bfam}[1]{\bigl(#1\bigr)}

\newcommand{\AB}[1]{\langle#1\rangle}
\newcommand{\bAB}[1]{\bigl\langle#1\bigr\rangle}
\newcommand{\BAB}[1]{\Bigl\langle#1\Bigr\rangle}
\newcommand{\CB}[1]{\{#1\}}
\newcommand{\bCB}[1]{\bigl\{#1\bigr\}}
\newcommand{\BCB}[1]{\Bigl\{#1\Bigr\}}
\newcommand{\SB}[1]{[#1]}

\newcommand{\RO}[1]{[#1)}

\newcommand{\Matrix}[1]{\begin{pmatrix}#1\end{pmatrix}}

\newcommand{\fMatrix}[1]{\mfootnotesize{\Matrix{#1}}}
\newcommand{\sMatrix}[1]{\mscriptsize{\Matrix{#1}}}
\newcommand{\tMatrix}[1]{\mtiny{\Matrix{#1}}}

\newcommand{\sbars}[1]{\:\bar{#1}^s\:}
\newcommand{\sodots}{\sbars{\odot}}
\newcommand{\set}[2][]{
    \ifthenelse{\equal{#1}{}}{
        \CB{#2}}{
        \CB{#1~|~#2}}}
\newcommand{\bset}[2][]{
    \ifthenelse{\equal{#1}{}}{
        \bCB{#2}}{
        \bCB{#1~|~#2}}}
\newcommand{\Bset}[2][]{
    \ifthenelse{\equal{#1}{}}{
        \BCB{#2}}{
        \BCB{#1~\big|~#2}}}
\newcommand{\zero}{\CB{0}}

\DeclareMathOperator{\ls}{\normalfont\msf{span}}
\DeclareMathOperator{\cls}{\ol{\ls}}

\DeclareMathOperator{\id}{\normalfont\msf{id}}

\renewcommand{\dim}{\operatorname{\msf{dim}}}

\newcommand{\C}{\bb{C}}

\newcommand{\E}{\bb{E}}

\newcommand{\K}{\bb{K}}

\newcommand{\N}{\bb{N}}

\newcommand{\R}{\bb{R}}
\newcommand{\bS}{\bb{S}}

\newcommand{\cA}{\cmc{A}}
\newcommand{\cB}{\cmc{B}}
\newcommand{\cC}{\cmc{C}}
\newcommand{\cD}{\cmc{D}}

\newcommand{\cH}{\cmc{H}}

\newcommand{\cL}{\cmc{L}}

\newcommand{\sB}{\eus{B}}

\newcommand{\sD}{\eus{D}}

\newcommand{\sF}{\eus{F}}

\newcommand{\sK}{\eus{K}}

\newcommand{\sN}{\eus{N}}

\newcommand{\sS}{\eus{S}}

\newcommand{\ek}{\euf{k}}
\newcommand{\el}{\euf{l}}
\newcommand{\en}{\euf{n}}

\newcommand{\ep}{\euf{p}}

\newcommand{\es}{\euf{s}}

\newcommand{\eB}{\euf{B}}

\newcommand{\eH}{\euf{H}}

\newcommand{\eR}{\euf{R}}

\newcommand{\U}{\mbf{1}}
\newcommand{\F}{{\mrm{F}}}
\newcommand{\G}{\Gamma}

\newcommand{\DG}{{\mrm{I}\hspace{-0.3ex}\G}}

\newcommand{\I}{{I\!\!\!\;I}}

\newcommand{\s}{\text{\scriptsize$\sS$}}

\ifthenelse{\boolean{nosectionequn}}{}{
    \numberwithin{equation}{section}
    }

\ifthenelse{\boolean{article}\or\boolean{letter}\or\boolean{nosectionequn}}{
    \setboolean{nosectionappendix}{true}}{}
\ifthenelse{\boolean{nosectionappendix}}{}{
    \renewcommand{\appendix}{
        \chapter*{\appendixname}
        \addcontentsline{toc}{chapter}{\appendixname}
        \renewcommand{\thesection}{\Alph{section}}
        \setcounter{section}{0}}}
   
\ifthenelse{\boolean{report}\or\boolean{book}}{
    }{}

\ifthenelse{\boolean{notheorem}}{}{
        \newcommand{\mnname}{Mathematical note.}
        \newcommand{\enname}{End of the note.}
        \newcommand{\definame}{Definition.}
        \newcommand{\propname}{Proposition.}
        \newcommand{\lemname}{Lemma.}
        \newcommand{\exname}{Example.}
        \newcommand{\exername}{Exercise.}
        \newcommand{\remname}{Remark.}
        \newcommand{\obname}{Observation.}
        \newcommand{\thmname}{Theorem.}
        \newcommand{\corname}{Corollary.}
        \newcommand{\proofname}{Proof.}
    \ifthenelse{\boolean{german}}{
        \renewcommand{\mnname}{Mathematische Notiz.}
        \renewcommand{\enname}{Ende der Notiz.}
        \renewcommand{\exname}{Beispiel.}
        \renewcommand{\exername}{Übung.}
        \renewcommand{\remname}{Bemerkung.}
        \renewcommand{\obname}{Beobachtung.}
        \renewcommand{\thmname}{Satz.}
        \renewcommand{\corname}{Korollar.}
        \renewcommand{\proofname}{Beweis.}}{}
    \ifthenelse{\boolean{italian}}{
        \renewcommand{\mnname}{Nota matematica.}
        \renewcommand{\enname}{Fina della nota.}
        \renewcommand{\definame}{Definizione.}
        \renewcommand{\propname}{Proposizione.}
        \renewcommand{\exname}{Esempio.}
        \renewcommand{\exername}{Esercizio.}
        \renewcommand{\remname}{Nota.}
        \renewcommand{\obname}{Osservazione.}
        \renewcommand{\thmname}{Teorema.}
        \renewcommand{\corname}{Corollario.}
        \renewcommand{\proofname}{Dimostrazione.}

       \renewcommand{\appendixname}{Appendice}

       }{}
    \theoremheaderfont{\normalfont\bfseries}
    \theoremstyle{change}
        \theorembodyfont{\rmfamily}
            \newtheorem{emp}{}[section]
                \newcommand{\bemp}[1][]{
                    \begin{emp}\hskip-\labelsep\bf{#1}\hskip\labelsep}
                \newcommand{\eemp}{\end{emp}}
\newtheorem{itemp}[emp]{}
                \newcommand{\bitemp}[1][]{
                    \begin{itemp}\hskip-\labelsep\bf{#1}\hskip\labelsep\normalfont\itshape}
                \newcommand{\eitemp}{\end{itemp}}
            \newtheorem{mn}[emp]{\mnname}
                \newcommand{\bnm}{\begin{mn}~\begin{quotation}\renewcommand{\baselinestretch}{1}\small\noindent\ignorespaces}
                \newcommand{\enm}{\end{quotation}\hfill\bf{\enname}\end{mn}}
            \newtheorem{ex}[emp]{\exname}
                \newcommand{\bex}{\begin{ex}}
                \newcommand{\eex}{\end{ex}}
            \newtheorem{exer}[emp]{\exername}
                \newcommand{\bexer}{\begin{exer}}
                \newcommand{\eexer}{\end{exer}}
            \newtheorem{defi}[emp]{\definame}
                \newcommand{\bdefi}{\begin{defi}}
                \newcommand{\edefi}{\end{defi}}
            \newtheorem{rem}[emp]{\remname}
                \newcommand{\brem}{\begin{rem}}
                \newcommand{\erem}{\end{rem}}
            \newtheorem{ob}[emp]{\obname}
                \newcommand{\bob}{\begin{ob}}
                \newcommand{\eob}{\end{ob}}
        \theorembodyfont{\normalfont\itshape}
            \newtheorem{thm}[emp]{\thmname}
                \newcommand{\bthm}{\begin{thm}}
                \newcommand{\ethm}{\end{thm}}
            \newtheorem{prop}[emp]{\propname}
                \newcommand{\bprop}{\begin{prop}}
                \newcommand{\eprop}{\end{prop}}
            \newtheorem{cor}[emp]{\corname}
                \newcommand{\bcor}{\begin{cor}}
                \newcommand{\ecor}{\end{cor}}
            \newtheorem{lem}[emp]{\lemname}
                \newcommand{\blem}{\begin{lem}}
                \newcommand{\elem}{\end{lem}}
\newenvironment{empn}[1]{\lf\noindent\bf{#1}\ignorespaces\hskip\labelsep}{\lf}
		\newcommand{\bempn}[1]{\begin{empn}{#1}}
		\newcommand{\eempn}{\end{empn}}
		\newcommand{\bitempn}[1]{\begin{empn}{#1}\normalfont\itshape}
		\newcommand{\eitempn}{\end{empn}}
                \newcommand{\bnmn}{\begin{empn}{\mnname}~\begin{quotation}\renewcommand{\baselinestretch}{1}\small\noindent\ignorespaces}
                \newcommand{\enmn}{\end{quotation}\hfill\bf{\enname}\end{empn}}
		\newcommand{\bexn}{\begin{empn}{\exname}}
		\newcommand{\eexn}{\end{empn}}
		\newcommand{\bexern}{\begin{empn}{\exername}}
		\newcommand{\eexern}{\end{empn}}
		\newcommand{\bdefin}{\begin{empn}{\definame}}
		\newcommand{\edefin}{\end{empn}}
		\newcommand{\bremn}{\begin{empn}{\remname}}
		\newcommand{\eremn}{\end{empn}}
		\newcommand{\bobn}{\begin{empn}{\obname}}
		\newcommand{\eobn}{\end{empn}}

\newcommand{\qedsymbol}{~\rule[-0.35mm]{2mm}{2mm}}
    \newcounter{proof}[emp]
    \newenvironment{Proof}[1]{
        \vspace{1ex}
        \renewcommand{\item}[1][\stepcounter{proof}(\roman{proof})]%
            {##1\hskip\labelsep}
        \noindent\textsc{#1\hskip\labelsep}}{
        \nolinebreak\qedsymbol}
    \newcommand{\proof}[1][\proofname]{
        \begin{Proof}{#1}\ignorespaces}
    \newcommand{\qed}{\end{Proof}}
    \newcommand{\noqed}{
        \renewcommand{\qedsymbol}{}
        \end{Proof}}}
    \ifthenelse{\boolean{italian}}{
        \renewcommand{\proofname}{Dimostrazione.}}{}

\usepackage[varg]{txfonts}

\usepackage[matrix,arrow,curve]{xy}

\DeclareMathOperator{\sid}{\scriptstyle\sf{id}}

\begin{document}

\bibliographystyle{amsalpha}

\title{\vspace{-4ex}Classification of $E_0$--Semigroups\\by Product Systems}
\author{Michael Skeide\thanks{This work is supported by research funds of University of Molise and Italian MIUR.}}
\date{January 2009, this revision May 2014;\\to appear in Memoirs of the AMS}

{
\renewcommand{\baselinestretch}{1}
\maketitle

\vfill

\begin{abstract}
\noindent
In his Memoir from 1989, Arveson started the modern theory of product systems. More precisely, with each \nbd{E_0}semigroup (that is, a unital endomorphism semigroup) on $\sB(H)$ he associated a product system of Hilbert spaces (Arveson system, henceforth). He also showed that the Arveson system determines the \nbd{E_0}semigroup up to cocycle conjugacy. In three successor papers, Arveson showed that every Arveson system comes from an \nbd{E_0}semigroup. There is, therefore, a one-to-one correspondence between \nbd{E_0}semigroups on $\sB(H)$ (up to cocycle conjugacy) and Arveson systems (up to isomorphism).

In the meantime, product systems of correspondences (or Hilbert bimodules) have been constructed from Markov semigroups on general unital \nbd{C^*}algebras or on von Neumann algebras. These product systems showed to be an efficient tool in the construction of dilations of Markov semigroups to \nbd{E_0}semigroups and to automorphism groups. In particular, product systems over correspondences over commutative algebras (as they arise from classical Markov processes) or other algebras with nontrivial center, show surprising features that can never happen with Arveson systems.

A dilation of a Markov semigroup constructed with the help of a product system always acts on $\sB^a(E)$, the algebra of adjointable operators on a Hilbert module $E$. (If the Markov semigroup is on $\sB(H)$ then $E$ is a Hilbert space.) Only very recently, we showed that every product system can occur as the product system of a dilation of a nontrivial Markov semigroup. This makes it necessary to extend the theory to the relation between \nbd{E_0}semigroups on $\sB^a(E)$ and product systems of correspondences.

In these notes we present a complete theory of classification of \nbd{E_0}semigroups by product systems of correspondences. As an application of our theory, we answer the fundamental question if a Markov semigroup admits a dilation by a cocycle perturbations of noise: It does if and only if it is spatial.
\end{abstract}

\clearpage

\small
\tableofcontents
}


{\parskip0.5ex plus 0.5ex minus 0.5ex

\section{Introduction} \label{intro}

An \hl{\nbd{E_0}semigroup} is a semigroup of unital endomorphisms of a unital \nbd{*}algebra. Our scope in these notes is to present a complete classification of \nbd{E_0}semigroups on $\sB^a(E)$, the \nbd{*}algebra of all adjointable operators on some Hilbert module over a unital $C^*$ or von Neumann algebra $\cB$, by \hl{product systems} of \nbd{\cB}correspondences. This is a far-reaching generalization of Arveson's classification of \nbd{E_0}semigroups on $\sB(H)$ by product systems of Hilbert spaces (henceforth, \hl{Arveson systems}). Our motivation to have this generalization, is its (maybe, surprising) importance in the study of properly irreversible dynamical systems, that is, of (quantum, or not) \hl{Markov semigroups} that are not \nbd{E_0}semigroups, and their dilations. (A Markov semigroup is a semigroup of unital completely positive maps on a unital $C^*$ or von Neumann algebra.) For instance, as an application of our classification, we determine the structure of those Markov semigroups that admit a dilation to a \it{cocycle perturbation} of a \it{noise}.

The theory has many ramifications. (For instance, we always have to distinguish the case of unital \nbd{C^*}algebras, with \nbd{C^*}modules and \nbd{C^*}correspondences, and the case of von Neumann algebras with von Neumann (or $W^*$) modules and correspondences. Also, the approach that works in the \nbd{C^*}case \bf{and} in the von Neumann case, is quite different from Arveson's approach. Only in the von Neumann case there is a ``dual'' approach, using \it{commutants} of von Neumann correspondences, that resembles Arveson's original approach.) As we wish, to keep track with these ramifications, we split this introduction into a `historical introduction', a more detailed `motivation', a `methodological introduction', and a brief description of the `results'. Many sections have thematic introductions on their own.

\lf\noindent\bf{Historical introduction.~}
One-parameter groups of automorphisms of the algebra $\sB(H)$ of all adjointable operators on a Hilbert space $H$ are well understood. At least since the work of Wigner \cite{Wig39} (see Sections 1.1, 2.4 and 3.4 of Arveson's monograph \cite{Arv03}) it is well known that every \nbd{\sigma}weakly continuous automorphism group $\alpha=\bfam{\alpha_t}_{t\in\R}$ on $\sB(H)$ is implemented as $\alpha_t=u_t\bullet u^*_t$ where $u=\bfam{u_t}_{t\in\R}$ is a group of unitaries in $\sB(H)$. Unitary groups, in turn, are understood in terms of their generators by \it{Stone's theorem}: $u_t=e^{it\cH}$, where $\cH$ can be any (possibly unbounded) self-adjoint operator on $H$. The generator of $\alpha$ is, therefore, (modulo domain questions) a derivation given as a commutator with the Stone generator. This is also the setting of ``standard'' quantum mechanics, in which one-parameter groups of automorphisms on $\cB(H)$ are models for the time flow of a quantum mechanical system.

Generalizations in several directions are possible. Firstly, one may allow more general algebras than $\sB(H)$. More advanced examples from quantum physics require this; but also if we wish to include the ``classical'' theories, then we should allow for commutative algebras. Automorphism groups (one-parameter or not, on $C^*$ or von Neumann algebras), so-called \nbd{C^*} or \it{\nbd{W^*}dynamical systems}, are a vast area and far from being understood completely. Secondly, one may relax from groups of automorphisms to \nbd{E_0}semigroups and further to Markov semigroups. (General contraction CP-semigroups that are non-Markov do occur, too. In particular, \nbd{n}parameter semigroups have applications to commuting tuples of operators in multivariate operator theory; see, for instance, the introduction of Shalit and Solel \cite{ShaSo09}. Here, we shall completely ignore these ramifications; however, the work Shalit and Skeide \cite{ShaSk10p} is in preparation.) Markov semigroups acting on commutative algebras, indeed, correspond to classical Markovian semigroups of transition probabilities, while Markov semigroups on noncommutative $C^*$ or von Neumann algebras are models for irreversible quantum dynamics.

Powers \cite{Pow88} initiated studying \nbd{E_0}semigroups on $\sB(H)$ for the sake of classifying general unbounded derivations of $\sB(H)$ by classifying the \nbd{E_0}semigroups generated by them; and he wishes to classify \nbd{E_0}semigroups up to \it{cocycle conjugacy}. Arveson, in his Memoir \cite{Arv89}, associated with every \nbd{E_0}semigroup on $\sB(H)$ ($H$ separable and infinite-dimensional) an Arveson system, and he showed that the \nbd{E_0}semigroup is determined by its Arveson system up to cocycle conjugacy. In the three successor papers \cite{Arv90a,Arv89a,Arv90} he showed that every Arveson system is the Arveson system associated with an \nbd{E_0}semigroup. In conclusion, Arveson systems form a complete cocycle conjugacy invariant for  \nbd{E_0}semigroups on $\sB(H)$.

The construction of product systems from dynamics has been generalized in various ways. In historical order: Bhat \cite{Bha96} associated with a Markov semigroup on $\sB(H)$ an Arveson system. He did this by, first, constructing a so-called (unique!) \it{minimal dilation} of the Markov semigroup to an \nbd{E_0}semigroup and, then, taking the Arveson system of that \nbd{E_0}semigroup. Bhat and Skeide \cite{BhSk00} took a completely different approach. For each Markov semigroup on a unital $C^*$ or von Neumann algebra $\cB$ they construct directly a product system of correspondences over $\cB$. This product system allows to get quite easily the minimal dilation. This dilation acts on the algebra $\sB^a(E)$ of adjointable operators on a Hilbert \nbd{\cB}module $E$. In Skeide \cite{Ske02} we presented the first construction of a product system (of \nbd{\cB}correspondences) from \nbd{E_0}semigroups on $\sB^a(E)$. The construction from \cite{Ske02} has been refined in various ways, which will be addressed throughout these notes. In this light, Bhat's Arveson system would generalize to the product system of \nbd{\cB}correspondences constructed from the \nbd{E_0}semigroup dilating a Markov semigroup that acts on $\sB^a(E)$, not on $\cB$.

Arveson systems classify \nbd{E_0}semigroups on $\sB(H)$ one-to-one up to cocycle conjugacy. However, one should note that this depends on the hidden assumption that $H$ is infinite-dimen\-sion\-al and separable; that is, the $\sB(H)$ are all the same. For the classification of \nbd{E_0}semigroups on $\sB^a(E)$ there are only quite rudimentary results. We know from \cite{Ske02} that two \nbd{E_0}semi\-groups acting on the same $\sB^a(E)$ are cocycle \it{equivalent} if and only if they have the same product system.\footnote{\label{FNSke02} %
Actually, \cite[Theorem 2.4]{Ske02} is in the situation when $E$ has a unit vector. But a quick look at that proof shows that, actually, this proof is not using the unit vector but only what will be called a \it{left dilation} after Convention \ref{fconv}. The simple generalization is formulated in Theorem \ref{ucethm} and furnished with a streamlined proof. 
}
As $\sB^a(E)$, unlike $\sB(H)$, allows in general for more than just inner automorphisms, cocycle \it{equivalence} and cocycle \it{conjugacy} need not coincide. (In fact, we shall see that cocycle conjugacy amounts to \it{Morita equivalence} of the product systems; see Theorem \ref{conjthm}.) The question when \nbd{E_0}semigroups acting on different $\sB^a(E)$ have the same (or possibly only Morita equivalent) product systems is open. Also, it is not known if the cocycle in \cite{Ske02} can be chosen continuous. In Skeide \cite{Ske07} we constructed for each continuous full product system of \nbd{C^*}correspondences over a unital \nbd{C^*}algebra a strongly continuous \nbd{E_0}semigroup. The case of von Neumann (or $W^*$) correspondences is completely open.

Arveson's theory reduces the classification of \nbd{E_0}semigroups up to cocycle conjugacy to the classification of Arveson systems, and we will answer the question to what extent product systems classify \nbd{E_0}semigroups on $\sB^a(E)$. In a first step, Arveson system are classified into \it{type I}, \it{II}, and \it{III}, depending on how many \it{units} they have. In a second step, type I and II systems (the so-called \it{spatial} ones) obtain a numerical \it{index}. Type I systems are Fock spaces and classified completely by their index; see Arveson \cite{Arv89}. We should say that, although the classification of Arveson systems has made notable progress starting with Tsirelson's preprints \cite{Tsi00p1,Tsi00p2} in 2000, which led to new large classes of type II and III Arveson systems, it seems hopeless to classify just those new classes. For general product systems of correspondences the situation is even worse. Type I and II systems need no longer be spatial, and only spatial product systems get a (no longer numerical) index; only spatial type I systems are Fock and classified by their index; see Skeide \cite{Ske06d} for spatial product systems and their index and Bhat, Liebscher and Skeide \cite{BLS10} for the counterexamples illustrating that spatiality is necessary. It is unclear which of the known constructions of non-type I Arveson systems generalize in a substantial way to modules. (`Substantial' means examples not obtained by tensoring a Fock product system with a non-type I Arveson system.) On the other hand, we know from Fagnola, Liebscher, and Skeide \cite{FLS05p} that product systems of classical Markov semigroups like those of Brownian motion and of Ornstein-Uhlenbeck processes are type III. Much has to be done in the classification of product systems; but this comes after the classification of \nbd{E_0}semigroups by product systems, to which we restrict in these notes.

\lf\noindent\bf{Motivation.~}
We think it is worth to say a few words, why it is important---beyond the scope of mere generalization---to push forward Arveson's theory to \nbd{E_0}semigroups on $\sB^a(E)$. Actually, we are interested in \it{proper} irreversible dynamics, and its dilations to \nbd{E_0}semigroups and further to automorphism groups. Although \nbd{E_0}semigroups \bf{are} Markov semigroups and irreversible unless they are automorphism semigroups, we do not really consider them typical examples of irreversible dynamics. For instance, faithful \nbd{E_0}semigroups are restrictions of automorphism semigroups on $\sB(H)$ to subalgebras. (Arveson and Kishimoto \cite{ArKi92} showed this for von Neumann algebras. In Skeide \cite{Ske11a}, we proved this for (not necessarily unital) \nbd{C^*}algebras. As a corollary of Appendix B.2 we will also give a new proof of \cite{ArKi92}; see Theorem \ref{ArKiW*thm}.) In the light of \it{paired \nbd{E_0}semigroups} in Powers and Robinson \cite{PoRo89} and \it{histories} in Arveson \cite{Arv96}, \nbd{E_0}semigroups are rather something that gives to automorphism groups a \it{causal structure}. But, \hl{proper} Markov semigroups are not \nbd{E_0}semigroups. So, why insist in classifying \nbd{E_0}semigroups?

The answer is (at least) two-fold. Firstly, the classification theorems we prove will allow us later on in these notes to give a complete answer to the fundamental question, when a Markov semigroup allows for a dilation to a cocycle perturbation of a \it{noise}: Namely, if and only if it is spatial. (Spatiality, introduced for $\sB(H)$ by Arveson \cite{Arv97a} and considerably wider than Powers' definition in \cite{Pow03}, is a property of the set of CP-semigroups \it{dominated} by the Markov semigroup. It is equivalent to the property of the product system being \it{spatial}, and this property together with our classification allows to construct noise and cocycle.) Secondly, quite recently in Skeide \cite{Ske12} we could show that every Arveson system is Bhat's Arveson system of a proper Markov semigroup. Floricel \cite{Flo08} obtained this Markov semigroup by compressing the \nbd{E_0}semigroup we constructed for a given Arveson system in Skeide \cite{Ske06}, and in \cite{Ske12} we we showed it is proper. This compression to a Markov semigroup works for all our constructions of \nbd{E_0}semigroups in these notes, and we dare to conjecture that also these Markov semigroups are proper. In other words, we conjecture that product systems of \nbd{\cB}correspondences do classify in some sense (proper) Markov semigroups on $\sB^a(E)$ and that all product systems do occur in that way as product systems of proper Markov semigroups. We leave this interesting problem to future work.

\lf\noindent\bf{Methodological introduction.~}
The emergence of \it{tensor product systems} of bimodules over a unital ring $\cA$ from endomorphism semigroups on that ring is, actually, a quite simple issue. The situation gets more interesting, when $\cA$ is the algebra of \bf{all} endomorphisms of a right module $E$ over another ring $\cB$. The \nbd{\cA}\nbd{\cB}bimodule $E$ (with the natural left action of $\cA$) plays, then, the role of a Morita equivalence from $\cA$ to $\cB$. This allows to ``transform'' the product system of \nbd{\cA}bimodules into a product system of \nbd{\cB}bimodules; confer Section \ref{E0MeSEC} and Footnote \ref{repMEfn} in Section \ref{vNalgSEC}.

Let us describe these algebraic ideas in more detail. By $\bS$ we denote one of the additive semigroups $\R_+=\RO{0,\infty}$ or $\N_0=\CB{0,1,2,\ldots}$. An \hl{\nbd{E_0}semi\-group} is a semigroup $\vt=\bfam{\vt_t}_{t\in\bS}$ of unital endomorphisms $\vt_t$ of a unital \nbd{*}algebra $\cA$, fulfilling $\vt_0=id_\cA$. Every \nbd{E_0}semigroup gives rise to a \it{product system} of bimodules $E_t$ over $\cA$ under tensor product over $\cA$ in the following way: Simply put $E_t={_{\vt_t}}\cA$, that is, the right module $\cA$ with left action $a.x_t:=\vt_t(a)x_t$ of $\cA$ via $\vt_t$. Denote by $\odot$ the tensor product over $\cA$.\footnote{%
Recall that the algebraic tensor product of a right \nbd{\cA}module $E$ and a left \nbd{\cA}module $F$ can be obtained as $E\odot F:=(E\otimes F)/\CB{xa\otimes y-x\otimes ay}$. It is determined by the universal property that, with the embedding $(x,y)\mapsto x\odot y$, it turns \hl{balanced} bilinear maps $j\colon E\times F\rightarrow V$ (that is, $j(xa,y)=j(x,ay)$) into unique linear maps $\breve{j}\colon E\odot F\rightarrow V$ fulfilling $\breve{j}(x\odot y)=j(x,y)$.
}
For every $s,t\in\bS$ we define an isomorphism $u_{s,t}\colon E_s\odot E_t\rightarrow E_{s+t}$ of bimodules by $x_s\odot y_t\mapsto \vt_t(x_s)y_t$, and these isomorphisms iterate associatively. Moreover, $E_0$ is $\cA$, the \hl{trivial} \nbd{\cA}bimodule, and for $s=0$ and $t=0$ the isomorphisms $u_{0,t}$ and $u_{s,0}$ reduce to left and right action of elements of $E_0=\cA$ on $E_t$ and $E_s$, respectively. (So far, this works even if $\cA$ is just a unital ring, not necessarily a \nbd{*}algebra.) If $\cA$ is a $C^*$ or a von Neumann algebra, then each $E_t$ is also a Hilbert \nbd{\cA}module with inner product $\AB{x_t,y_t}:=x_t^*y_t$. In fact, $E_t$ with its bimodule structure is a correspondence over $\cA$ and the $u_{s,t}$ are also isometric for the tensor product of correspondences.

We see that the family $E^\odot=\bfam{E_t}_{t\in\bS}$ forms a tensor product system in the sense of Bhat and Skeide \cite{BhSk00}. We call such a product system a \hl{one-dimensional} product system, because all right \nbd{\cA}modules $E_t$ are one-dimensional. Every one-dimensional product system of \nbd{\cA}correspondences arises in that way from an \nbd{E_0}semigroup on $\cA$. The \hl{trivial} product system is that one-dimensional product system where also the left action on each $\cA_t$ is the trivial one and where the product system operation is just multiplication in $\cA$. The product system of an \nbd{E_0}semigroup $\vt$ on $\cA$ is isomorphic to the trivial one if and only if $\vt$ is a semigroup of inner automorphisms. So far, all statements in this methodological introduction are easy exercises; see \cite[Example 11.1.3]{Ske01}. In this example it is also proved that the \nbd{E_0}semigroups on $\cA$ are determined by their one-dimensional product systems up to unitary cocycle equivalence (Definition \ref{unicocdef}).

If $\cA=\sB^a(E)$ is the ($C^*$ or von Neumann) algebra of all adjointable mappings on a Hilbert \nbd{\cB}module $E$, and if all $\vt_t$ are sufficiently continuous (strict in the \nbd{C^*}case and normal in the von Neumann case), then the situation gets more interesting. The representation theory of $\sB^a(E)$ asserts that for each $\vt_t$ there is a \hl{multiplicity correspondence} $E_t$ over $\cB$ and an identification $E=E\odot E_t$ in such a way that $\vt_t(a)=a\odot\id_t$ for all $a\in\sB^a(E)$; see Theorem \ref{strirepthm}. Moreover, the $E_t$ compose associatively as $E_s\odot E_t=E_{s+t}$ under tensor product (of $C^*$ or of von Neumann correspondences over $\cB$), that is, they form a product system $E^\odot$; see Section \ref{E0cocSEC}.

The representation theory, \it{cum grano salis}, may be viewed as an operation of \it{Morita equivalence} from the correspondence $_{\vt_t}\sB^a(E)$ over $\sB^a(E)$ to the correspondence $E_t$ over $\cB$, where the Hilbert \nbd{\cB}module $E$ plays, again \it{cum grano salis}, the role of the Morita equivalence from $\sB^a(E)$ to $\cB$; again, see Footnote \ref{repMEfn} in Section \ref{vNalgSEC}. In this sense, that is \it{cum grano salis}, the product system $E^\odot$ of correspondences over $\cB$ is \it{Morita equivalent} to the one-dimensional product system $\bfam{{_{\vt_t}}\sB^a(E)}_{t\in\bS}$ of correspondences over $\sB^a(E)$. (If $E$ is a full Hilbert \nbd{\cB}module, then it is a Morita equivalence from the \it{compacts} $\sK(E)$ to $\cB$, and, by Kaspaprov \cite{Kas80}, $\sB^a(E)$ is only the \it{multiplier algebra} of $\sK(E)$. This is what \it{cum grano salis} is referring to. In the von Neumann case the statements are exact.)

As just explained, each strict (or normal) \nbd{E_0}semigroup $\vt$ on $\sB^a(E)$ where $E$ is a Hilbert (or von Neumann) \nbd{\cB}module gives rise to a product system $E^\odot$ of correspondences $E_t$ over $\cB$. We wish to classify \nbd{E_0}semigroups by their product systems, so, the first thing is to look for obvious conditions product systems coming from \nbd{E_0}semigroups do fulfill. There is no harm in assuming $E$ is full, and under this assumption the isomorphism $E=E\odot E_t$ tells us that $E_t$ must be full, too.%
\footnote{Recall that $E$ is \hl{full} if the range ideal $\cB_E:=\cls\AB{E,E}$ of $E$ coincides with $\cB$. (For \hl{strongly full}, take strong closure.) If $E$ is not full, then for $t>0$, $E_t$ from the representation theory has to be full over $\cB_E$, while $E_0$ is put to be $\cB$ by hand. If we expect expect continuity properties to hold, then it is unreasonable not to restrict to full $E$, as for nonfull $E$ the product system will not have enough sections being continuous as $t=0$.}%
~So, we wish to classify \nbd{E_0}semigroups on $\sB^a(E)$ for full $E$ by full product systems. We have to do two things. Firstly, showing that every (strongly) full product system comes from an \nbd{E_0}semigroup. Secondly, indicating the equivalence relation induced by product systems on the class of \nbd{E_0}semigroups.

For the first issue, existence of \nbd{E_0}semigroups for product systems, it turns out that in all relevant cases (and most not so relevant cases), being full (respectively, being strongly full in the von Neumann or \nbd{W^*}case) is the only condition a product system must fulfill in order to be the product system of an \nbd{E_0}semigroup. Some results were known. But we add here the strongly continuous von Neumann case, which is the most difficult and in some sense also the most important case. We should say that in these notes we worry mainly about the \hl{continuous time case} $\bS=\R_+$. The \hl{discrete case} $\bS=\N_0$ has been settled in Skeide \cite{Ske09} (preprint 2004); it is a major ingredient for the solution of the continuous time case.

We recall that for Hilbert spaces, existence of \nbd{E_0}semigroups has been settled by Arveson \cite{Arv90a,Arv89a,Arv90}. Approximately 15 years later, the result has been reproved by Liebscher \cite{Lie09} (preprint 2003). Liebscher also showed that the \nbd{E_0}semigroup may be chosen \it{pure}. In Skeide \cite{Ske06}, we provided a short and elementary proof, and shortly later Arveson \cite{Arv06} presented a different short construction; in \cite{Ske06a} we showed that his construction leads to an \nbd{E_0}semigroup unitarily equivalent (that is, \it{conjugate}) to that in \cite{Ske06}. On principal reasons, Arveson's first proof \cite{Arv90} and Liebscher's proof \cite{Lie09} cannot work for modules. The simple algebraic idea of the proof in \cite{Ske06} generalizes easily to the algebraic situation (without continuity conditions) and resolves it if full generality; we repeat this idea, basically Equation \eqref{idea}, briefly in Section \ref{vNalgSEC}. Although our idea from \cite{Ske06} would work, in principle, also under continuity conditions, Arveson's method \cite{Arv06} runs much more smoothly; we use only his version.

The solution of the second issue, when do two \nbd{E_0}semigroups have the same product system?, which we present here, is surprisingly simple. We know when two \nbd{E_0}semigroups acting on the \bf{same} $\sB^a(E)$ have the same product system, namely, if and only if they are cocycle equivalent. When they act on possibly different $\sB^a(E)$, the problem is how to compare those different $\sB^a(E)$. Brown, Green, and Rieffel \cite{BGR77} have resolved a comparable problem for the \nbd{C^*}algebra $\sK(E)$ of \it{compact} operators on Hilbert modules $E$. If $E$ is full over a \nbd{\sigma}unital \nbd{C^*}algebra $\cB$ and countably generated, then the column space $E^\infty$ is isomorphic (as a Hilbert \nbd{\cB}module) to the \it{standard} Hilbert \nbd{\cB}module $\cB^\infty$. So, for all such $E$, we have $\sK(E)\otimes\sK=\sK(E^\infty)\cong\sK(\cB^\infty)=\cB\otimes\sK$ (where $\sK:=\sK(\K)$ denotes the compact operators on the infinite-dimensional separable Hilbert space $\K:=\C^\infty$). In other words, all such $\sK(E)$ are \it{stably isomorphic}. Consequently, for all such $E$, the algebras $\sB^a(E^\infty)$ are \it{strictly} isomorphic (the isomorphisms sends the compacts onto the compacts), or in other words, all such $\sB^a(E)$ are what we shall call \it{stably} strictly isomorphic. Strict \nbd{E_0}semigroups on $\sB^a(E)$ may be \it{amplified} to strict \nbd{E_0}semigroups on $\sB^a(E^\infty)=\sB^a(\cB^\infty)$. There, we may ask if two of these \nbd{E_0}semigroups are cocycle equivalent. And that is what we do: We classify (strict) \nbd{E_0}semigroups on $\sB^a(E)$ (with possibly different full countably generated Hilbert modules $E$ over a fixed \nbd{\sigma}unital \nbd{C^*}algebra $\cB$) by full product systems up to \it{stable} unitary cocycle equivalence.

Recall that we are interested rather in unital $\cB$. We discuss \nbd{\sigma}unital $\cB$, simply because (thanks to \cite{BGR77}) it is possible. (The existence of \nbd{E_0}semigroups in the case of nonunital $\cB$ is dealt with in Skeide \cite{Ske09a}. We do not tackle this here, though, a number of ideas from \cite{Ske09a} reappear in the von Neumann case. Passing to unital $\cB$ simplifies quite bit the proof of $E^\infty\cong\cB^\infty$ in \cite{BGR77}, but it does not free it from using Kasparov's \it{absorption theorem}.) In Section \ref{contSEC} we explain that in order to get by the construction in \cite{Ske07} from a product system $E^\odot$ an \nbd{E_0}semigroup on $\sB^a(E)$ where $E$ fulfills the countability hypotheses, it is indispensable that the product system $E^\odot$ be continuous and \it{countably generated}.

The von Neumann case runs much more smoothly, though the methods are quite different. Essentially, one needs the well-known fact that a suitable amplification of a faithful normal representation of a von Neumann algebra is unitarily equivalent to an amplification of the identity representation. Without assumptions like separable pre-duals, the cardinality of the amplification may vary. But this is a price we are happy to pay for having the classification up to stable cocycle conjugacy in full generality both with and without continuity conditions.

In either case, where there is \it{stable} cocycle \it{equivalence} for \nbd{E_0}semigroups on possibly different $\sB^a(E)$ but the modules $E$ still full over the same algebra, we can weaken that question to stable cocycle \it{conjugacy}. For conjugacy of \nbd{E_0}semigroups on $\sB^a(E)$ and on $\sB^a(F)$ it is no longer necessary that $E$ and $F$ are modules over the same algebra. (In Example \ref{nonidrem} we provide two \nbd{E_0}semigroups on the same $\sB^a(E)$ that are conjugate by a noninner automorphism of $\sB^a(E)$, but there is no unitary cocycle with respect to one semigroup such that conjugation with that cocycle gives the other semigroup.) The algebras have to be just Morita equivalent, and in a classification up to stable unitary cocycle conjugacy the product systems must be Morita equivalent via the same Morita equivalence. (Consequently, the two product systems in Example \ref{nonidrem} are Morita equivalent but not isomorphic.) In between, there is the notion of ternary isomorphism of $E$ and $F$ (which implies that the algebras are isomorphic), and the classification is by ternary isomorphic product systems.

We said we determine the class of Markov semigroups that allow for a dilation which is a cocycle perturbation of a noise. We refer to the introduction of Section \ref{HuPaSEC} for a discussion of the method.

We should say that in the von Neumann case, there is a ramification into the construction of a product system of correspondences over the commutant of $\cB$, $\cB'$, from an \nbd{E_0}semigroup. That construction resembles much more Arveson's original construction from \cite{Arv89}, but it does not work for \nbd{C^*}modules. In Arveson's approach the relation between Arveson system and \nbd{E_0}semigroup is given in terms of a \it{representation} of the Arveson system. As pointed out in Skeide \cite{Ske03c}, this connection remains true also for von Neumann modules. Alevras \cite{Ale04} investigated, as a special case, \nbd{E_0}semigroups on a fixed type II$_1$ factor $\cA$, and classified them up to unitary cocycle conjugacy by the \it{commutant} (in the sense of \cite{Ske03c}) of their product systems. (So, together with \cite{Ske03c}, Alevras case is included in \cite[Theorem 2.4]{Ske02} as $E=\cA=\sB^a(\cA)$.) No result on existence of \nbd{E_0}semigroups (or, better, of representations) is given in \cite{Ale04}. In Appendix B.2, we will prove in full generality existence of faithful nondegenerate representations for \it{faithful} product systems of von Neumann (or $W^*$) correspondences (known from \cite{Ske11a} in the \nbd{C^*}case). Faithful means that all $E_t$ have faithful left action. Faithfulness is a condition \it{dual} (under commutant) to strong fullness, and it is a necessary condition for having a faithful representation.

In a discussion of methodology the following subtlety should be mentioned. Arveson systems are \it{measurable} product systems; Alevras' product systems are \it{measurable} product systems; a general definition of \it{weakly measurable} product systems of \nbd{W^*}correspondences (with separable pre-dual) has been proposed by Muhly and Solel \cite{MuSo07}; and Hirshberg \cite{Hir04} discussed a version of countably generated \it{measurable} product systems of \nbd{C^*}correspondences. Except for our construction of an \nbd{E_0}semigroup for an Arveson system in \cite{Ske06} (where we wanted to be compatible with Arveson's original result), we always worked with \it{continuous} product systems instead of measurable ones. The basic idea how to define continuous product systems is from Skeide \cite[Section 7]{Ske03b}, and remains basically unchanged. Here we add the version for strongly continuous product systems of von Neumann (or $W^*$) correspondences. We do not see a need to work with measurable versions. Continuous versions simply work: Every continuous \nbd{E_0}semigroup gives a continuous product system; every continuous product system gives a continuous \nbd{E_0}semigroup; cocycle equivalence is in terms of continuous cocycles. Additionally, unlike with measurability, most results four continuous product systems do not require countability assumptions. (However, continuity properties do, as usual, allow to prove that our constructions preserve countability hypotheses; see Theorem \ref{existcthm} and Remark \ref{countvNrem}.)

\lf\noindent\bf{Results.~}
Morita equivalence is crucial. In Section \ref{MeSEC} we give a definition of Morita equivalence (Definition \ref{Medefi}) that, in our opinion, runs more smoothly than others, and we prove (Theorem \ref{Meinvthm}) that it is equivalent to one of the usual definitions. (We think this is \it{folklore} but we did not find any proof.) We also review the results from Muhly, Skeide and Solel \cite{MSS06} about strict representations of $\sB^a(E)$ and from Skeide \cite{Ske09} about strict isomorphisms.

In Section \ref{stabMeSEC} we explain the basic idea based on the well-known result (in the form of Lance \cite[Proposition 7.4]{Lan95}) that asserts that for a full countably generated Hilbert module $E$ over a \nbd{\sigma}unital \nbd{C^*}algebra the multiple $E^\infty$ is isomorphic to the column space $\cB^\infty$. In particular, all such modules $E$ are \it{stably} isomorphic. A version for modules over different algebras gives rise to the notion of \it{stable Morita equivalence} for Hilbert modules (Definition \ref{sMedefi}) with the basic Theorem \ref{sMethm} about when stable Morita equivalence happens. In Section \ref{terSEC} we discuss ternary isomorphisms and derive analogue results for the corresponding subclass of Morita equivalences among Hilbert modules.

Starting from Section \ref{ccSEC} we come to cocycle conjugacy of \nbd{E_0}semigroups. While for algebras isomorphic to $\sB(H)$ there is essentially only one notion of cocycle conjugacy, for general \nbd{*}algebras this is not so. In Section \ref{ccSEC} we discuss cocycle conjugacy of \nbd{E_0}semigroups in the completely algebraic setting of the beginning of the methodological introduction. In Section \ref{E0cocSEC} we switch to \nbd{E_0}semigroups on $\sB^a(E)$. We explain how the product system of an \nbd{E_0}semigroup on $\sB^a(E)$ is defined, and how it gives back the \nbd{E_0}semigroup in terms of a \it{left dilation}. We explain what it means that two \nbd{E_0}semigroups acting on the same $\sB^a(E)$ have isomorphic product systems (Proposition \ref{conisoprop} and Theorem \ref{ucethm}). In Theorem \ref{Bclassthm} we give a necessary and sufficient criterion for that a family of \nbd{E_0}semigroups acting on different full Hilbert \nbd{\cB}modules all have the same product system. This theorem is, however, not in terms of cocycle conjugacy and, therefore, gives only a partial solution to our classification problem. Section \ref{E0MeSEC} classifies \nbd{E_0}semigroups acting on strictly isomorphic $\sB^a(E)$ and $\sB^a(F)$ in terms of \it{Morita equivalence} of product systems (Theorem \ref{conjthm}) and discusses the restriction to the special case of \it{ternary equivalence}. In Section \ref{sccSEC} we apply the amplification idea from Section \ref{stabMeSEC} to classify (under suitable countability assumptions) all \nbd{E_0}semigroups in terms of \it{stable cocycle conjugacy} either by isomorphism classes, or by Morita equivalence classes, or by ternary equivalence classes of their product systems (Theorem \ref{sclthm}). Section \ref{contSEC}, finally, takes into account questions of continuity (strong continuity on the \nbd{E_0}semigroup side and continuity on the product systems side). Provided we have a \it{countably generated} product system $E^\odot$, only with continuity conditions and for unital \nbd{C^*}algebras we are able to guarantee that the \nbd{E_0}semigroup constructed for $E^\odot$ in Skeide \cite{Ske07} acts on a $\sB^a(E)$ where $E$ is countably generated. Under these hypotheses we obtain a one-to-one classification (Theorem \ref{C*contclassthm}) of \nbd{E_0}semigroups by product systems with continuity conditions. (The case of not necessarily unital \nbd{C^*}algebras can be included using the existence result in Skeide \cite{Ske09a}.)

One of the basic problems of quantum probability is to obtain a \it{dilation} of a \it{Markov semigroup} to a unitary cocycle perturbation of a \it{noise}, a so-called \it{Hudson-Parthasarathy dilation}. In Section \ref{HuPaSEC} we apply our classification results to prove existence of Hudson-Parthasarathy dilations for \it{spatial} Markov semigroups acting on a unital separable \nbd{C^*}algebra (Theorem \ref{C*HPthm}).

In the remainder we tackle the von Neumann case and resolve it, unlike the \nbd{C^*}case where countability hypotheses are required, in full generality. In Section \ref{vNalgSEC} we obtain the analogue results about algebraic classification of Sections \ref{MeSEC} -- \ref{sccSEC} and also parts of Section \ref{contSEC}. In Section \ref{vNcontSEC} we obtain the analogue results in the continuous case. Since continuity in that case means strong continuity in the von Neumann sense (which is much weaker a condition), considerable work has to be done. In particular, for the first time we give a concise definition of \it{strongly continuous product system} of von Neumann correspondences (Definition \ref{scPSdefi}) and we prove that every such (strongly full) product system is the strongly continuous product system of a strongly continuous normal \nbd{E_0}semigroup (Theorem \ref{vNexithm}). The proofs in Section \ref{vNcontSEC} are limited to those things that may be adapted easily from the corresponding \nbd{C^*}proofs. The proof of  Theorem \ref{vNexithm} requires considerable technical extra work and new ideas; it is postponed to Appendix B. We obtain results in full analogy with those of Section \ref{contSEC} and without any countability requirement. Section \ref{vNHPSEC}, finally, deals with Hudson-Parthasarathy dilations of spatial Markov semigroups on von Neumann algebras. (Also here, within Section \ref{vNHPSEC}, we limit ourselves to those parts of the proofs that generalize more or less directly from the \nbd{C^*}case. But even more new ideas and techniques are required; these are discussed in Appendix A.2.) Once more, the classification up to stable cocycle conjugacy allows to prove that a Markov semigroup on a von Neumann algebra admits a Hudson-Parthasarathy dilation if (and, in a sense, only if) it is spatial (Theorem \ref{vNHPthm}). It is here in the von Neumann case, where this result shows its full power. Indeed, while in the \nbd{C^*}case the spatial Markov semigroups are limited to a subset of the uniformly continuous Markov semigroups, in the von Neumann case most known strongly operator continuous Markov semigroups are spatial; see the introduction of Section \ref{vNHPSEC}.

Appendices A and B may be considered as the beginning of a systematic theory of strongly continuous product systems. But, the theory is developed only as far as necessary for our applications: Appendix A deals with \it{strong type I} product systems (systems generated by their strongly continuous units), and develops what is necessary to prove the results about spatial Markov semigroups. Appendix A also fills a gap in the proofs of \cite[Theorems 10.2 and 12.1]{BhSk00} (see acknowledgments), generalizing them considerably (Theorems \ref{cdilthm} and \ref{scdilthm}). As Corollary \ref{MScor} we reprove the result due to Markiewicz and Shalit \cite{MaSha10} that every weakly operator continuous and contractive CP-semigroup on a von Neumann algebra is strongly operator continuous. We also prove that spatial Markov semigroups on a von Neumann algebra have spatial GNS-systems (Theorem \ref{vNspthm}) contrasting the negative statement in the \nbd{C^*}case (see Bhat, Liebscher, and Skeide \cite{BLS10}), and we prove that elementary CP-semigroups have trivial GNS-systems (Corollary \ref{elemtrivcor}). Appendix B.1 provides the proof of existence of \nbd{E_0}semigroups in the von Neumann case (Theorem \ref{vNexithm}), utilizing also some results from Appendix A. In Appendix B.2 we prove that every \it{faithful} strongly continuous product systems admits a \it{strongly continuous nondegenerate normal representations} by operators on a Hilbert space. This result is connected rather indirectly with the classification by product systems. (See the comments in the methodological introduction.) However, it is important in the theory of product systems, and the technical tools developed in both appendices allow easily for the necessary modification of the proof of the \nbd{C^*}case in Skeide \cite{Ske11a}. Therefore, we include it here. As applications we provide a different proof of a result by Arveson and Kishimoto \cite{ArKi92} about embeddability into an automorphism group, and a new result about existence of \it{elementary dilations} of strongly continuous CP-semigroups on von Neumann algebras.

Several sections contain more detailed introductions to special aspects that we did not deal with properly in the present general introduction: Section \ref{HuPaSEC} (dilations, noises, and Hudson-Parthasarathy dilations); Section \ref{vNalgSEC} (von Neumann and \nbd{W^*}modules); Appendix B.2 (representations and commutants of product systems). Section \ref{MeSEC} (Morita equivalence and representations of $\sB^a(E)$) and Section \ref{ccSEC} (algebraic cocycle conjugacy) are introductory in their own right.
}


\lf\noindent
\bf{To the reader.~}
We wish to mention some (partly non-standard) conventions and a principle regarding the organization of these notes. We hope clarifying these will help the reader. The principle is the following: Remarks may be ignored where they stand. (They may contain additional information that is not needed anywhere in these notes, but that help to put them into perspectives or outline further developments. They also may contain information that is referenced to, later on. But it is always---we hope---safe to ignore them at first reading.) Observations, on the other hand, have the rank of propositions, just that we found it more convenient to integrate the argument that proves an observation into its text. (Like propositions, observations might be there, because they are needed immediately, or because they are referenced to later on. Always, an observation is---we hope---in a place where it logically belongs to.)

As for the conventions: By \hl{semigroup} we always mean \it{one-parameter} semigroup. \hl{One-parameter semigroup} means we have a semigroup $T=\bfam{T_t}_{t\in\bS}$ over one of the additive monoids $\bS=\N_0=\CB{0,1,\ldots}$ (the \hl{discrete} case) or $\bS=\R_+=\RO{0,\infty}$ (the \hl{continuous time} case). Continuous time just refers to that the parameter $t$ is from the \it{continuum} $\R_+$, not from the discrete set $\N_0$. It does not mean that the semigroup depends continuously on $t$. If we wish to emphasize that a semigroup need not be continuous, we sometimes add the adjective \hl{algebraic}. (Be aware that in the literature there is some confusion about what \it{discrete} is referring to.) One-parameter semigroups are, actually, always assumed to be one-parameter \it{monoids} in that the family of maps $T_t$ on a set $\cB$, apart from the semigroup property $T_s\circ T_t=T_{s+t}$, also fulfills $T_0=\id_\cB$.

The sign `$\otimes$' stands exclusively for the tensor product of vector spaces (over $\C$) or spaces that have been obtained by completion or closure of the vector space tensor product. (An example is the \it{external} tensor product of Hilbert modules.) The sign `$\odot$' stands for spaces that have been obtained from the vector space tensor product by (possibly) taking a quotient. (For instance, the \it{internal} tensor product of correspondences. The element $x\odot y$ is the equivalence class $x\otimes y+\sN$ where $\sN$ is the subspace that is divided out.) This convention has proved to be very useful to produce readable and self-explanatory formulae. (For instance, the formula $(E\odot F)\otimes(E'\odot F')=(E\otimes E')\odot(F\otimes F')$ involving internal and external tensor products, makes sense via \it{canonical identification} without any further explanation.)

Homomorphism always means \nbd{*}homomorphism. Likewise, for representations.

$\sB^a(E)$ stands for the (bounded) adjointable operators on a (pre-)Hilbert module $E$. We use \it{Dirac's bra(c)ket notation}---but without brackets. That is, without mention $x\in E$ is interpreted as map $b\mapsto xb$ in $\sB^a(\cB,E)$, and $x^*\colon y\mapsto\AB{x,y}$ its adjoint. A \hl{rank-one operator} is written as $xy^*\colon z\mapsto x\AB{y,z}$.

\lf\noindent
\bf{Acknowledgments.~}
This work originates in a question we asked with Rajarama Bhat and Volkmar Liebscher in the program of a two weeks \it{Research in Pairs} project in Oberwolfach. The decisive idea, surprisingly different from what we expected in Oberwolfach, I had through a very pleasant three months stay at ISI Bangalore, for which I express my gratitude to Rajarama Bhat. I am grateful to Orr Shalit, who pointed out to us a gap in \cite{BhSk00}, and for discussions about Appendix A where that gap is fixed.

Last but, surely, not least, it is my urgent wish to express my deep gratitude to the referee for a \it{superb} job in record time! Her/his efforts improved this work considerably.

\clearpage

\section{Morita equivalence and representations}\label{MeSEC}

The relation between \nbd{E_0}semigroups and product systems goes via the representation theory of $\sB^a(E)$ for a Hilbert \nbd{\cB}module $E$. The representation theory has been discussed first in Skeide \cite{Ske02} in the case when $E$ has a \hl{unit vector} $\xi$ (that is, $\AB{\xi,\xi}=\U\in\cB$) and in Muhly, Skeide and Solel \cite{MSS06} for the general case. In particular the approach in \cite{MSS06}, a slight extension of Rieffel's \cite{Rie74} discussion of the representations of the \it{imprimitivity algebra} (that is, the finite-rank operators; see below) underlines the role played by Morita equivalence. We use this section to introduce some notation and to review the relation between Morita equivalence and the theory of strict representations of $\sB^a(E)$. The definition of Morita equivalence we use here is different from standard definitions. Although it is probably folklore that it is equivalent to standard definitions, we do not know any reference. Therefore, we include a proof (Theorem \ref{Meinvthm}).

Let $E$ be Hilbert module over a \nbd{C^*}algebra $\cB$. We say, $E$ is \hl{full} if the \hl{range ideal} $\cB_E:=\cls\AB{E,E}$ in $\cB$ coincides with $\cB$. By $\sB^a(E)$ we denote the algebra of all adjointable operators on $E$. Often, we consider an element $x\in E$ as mapping $x\colon b\mapsto xb$ from $\cB$ to $E$. The adjoint of that mapping is $x^*\colon y\mapsto\AB{x,y}$. The linear hull of the \hl{rank-one operators} $xy^*$ is the algebra $\sF(E)$ of \hl{finite-rank operators}. The completion $\sK(E)$ of $\sF(E)$ is the algebra of \hl{compact operators}. We use similar notations for operators between Hilbert \nbd{\cB}modules $E$ and $F$.

A \hl{correspondence} from a \nbd{C^*}algebra $\cA$ to a \nbd{C^*}algebra $\cB$ (or \hl{\nbd{\cA}\nbd{\cB}correspondence}) is a Hilbert \nbd{\cB}module $E$ with a left action of $\cA$ that defines a \hl{nondegenerate}(\bf{!}) representation of $\cA$ by adjointable operators on $E$. A correspondence is \hl{faithful} if the left action defines a faithful representation. The (internal) \hl{tensor product} of an \nbd{\cA}\nbd{\cB}correspondence $E$ and a \nbd{\cB}\nbd{\cC}correspondence $F$ is the unique (up to isomorphism) \nbd{\cA}\nbd{\cC}correspondence $E\odot F$ generated by elementary tensors $x\odot y$ with inner product $\AB{x\odot y,x'\odot y'}=\AB{y,\AB{x,x'}y'}$ and left action $a(x\odot y)=(ax)\odot y$.

Every \nbd{C^*}algebra $\cB$ is a \nbd{\cB}\nbd{\cB}correspondence with the natural bimodule operation and inner product $\AB{b,b'}:=b^*b'$. We refer to $\cB$ with this structure as the \hl{identity \nbd{\cB}correspondence}. For every \nbd{\cA}\nbd{\cB}cor\-re\-spond\-ence $E$ we will always identify both correspondences $E\odot\cB$ and $\cA\odot E$ with $E$ via the \hl{canonical identifications} $x\odot b\mapsto xb$ and $a\odot x\mapsto ax$, respectively. (Note that the second identification is possible only, because we require the left action to be nondegenerate. Nondegeneracy of the right action is automatic.)

\bdefi\label{Medefi}
A \hl{Morita equivalence} from $\cA$ to $\cB$ is an \nbd{\cA}\nbd{\cB}correspondence $M$ for which there exists a \nbd{\cB}\nbd{\cA}correspondence $N$ such that
\bal{\label{Medefeq}
N\odot M
&
~\cong~
\cB,
&
M\odot N
&
~\cong~
\cA,
}\eal
as correspondences over $\cB$ and over $\cA$, respectively. We call $N$ an \hl{inverse} of $M$ under tensor product.

Following Rieffel \cite{Rie74a}, two \nbd{C^*}al\-ge\-bras $\cA$ and $\cB$ are \hl{strongly Morita equivalent} if there exists an \nbd{\cA}\nbd{\cB}Morita equivalence. We use nowadays convention and speak just of \hl{Morita equivalent} \nbd{C^*}algebras.
\edefi

We observe that a Morita equivalence is necessarily faithful and full. (If $M$ is not full, then $N\odot M$ is not full, too; if $M$ is not faithful, then $M\odot N$ is not faithful, too.)
\bprop
\begin{enumerate}
\item
The correspondence $N$ in \eqref{Medefeq} is unique up to isomorphism.

\item
Morita equivalence of \nbd{C^*}algebras is an equivalence relation.
\end{enumerate}
\eprop

\proof
1.)~
Suppose $N'$ is another \nbd{\cB}\nbd{\cA}correspondence fulfilling \eqref{Medefeq}. Then $N\cong\cB\odot N\cong N'\odot M\odot N\cong N'\odot\cA\cong N'$.

2.)~
$\cB$ is a \nbd{\cB}\nbd{\cB}Morita equivalence (with $N=\cB$). So, Morita equivalence is reflexive.

If $M$ is an \nbd{\cA}\nbd{\cB}Morita equivalence, then $N$ is a \nbd{\cB}\nbd{\cA}Morita equivalence. So, Morita equivalence is symmetric.

If $M_1$ is an \nbd{\cA}\nbd{\cB}Morita equivalence (with inverse $N_1$, say) and if $M_2$ is a \nbd{\cB}\nbd{\cC}Morita equivalence (with inverse $N_2$, say), then $M_1\odot M_2$ is an \nbd{\cA}\nbd{\cC}Morita equivalence with inverse $N_2\odot N_1$. So, Morita equivalence is transitive.\qed

\bprop\label{assdiagprop}
The identifications in \eqref{Medefeq} can be chosen  such that diagrams
\beqn{
\xymatrix{
				&M\odot N\odot M	\ar[dl]	\ar[dr]	&
\\
\cA\odot M	\ar@{=}[r]	&M							&M\odot\cB	\ar@{=}[l]
}
~~~~~~~~~~~~
\xymatrix{
				&N\odot M\odot N	\ar[dl]	\ar[dr]	&
\\
\cB\odot N	\ar@{=}[r]	&N							&N\odot\cA	\ar@{=}[l]
}
}\eeqn
commute.
\eprop

\proof
Fix two isomorphisms (that is, bilinear unitaries) $u\colon N\odot M\rightarrow\cB$ and $v\colon M\odot N\rightarrow\cA$. To begin with suppose that the left diagram commutes, that is, $v(m\odot n)m'=mu(n\odot m')$ for all $m,m'\in M$ and all $n\in N$. Since $M$ is faithful, the right diagram commutes if and only if it commutes also when tensored with $M$ from the right. Evaluating the left hand path on an elementary tensor $n\odot m\odot n'\odot m'$ we find
\beqn{
n\odot m\odot n'\odot m'
\longmapsto
u(n\odot m)n'\odot m'.
}\eeqn
Evaluating the right hand path we find
\beqn{
n\odot m\odot n'\odot m'
\longmapsto
nv(m\odot n')\odot m'
~=~
n\odot v(m\odot n')m'
~=~
n\odot mu(n'\odot m').
}\eeqn
Applying the isomorphism $u$ to both elements, by bilinearity we find
\beqn{
u(u(n\odot m)n'\odot m')
~=~
u(n\odot m)u(n'\odot m')
~=~
u(n\odot mu(n'\odot m')).
}\eeqn
In conclusion: If the left diagram commutes then so does the right diagram. By symmetry, of course, also the converse statement is true.

Now suppose that the left diagram does not necessarily commute. Then, still, the map
\beqn{
w\colon
M
~=~
M\odot\cB
~\xrightarrow{~\sid_M\odot u^*~}~
M\odot N\odot M
~\xrightarrow{~v\odot\sid_M~}~
\cA\odot M
~=~
M
}\eeqn
defines an automorphism $w$ of $M$ that sends $mu(n\odot m')$ to $v(m\odot n)m'$. If we replace $v$ with $v':=v(w^*\odot\id_N)$, then the automorphism $w'\colon mu(n\odot m')\mapsto v'(m\odot n)m'$ of $M$ corresponding to the new pair $u,v'$ satisfies
\beqn{
w'((wm)u(n\odot m'))
~=~
v'((wm)\odot n)m'
~=~
v(m\odot n)m'
~=~
w(mu(n\odot m'))
~=~
(wm)u(n\odot m').
}\eeqn
Therefore, $w'$ is the identity. Equivalently, for the pair $u,v'$ the left and, therefore, both diagrams commute.\qed

\bemp[Convention.~]\label{Massconv}
After this proposition we shall always assume that the diagrams commute. This allows us to identify $\cB$ with $N\odot M$ and $\cA$ with $M\odot N$ without having to worry about brackets in tensor products.
\eemp

The following example is basic for everything about Morita equivalence.

\bex\label{fundMeex}
Every Hilbert \nbd{\cB}module $E$ may be viewed as Morita equivalence from $\sK(E)$ to $\cB_E$. In fact, the space $E^*=\CB{x^*\colon x\in E}$ becomes a correspondence from $\cB_E$ to $\sK(E)$ if we define the inner product $\AB{x^*,y^*}:=xy^*\in\sK(E)$ and the bimodule operation $bx^*a:=(a^*xb^*)^*$. Clearly, $E^*\odot E=\cB_E$ (via $x^*\odot y=\AB{x,y}$), and $E\odot E^*=\sK(E)$ (via $x\odot y^*=xy^*$). Moreover, since
\beqn{
(x\odot y^*)\odot z
~=~
(xy^*)z
~=~
x\AB{y,z}
~=~
x\odot(y^*\odot z),
}\eeqn
these identifications also satisfy Convention \ref{Massconv}.
\eex

\bcor\label{KE*cor}
$\sK(E^*)=\cB_E$ and $E^{**}=E$ as correspondence from $\sK(E)$ to $\sK(E^*)=\cB_E$.
\ecor

\proof
Since $E$ is a full Hilbert \nbd{\cB_E}module, the left action of $\cB_E$ on $E^*$ is faithful. It follows that $\AB{x,y}\mapsto x^*y^{**}$ defines an injective homomorphism from $\cB_E$ onto $\sK(E^*)$.\qed

\lf
The following result makes the connection with the definition of Morita equivalence in Lance \cite[Chapter 7]{Lan95}. It also shows that Example \ref{fundMeex} captures, in a sense, the most general situation of Morita equivalence.

\bthm\label{Meinvthm}
An \nbd{\cA}\nbd{\cB}correspondence $M$ is a Morita equivalence if and only if $M$ is full and the left action defines an isomorphism $\cA\rightarrow\sK(M)$.
\ethm

\proof
We already know that for being a Morita equivalence, $M$ must be full and faithful. So the only question is whether or not the injection $\cA\rightarrow\sB^a(M)$ is onto $\sK(M)$.

By Example \ref{fundMeex}, $M$ with the canonical action of $\sK(M)$ is a Morita equivalence from $\sK(M)$ to $\cB$. So, if the left action of $\cA$ on $M$ defines an isomorphism $\alpha\colon\cA\rightarrow\sK(M)$, then we turn the \nbd{\cB}\nbd{\sK(M)}correspondence $M^*$ into a \nbd{\cB}\nbd{\cA}correspondence $N$ with inner product $\AB{m^*,m'^*}_N:=\alpha^{-1}(mm'^*)$ and right action $m^*a=m^*\alpha(a)$. Clearly, $N$ is an inverse of $M$ under tensor product.

Conversely, suppose that $M$ is a Morita equivalence with inverse $N$, say. The idea is to establish a map $u\colon N\rightarrow M^*$ such that $\AB{n,n'}\mapsto\AB{un,un'}$ extends as an isomorphism $\alpha\colon\cA\rightarrow\sK(M)$. A look at how to resolve $N=\cB\odot N=M^*\odot M\odot N$ in the other direction to give $M^*$ by ``bringing somehow $M\odot N$ under the $*$'', reveals
\beqn{
u
\colon
\AB{m,m'}n
~\longmapsto~
((m'\odot n)^*m)^*
}\eeqn
as the only reasonable attempt. We do not worry, at that point, about whether $u$ is well-defined. (See, however, Remark \ref{urem} below.) What we wish to show is that the map
\beqn{
\alpha
\colon
\BAB{\AB{m_1,m_1'}n_1,\AB{m_2,m_2'}n_2}
~\longmapsto~
((m'_1\odot n_1)^*m_1)((m'_2\odot n_2)^*m_2)^*
}\eeqn
(the right-hand side is, $\bAB{u(\AB{m_1,m_1'}n_1),u(\AB{m_2,m_2'}n_2)}$, once $u$ showed to be well-defined) is nothing but the canonical homomorphism $\cA\rightarrow\sB^a(M)$ when applied to
\beqn{
a
~:=~
\bAB{\AB{m_1,m_1'}n_1,\AB{m_2,m_2'}n_2}
~\in~
\cA.
}\eeqn
(From this everything follows: Well-definedness, because the canonical homomorphism is well-defined. Injectivity, because $M$ is faithful. Surjectivity onto $\sK(M)$, because we obtain a dense subset of the rank-one operators.) To achieve our goal we calculate the matrix element $\AB{\wt{m}_1,a\wt{m}_2}$ and convince ourselves that it coincides with the corresponding matrix element of the operator on the right-hand side. We find
\beqn{
\AB{\wt{m}_1,a\wt{m}_2}
~=~
\BAB{\AB{m_1,m_1'}n_1\odot\wt{m}_1,\AB{m_2,m_2'}n_2\odot\wt{m}_2}
~=~
\bfam{\AB{m_1,m_1'}n_1\odot\wt{m}_1}^*\bfam{\AB{m_2,m_2'}n_2\odot\wt{m}_2}
}\eeqn
and
\bmun{
\AB{\wt{m}_1,((m'_1\odot n_1)^*m_1)((m'_2\odot n_2)^*m_2)^*\wt{m}_2}
~=~
\BAB{((m'_1\odot n_1)^*m_1)^*\wt{m}_1,((m'_2\odot n_2)^*m_2)^*\wt{m}_2}
\\
~=~
\bAB{(m'_1\odot n_1)^*m_1,\wt{m}_1}^*\bAB{(m'_2\odot n_2)^*m_2,\wt{m}_2}
~=~
\bAB{m_1,(m'_1\odot n_1)\wt{m}_1}^*\bAB{m_2,(m'_2\odot n_2)\wt{m}_2}
\\
~=~
\AB{m_1,m'_1(n_1\odot\wt{m}_1)}^*\AB{m_2,m'_2(n_2\odot\wt{m}_2)}
~=~
\bfam{\AB{m_1,m_1'}n_1\odot\wt{m}_1}^*\bfam{\AB{m_2,m_2'}n_2\odot\wt{m}_2},
}\emun
where in the step from the second to the last line we applied Convention \ref{Massconv}.\qed

\brem
Theorem \ref{Meinvthm} is most probably folklore. But we do not know any reference. In fact, we used the statement of Theorem \ref{Meinvthm} in the proof of \cite[Corollary 1.11]{MSS06}. Since that result is too important for these notes, we decided to include a formal proof of Theorem \ref{Meinvthm} and give also a formal proof of \cite[Corollary 1.11]{MSS06}; see Corollary \ref{Morcor} below.
\erem

\brem\label{urem}
The map $u$ in the proof, actually, is what we called a \it{ternary isomorphism} in Abbaspour and Skeide \cite{AbSk07} and the Hilbert modules $N$ and $M^*$ are \it{ternary isomorphic}. (We come back to this in Section \ref{terSEC}; see, in particular, Remark \ref{utrem}.) The preceding proof is inspired very much by \cite[Theorem 2.1]{AbSk07}, which asserts that a ternary homomorphism between Hilbert modules induces a homomorphism between their range ideals. However, it would have been (notationally) more complicated to prove that $u$ is a ternary homomorphism and that the isomorphism $\cA\rightarrow\sK(M)$ induced by \cite[Theorem 2.1]{AbSk07}, is just the canonical map. Here, by calculating matrix elements, we take advantage of the fact that both $\cA$ and $\sK(M)$ are represented faithfully as operators on $M$.
\erem

\bcor\label{A1A2cor}
Let $M$ be a full Hilbert \nbd{\cB}module. Suppose we can turn $M$ into \nbd{\cA_i}\nbd{\cB}Morita equivalences $_iM$ $(i=1,2)$ via homomorphisms $\cA_i\rightarrow\sB^a(M)$ (necessarily faithful and onto $\sK(M)$ by Theorem \ref{Meinvthm}). Then there exists a unique isomorphism $\alpha\colon\cA_1\rightarrow\cA_2$ fulfilling
\beqn{
\alpha(a_1)m
~=~
a_1m
}\eeqn
for all $a_1\in\cA_1$ and $m\in M$ (first, considered an element in $_2M$ and, then, considered an element in $_1M$). Moreover, the inverses $N_i$ of $_iM$ are ternary isomorphic (via a ternary isomorphism inducing $\alpha$ as explained in Remark \ref{urem}).
\ecor

\bemp[Convention.~] \label{N=M*conv}
After this corollary, we always identify the inverse $N$ of a Morita equivalence from $\cA_1:=\cA$ to $\cB$ with the Morita equivalence $M^*$ from $\cB$ to $\cA_2:=\sK(M)$. The \nbd{\cA}valued inner product of $N=M^*$ is recovered as $\AB{m'^*,m^*}=\alpha^{-1}(m'm^*)$, where $\alpha$ is the left action of $\cA$ as in Theorem \ref{Meinvthm}. (Whether we write $N$ or $M^*$ will be decided depending on the circumstances.)
\eemp

The following representation theory of $\sK(E)$ is a simple consequence of $E\odot E^*=\sK(E)$. The uniqueness statement follows from $E^*\odot E=\cB_E$.

\bcor\label{Krepcor}
Let $E$ be a Hilbert \nbd{\cB}module. Suppose $F$ is a correspondence from $\sK(E)$ to $\cC$ (that is, a nondegenerate representation of $\sK(E)$ by adjointable operators on the Hilbert \nbd{\cC}module $F$). Then
\beqn{
F
~=~
\sK(E)\odot F
~=~
E\odot E^*\odot F
~=~
E\odot\wt{F}
}\eeqn
as \nbd{\sK(E)}\nbd{\cC}correspondences, where we defined the \nbd{\cB}\nbd{\cC}correspondence $\wt{F}:=E^*\odot F$. (The identifications are the canonical ones.)

Moreover, $\wt{F}$ is also a \nbd{\cB_E}\nbd{\cC}correspondence, and as such it is the unique (up to isomorphism) \nbd{\cB_E}\nbd{\cC}correspondence for which $F\cong E\odot\wt{F}$ (as \nbd{\sK(E)}\nbd{\cC}correspondences).
\ecor

So far, this has been discussed already by Rieffel \cite{Rie74}. Actually, Rieffel discussed representations of the pre--\nbd{C^*}algebra $\sF(E)$ on a Hilbert space. The extension to Hilbert \nbd{\cC}modules as representation spaces is marginal. The observation that the representation of $\sF(E)$ extends uniquely not only to $\sK(E)$ but even to $\sB^a(E)$ is key! In fact, if we have a nondegenerate representation of the ideal $\sF(E)$ in $\sB^a(E)$ on $F$, then $fy\mapsto(af)y$  $(f\in\sF(E),y\in F)$ for $a\in\sB^a(E)$ induces a representation of the \nbd{C^*}algebra $\sB^a(E)$ on the dense pre-Hilbert \nbd{\cC}module $\ls\sF(E)F$ of $F$. Such a representation is by bounded operators, automatically. (Note that for this it is not even necessary to require that the representation of $\sF(E)$ is by bounded operators.)

We see that unital representations $\vt$ of $\sB^a(E)$ behave well as soon as the action of $\sF(E)$ or, equivalently, of $\sK(E)$ is already nondegenerate, so that the representation module becomes a correspondences with left action of $\sK(E)$. And if $\vt$ is nonunital, we simply restrict to $\vt(\id_E)F$. We define: A homomorphism $\vt\colon\sB^a(E)\rightarrow\sB^a(F)$ is \hl{strict} if $\cls\sK(E)F=\vt(\id_E)F$.

\brem
We avoid to use the usual definition of strictness. But we remind the reader that our definition is equivalent to usual strict continuity of the representation $\vt$ on bounded subsets. (For the curious reader, and because it is, though surely folklore, not actually easy to find: By Kasparov \cite{Kas80}, $\sB^a(E)$ is the multiplier algebra of $M(\sK(E))$, and the strict topology of $M(\sK(E))$ is defined precisely as the \nbd{*}strong topology of $\sB^a(\sK(E))$; see Lemmata 1.7.10, 1.7.13, Corollary 1.7.14, and Definition 1.7.15 in Skeide \cite{Ske01}. On bounded subsets, the strict topology of $\sB^a(E)$ coincides with the \nbd{*}strong topology; \cite[Proposition 1.7.16]{Ske01} or \cite[Proposition 8.1]{Lan95}. Now, a bounded approximate unit for the finite-rank operators on $E$, clearly, converges \nbd{*}strongly to the identity of $E$. So, a representation $\vt$ that is strictly continuous on bounded subsets, has to act nondegenerately on $\vt(\id_E)F$. Conversely, the amplification map $a\mapsto a\odot\id_{\wt{F}}$ is, clearly, \nbd{*}strongly continuous on bounded subsets.)
\erem

Fixing the isomorphism suggested by the canonical identifications in Corollary \ref{Krepcor}, we obtain the representation theorem \cite[Theorem 1.4]{MSS06}.

\bthm\label{strirepthm}
Let $E$ be a Hilbert \nbd{\cB}module, let $F$ be a Hilbert \nbd{\cC}module and let $\vt\colon\sB^a(E)\rightarrow\sB^a(F)$ be a strict unital homomorphism. (In other words, $F$ is a correspondence from $\sB^a(E)$ to $\cC$ with strict left action and, thus, also a correspondence from $\sK(E)$ to $\cC$.) Then $F_\vt:=E^*\odot F$ is a correspondence from $\cB$ to $\cC$ and the formula
\beqn{
u(x_1\odot(x_2^*\odot y))
~:=~
\vt(x_1x_2^*)y
}\eeqn
defines a unitary
\beqn{
u
\colon
E\odot F_\vt
~\longrightarrow~
F
}\eeqn
such that
\beqn{
\vt(a)
~=~
u(a\odot\id_{F_\vt})u^*.
}\eeqn
\ethm

\bcor\label{Morcor}
A full Hilbert \nbd{\cB}module $E$ and a full Hilbert \nbd{\cC}module $F$ have \hl{strictly isomorphic} operator algebras (the isomorphism and its inverse are strict mappings) if and only if there is a Morita equivalence $M$ from $\cB$ to $\cC$ such that $F\cong E\odot M$.
\ecor

\proof
Suppose $F\cong E\odot M$. Then there is exactly one isomorphism $F\odot M^*\cong E$ making $(E\odot M)\odot M^*\cong F\odot M^*\cong E$ compatible with the canonical identification $E\odot(M\odot M^*)\cong E\odot\cB=E$, namely, the one defined by the chain $F\odot M^*\cong(E\odot M)\odot M^*\cong E\odot(M\odot M^*)\cong E$. In these identification, $\vt:=\bullet\odot\id_M$ defines a strict homomorphism $\sB^a(E)$ into $\sB^a(F)$ and $\theta:=\bullet\odot\id_{M^*}$ its inverse.

For the converse direction, taking into account also Theorem \ref{Meinvthm}, we reproduce the proof of \cite[Corollary 1.11]{MSS06}: The two correspondences $F_\vt$ and $E_{\vt^{-1}}$ of a bistrict isomorphism $\vt\colon\sB^a(E)\rightarrow\sB^a(F)$ have tensor products $F_\vt\odot E_{\vt^{-1}}$ and $E_{\vt^{-1}}\odot F_\vt$ that induce the identities $\vt^{-1}\circ\vt=\id_E$ and $\vt\circ\vt^{-1}=\id_F$, respectively. By the uniqueness result in Corollary \ref{Krepcor}, $F_\vt$ and $E_{\vt^{-1}}$ are inverses under tensor product. (Identities are induced by the identity correspondences.)\qed

\lf(See also the beginning of Section \ref{E0MeSEC}, where explicit identifications for the ``only if'' direction are chosen.)

\bitemp[Corollary {\cite[Remark 1.13]{MSS06}}.~]\label{bistrcor}
An isomorphism $\vt\colon\sB^a(E)\rightarrow\sB^a(F)$ is bistrict if and only if both $\vt$ and $\vt^{-1}$ take the compacts into (and, therefore, onto) the compacts.
\eitemp

\bemp[Definition {\cite[Definition 5.7]{Ske09}}.~]
A Hilbert \nbd{\cB}module $E$ and a Hilbert \nbd{\cC}module $F$ are \hl{Morita equivalent} if there is a Morita equivalence $M$ from $\cB$ to $\cC$ such that $E\odot M\cong F$ (or $E\cong F\odot M^*$).
\eemp

With this definition Corollary \ref{Morcor} may be rephrased as follows.

\bcor\label{striMocor}
Two full Hilbert modules have strictly isomorphic operator algebras if and only if they are Morita equivalent.
\ecor

\bob
If, in the notation of Corollary \ref{Morcor}, $E$ and $F$ are not necessarily full, then strict isomorphism of $\sB^a(E)$ and $\sB^a(F)$ does not necessarily imply that  $\cB$ and $\cC$ Morita equivalent. (Only $\cB_E$ and $\cC_F$ are Morita equivalent. For instance, if $\cB$ is a commutative \nbd{C^*}algebra and $\cC$ an ideal in $\cB$ not isomorphic to $\cB$, then $E=\cC$ considered as Hilbert \nbd{\cB}module and the Hilbert \nbd{\cC}module $F=\cC$ have the same compact operators. So, $\sB^a(E)$ and $\sB^a(F)$ are strictly isomorphic. But, $\cB$ and $\cC$ are not Morita equivalent, because commutative \nbd{C^*}algebras are Morita equivalent if and only if they are isomorphic.)

However, if $E$ and $F$ are Morita equivalent via $M$, say, then still $\sB^a(E)$ and $\sB^a(F)$ are strictly isomorphic and the Morita equivalence from $\cB_E$ to $\cC_F$ inducing that isomorphism is simply $\cB_E\odot M\odot\cC_F=\cls\cB_EM\cC_F$.
\eob

\brem
Anoussis and Todorov \cite{AnTo05} show that for separable \nbd{C^*}algebras and countably generated Hilbert modules every isomorphism $\sB^a(E)\rightarrow\sB^a(F)$ takes the compacts onto the compacts; see, once more, \cite[Remark 1.13]{MSS06}. Therefore, in such a situation the requirement strict is unnecessary. (For nonseparable \nbd{C^*}algebras the statement may fail. Think of the Hilbert \nbd{\cB}module $E$ and the Hilbert \nbd{\sB^a(E)}module $F:=\sB^a(E)$, which have the same operators. But $\sK(F)=\sB^a(E)\ne\sK(E)$, unless $E$ is algebraically finitely generated; see \cite[Corollary 3.4]{Ske09}.)
\erem

\section{Stable Morita equivalence for Hilbert modules}\label{stabMeSEC}

Let $\K$ denote an infinite-dimensional separable Hilbert space and denote $\sK:=\sK(\K)$. Two \nbd{C^*}algebras $\cA$ and $\cB$ are \hl{stably isomorphic} if $\cA\otimes\sK\cong\cB\otimes\sK$. A \nbd{C^*}algebra is \hl{\nbd{\sigma}unital} if it has a countable approximate unit. The main result of Brown, Green and Rieffel \cite{BGR77} asserts that two \nbd{\sigma}unital \nbd{C^*}algebras $\cA$ and $\cB$ are stably isomorphic if and only if they are Morita equivalent.

The proof of the forward direction is simple and works for arbitrary \nbd{C^*}algebras. Indeed, for a Hilbert \nbd{\cB}module $E$ we denote by $E\otimes\K$ the \hl{external tensor product} (the Hilbert \nbd{\cB}module obtained by completion from the algebraic tensor product $E\,\ul{\otimes}\,\K$ with the obvious inner product). One easily checks that $\sK(E\otimes\K)=\sK(E)\otimes\sK$. In particular, if we put $\K_\cB:=\cB\otimes\K$, then $\sK(\K_\cB)=\cB\otimes\sK$. So $\K_\cB$ is a Morita equivalence from $\cB\otimes\sK$ to $\cB$, and if $\cA\otimes\sK$ and $\cB\otimes\sK$ are isomorphic, then $\cA$ and $\cB$ are Morita equivalent.

In the version of the proof of the backward direction as presented in Lance \cite[Chapter 7]{Lan95}, the following result is key.

\bitemp[Proposition {\cite[Proposition 7.4]{Lan95}}.~]\label{Lanprop}
\begin{enumerate}
\item\label{Lp1}
Suppose $E$ is a full Hilbert module over a \nbd{\sigma}unital \nbd{C^*}algebra $\cB$. Then $E\otimes\K$ has a direct summand $\cB$.

\item\label{Lp2}
Suppose $E$ is a countably generated Hilbert \nbd{\cB}module that has $\cB$ as a direct summand. Then $E\otimes\K\cong\K_\cB$.
\end{enumerate}
So, if $E$ is a countably generated full Hilbert module over a \nbd{\sigma}unital \nbd{C^*}algebra $\cB$, then $E\otimes\K\cong\K_\cB$.
\eitemp

\brem
Part \ref{Lp1} has a much simpler proof when $\cB$ is unital. In fact, in that case $\K$ may be replaced by a suitable finite-dimensional Hilbert space; see \cite[Lemma 3.2]{Ske09}.
\erem

\brem\label{Kasrem}
The proof of Part \ref{Lp2} relies on Kasparov's \it{stabilization theorem} \cite{Kas80}. In fact, if $E=\cB\oplus F$, then $E\otimes\K\cong\K_\cB\oplus(F\otimes\K)$. Since with $E$ also $F\otimes\K$ is countably generated, the stabilization theorem asserts $\K_\cB\oplus(F\otimes\K)\cong\K_\cB$.
\erem

\bdefi\label{sMedefi}
Let $E$ and $F$ denote Hilbert modules.
\begin{enumerate}
\item
$E$ and $F$ are \hl{stably Morita equivalent} if $E\otimes\K$ and $F\otimes\K$ are Morita equivalent.

\item
$\sB^a(E)$ and $\sB^a(F)$ are \hl{stably strictly isomorphic} if $\sB^a(E\otimes\K)$ and $\sB^a(F\otimes\K)$ are strictly isomorphic.
\end{enumerate}
\edefi
By Corollary \ref{bistrcor}, ~$\sB^a(E\otimes\K)$~ and ~$\sB^a(F\otimes\K)$~ are strictly isomorphic (if and) only if ~$\sK(E\otimes\K)=\sK(E)\otimes\sK$~ and ~$\sK(F\otimes\K)=\sK(F)\otimes\sK$~ are isomorphic, that is, if and only ~$\sK(E)$~ and ~$\sK(F)$~ are stably isomorphic. (Of course, isomorphic \nbd{C^*}algebras have strictly isomorphic multiplier algebras, so $\sK(E\otimes\K)\cong\sK(F\otimes\K)$ implies $\sB^a(E\otimes\K)\cong\sB^a(F\otimes\K)$ strictly.)

\bthm\label{sMethm}
Let $E$ and $F$ denote full Hilbert modules over \nbd{C^*}algebras $\cB$ and $\cC$, respectively.

\begin{enumerate}
\item\label{sM1}
$E$ and $F$ are stably Morita equivalent if and only if $\sB^a(E)$ and $\sB^a(F)$ are stably strictly isomorphic. Either condition implies that $\cB$ and $\cC$ are Morita equivalent.

\item\label{sM2}
Suppose $E$ and $F$ are countably generated and $\cB$ and $\cC$ are \nbd{\sigma}uni\-tal. Then the following conditions are all equivalent:
\begin{enumerate}
\item\label{sM2i}
$E$ and $F$ are stably Morita equivalent.

\item\label{sM2ii}
$\sB^a(E)$ and $\sB^a(F)$ are stably strictly isomorphic.

\item\label{sM2iii}
$\cB$ and $\cC$ are Morita equivalent.

\item\label{sM2iv}
$\cB$ and $\cC$ are stably isomorphic.

\end{enumerate}
\end{enumerate}
\ethm

\proof
Part \ref{sM1} is Corollary \ref{striMocor} and equivalence of \eqref{sM2i} and \eqref{sM2ii} is \ref{sM1} restricted to the special case. Equivalence of \eqref{sM2iii} and \eqref{sM2iv} is \cite{BGR77}. Clearly, \eqref{sM2i} $\Longrightarrow$ \eqref{sM2iii} directly from the definition, while \eqref{sM2iv} $\Longrightarrow$ \eqref{sM2i} follows from Proposition \ref{Lanprop} and the observation that $\K_\cB$ and $\K_\cC$ are Morita equivalent if $\cB$ and $\cC$ are, via the same Morita equivalence.\qed

\newpage

\section{Ternary isomorphisms}\label{terSEC}

The isomorphisms in the category of Hilbert \nbd{\cB}modules are the \hl{unitaries}, that is, the inner product preserving surjections. If $u\colon E\rightarrow F$ is a unitary, then the \hl{conjugation} $u\bullet u^*\colon\sB^a(E)\rightarrow\sB^a(F)$ defines a strict isomorphism. If $E$ and $F$ are isomorphic, we say $\sB^a(E)$ and $\sB^a(F)$ are \hl{inner conjugate}.

In the sequel, we shall say that strictly isomorphic $\sB^a(E)$ and $\sB^a(F)$ are \hl{strictly conjugate}. We know that a full Hilbert \nbd{\cB}module $E$ and a full Hilbert \nbd{\cC}module $F$ have strictly isomorphic operator algebras if and only if the modules are Morita equivalent. Isomorphic Hilbert \nbd{\cB}modules are Morita equivalent via the identity correspondence $\cB$. But, Morita equivalent full Hilbert \nbd{\cB}modules need not be isomorphic.

\bex\label{nuvMex}
Let $\cB:=\sMatrix{\C&0\\0&M_2}\subset M_3$ an let $E:=\sMatrix{0&{\C^2}^*\\\C^2&0}\subset M_3$ be the \nbd{\cB}correspondence obtained by restricting the operations of the identity \nbd{M_3}correspondence $M_3$ to the subsets $E$ and $\cB$. Then $E$ is a self-inverse Morita equivalence. From $E\odot E=\cB$ it follows that $E$ and $\cB$ are Morita equivalent as Hilbert \nbd{\cB}modules. Of course, they are not isomorphic. In fact, their dimensions as complex vector spaces differ, so that there is not even a linear bijection between them. Also, $E$ has no unit vector, but $\cB$, of course, has one.
\eex

In between isomorphism of Hilbert \nbd{\cB}modules and Morita equivalence of Hilbert modules there is another equivalence relation, based on ternary isomorphisms. A \hl{ternary homomorphism} from a Hilbert \nbd{\cB}module $E$ to a Hilbert \nbd{\cC}module $F$ is a linear map $u\colon E\rightarrow F$ (\it{a priori} not bounded; see, however, Skeide and Sumesh \cite[Footnote 2]{SkSu14} for why linearity cannot be dropped) that satisfies
\beqn{
u(x\AB{y,z})
~=~
(ux)\AB{uy,uz}
}\eeqn
for all $x,y,z\in E$. A \hl{ternary unitary} is a bijective ternary homomorphism. Clearly, if $u$ is a ternary unitary, then so is $u^{-1}$. If there is a ternary unitary from $E$ to $F$, we say $E$ and $F$ are \hl{ternary isomorphic}.

Ternary homomorphisms have the advantage that they do not refer in any way to the \nbd{C^*}al\-ge\-bras over which the modules are modules. (In fact, we may turn the class of all Hilbert modules, without fixing an algebra, into a category by choosing as morphisms the ternary homomorphisms.) The following notion takes into account the algebras more explicitly. A \hl{generalized isometry} from a Hilbert \nbd{\cB}module $E$ to a Hilbert \nbd{\cC}module $F$ is a map $u\colon E\rightarrow F$ (\it{a priori} neither linear nor bounded) such that there exists a homomorphism $\vp\colon\cB\rightarrow\cC$ fulfilling
\beq{\label{geniso}
\AB{ux,uy}
~=~
\vp(\AB{x,y})
}\eeq
for all $x,y\in E$. Once the homomorphism $\vp$ is fixed, we shall also speak of a \hl{\nbd{\vp}isometry} $u$; see Skeide \cite{Ske06c}.

The connection between ternary homomorphisms and generalized isometries is made by the following result.

\bitemp[Theorem {\cite[Theorem 2.1]{AbSk07}}.~]\label{gen=ter}
For a map $u$ from a full Hilbert \nbd{\cB}module $E$ to a Hilbert \nbd{\cC}module $F$ the following statements are equivalent:
\begin{enumerate}
\item
$u$ is a generalized isometry.

\item
$u$ is a ternary homomorphism.
\end{enumerate}
\eitemp

\brem\label{utrem}
Of course, the homomorphism $\vp$ turning a ternary homomorphism into a generalized isometry is the unique homomorphism satisfying \eqref{geniso}. This is essentially what we used in the proof of Theorem \ref{Meinvthm}. As mentioned in Remark \ref{urem}, the map $u$ in that proof is a ternary homomorphism. Just that it was easier in the particular case to establish $u$ as a generalized isometry. \cite[Theorem 2.1]{AbSk07} now assures that $u$ is, indeed, a ternary homomorphism and, therefore, a ternary unitary.

Recognizing a ternary homomorphism as a \nbd{\vp}isometry has more consequences. For instance, there is the notion of \nbd{\vp}adjointable operators with all results but also with all problems known from the usual adjointable operators; see \cite[Observation 1.9]{Ske06c}.
\erem

\bob
Clearly, a generalized isometry is linear and contractive (even completely contractive; see Theorem \ref{terhextthm}). Therefore, so is a ternary homomorphism.
\eob

Note that if $u\colon E\rightarrow F$ is a ternary homomorphism, then $u^*\colon x^*\mapsto(ux)^*$ is a ternary homomorphism from $E^*$ to $F^*$. The following theorem now follows easily from \cite[Theorem 2.1]{AbSk07}. We omit the proof.

\bthm\label{terhextthm}
For a map $u$ from a Hilbert \nbd{\cB}module $E$ to a Hilbert \nbd{\cC}module $F$ the following statements are equivalent:
\begin{enumerate}
\item
$u$ is a ternary homomorphism.

\item
$u$ extends as a (unique!) homomorphism between the \hl{reduced linking algebras}
\beqn{
\Phi_u
\colon
\sMatrix{\cB_E&E^*\\E&\sK(E)}
~\longrightarrow~
\sMatrix{\cC_F&F^*\\F&\sK(F)}
}\eeqn
respecting the corners.
\end{enumerate}
In either case, $\Phi_u$ is injective if and only if $u$ is, and $\Phi_u$ is surjective if and only if $u$ is.
\ethm

Let $u\colon E\rightarrow F$ be a ternary unitary. Then the conjugation map $u\bullet u^{-1}\colon a\mapsto uau^{-1}$ is, clearly, multiplicative. Note also
\bmun{
(ux)\AB{uy,(uau^{-1})uz}
~=~
(ux)\AB{uy,uaz}
~=~
u(x\AB{y,az})
\\
~=~
u(x\AB{a^*y,z})
~=~
(ux)\AB{ua^*y,uz}
~=~
(ux)\AB{(ua^*u^{-1})uy,uz}
}\emun
for all $x,y,z\in E$, so that $u\bullet u^{-1}$ is a \nbd{*}map. In other words, conjugation with $u$ still defines an isomorphism from $\sB^a(E)$ to $\sB^a(F)$. We call an isomorphism obtained by conjugation with a ternary unitary a \hl{ternary inner isomorphism} and we call $\sB^a(E)$ and $\sB^a(F)$ \hl{ternary conjugate}.

The restriction of a ternary inner isomorphisms induced by $u$ to $\sK(E)$ is precisely the restriction of $\Phi_u$ to $\sK(E)$ and, therefore, an isomorphism onto $\sK(F)$. It follows that ternary inner isomorphisms are bistrict. The Morita equivalence from $\cB_E$ to $\cC_F$ of such a ternary inner isomorphism is simply $_\vp\cC_F$, where $\vp$ is the restriction of $\Phi_u$ to an isomorphism from $\cB_E$ to $\cC_F$.

\brem
In Skeide \cite{Ske06c} we have analyzed the group of ternary inner automorphisms of $\sB^a(E)$ and how it is reflected in the Picard group of $\cB$. The \hl{Picard group} of $\cB$ is the group of isomorphism classes of Morita equivalences over $\cB$ under tensor product; see \cite{BGR77}. It contains the (opposite of the) group of the automorphisms of $\cB$ modulo the multiplier inner automorphisms. (Multiplier inner automorphisms are obtained by conjugation with a unitary in the multiplier algebra. They are called generalized inner automorphisms in \cite{BGR77,Ske06c}. We now prefer to follow the modern terminology in Blackadar \cite{Bla06}.) One main point of \cite{Ske06c} is, very roughly, that there are full Hilbert \nbd{\cB}modules $E$ such that not every automorphism of $\cB$ occurs as the automorphism $\vp$ induced by a ternary unitary on $u$ on $E$. In other words, not all automorphisms of $\cB$ extend to automorphisms of the linking algebra of $E$. Equivalently, not for every automorphism $\vp$ of $\cB$ the Hilbert \nbd{\cB}modules $E$ and $E\odot{_\vp\cB}$ are isomorphic.
\erem

\bdefi\label{terisodef}
Let $E$ and $F$ denote Hilbert modules.
\begin{enumerate}
\item
$E$ and $F$ are \hl{stably ternary isomorphic} if $E\otimes\K$ and $F\otimes\K$ are ternary isomorphic.

\item
$\sB^a(E)$ ~and~ $\sB^a(F)$ are \hl{stably ternary conjugate} if ~$\sB^a(E\otimes\K)$~ and ~$\sB^a(F\otimes\K)$~ are ternary conjugate.
\end{enumerate}
By definition of ternary conjugate, the two properties are equivalent.
\edefi

\bthm\label{ter=isothm}
Let $E$ and $F$ denote full countably generated Hilbert modules over \nbd{\sigma}unital \nbd{C^*}algebras $\cB$ and $\cC$, respectively. Then either of the conditions in Definition \ref{terisodef} holds if and only if $\cB$ and $\cC$ are isomorphic.
\ethm

\proof
$E\otimes\K$ and $F\otimes\K$ are full Hilbert modules over $\cB$ and $\cC$, respectively. Suppose the first condition of Definition \ref{terisodef} holds. Then $E\otimes\K$ and $F\otimes\K$ are ternary isomorphic, so that $\cB$ and $\cC$ are isomorphic. (This does not depend on countability hypotheses.) Suppose, on the other hand, $\cB$ and $\cC$ are isomorphic via an isomorphism $\vp$, say. Then $E_\vp:=E\odot{_\vp}\cC$ is a full countably generated Hilbert \nbd{\cC}module ternary isomorphic to $E$ via $x\odot c\mapsto x\vp^{-1}(c)$. By Proposition \ref{Lanprop} the Hilbert \nbd{\cC}modules $E_\vp\otimes\K$ and $F\otimes\K$ are isomorphic, so that $E\otimes\K$ and $F\otimes\K$ are ternary isomorphic.\qed

\section{Cocycle conjugacy of $E_0$--semigroups}\label{ccSEC}

In this section we discuss several notions of cocycle conjugacy in an algebraic context. (We use \nbd{C^*}alge\-bras just for convenience. General unital \nbd{*}algebras, like in the beginning of the methodological introduction, would do as well.) We put particular emphasis on the fact that, unlike the case $\sB(H)$ where all cocycles are unitarily implemented, here the character of the automorphisms forming the cocycles may vary. In the end, also in these notes we shall concentrate on unitary cocycles. (See, however, Remark \ref{nunicocrem}.) But, this is our choice, and it is \bf{this} choice that does the job of resolving our classification problem. The differences in  the several notions of \it{unitary cocycle conjugacy} we employ, lie in what sort of conjugacies we will allow.

Let $\vt$ and $\theta$ denote unital endomorphisms of unital \nbd{C^*}algebras $\cA$ and $\cB$, respectively. (For nonunital algebras one would replace unital with nondegenerate in the sense that $\vt(\cA)\cA$ should be total in $\cA$. We do not tackle these problems. Though, interesting phenomena may happen, worth of a separate investigation.) $\vt$ and $\theta$ are \hl{conjugate} if there exists an isomorphism $\alpha\colon\cA\rightarrow\cB$ such that $\alpha\circ\vt=\theta\circ\alpha$. If $\vt$ and $\theta$ are conjugate, then for every $n\in\N_0$ the members $\vt_n:=\vt^n$ and $\theta_n:=\theta^n$ of the \nbd{E_0}semigroups generated by $\vt$ and $\theta$, respectively, are conjugate via the same isomorphism $\alpha$. In general, we say two \nbd{E_0}semigroups $\vt$ and $\theta$ are \hl{conjugate} if there is an isomorphism $\alpha$ such that $\alpha\circ\vt_t=\theta_t\circ\alpha$ for all $t\in\bS$. Of course, conjugacy of \nbd{E_0}semigroups is an equivalence relation.

If $\vt$ and $\theta$ are two unital endomorphisms of the \bf{same} unital \nbd{C^*}algebra $\cA$, then we may ask whether there is an automorphism $\beta$ of $\cA$, such that $\theta=\beta\circ\vt$. In this case, we may not expect that $\theta_n=(\beta\circ\vt)^n$ would be equal to $\beta\circ\vt_n$. In the (rare) case when $\beta$ and $\vt$ commute, we find $\theta_n=\beta^n\circ\vt_n$. In general, we may not even expect that there exist automorphisms $\beta_n$ such that $\theta_n=\beta_n\circ\vt_n$.

\bex
Let $\cA=\C^3$, the diagonal subalgebra of $M_3$. Define the one-sided shift $\vt\tMatrix{a\\b\\c}:=\tMatrix{a\\a\\b}$ and the cyclic permutation $\beta\tMatrix{a\\b\\c}:=\tMatrix{c\\a\\b}$. Put $\theta:=\beta\circ\vt$, so that $\theta\tMatrix{a\\b\\c}=\tMatrix{b\\a\\a}$. Then $\vt^2\tMatrix{a\\b\\c}=\tMatrix{a\\a\\a}$, while $\theta^2\tMatrix{a\\b\\c}=\tMatrix{a\\b\\b}$. There is no automorphism $\beta_2$ such that $\beta_2\tMatrix{a\\a\\a}=\tMatrix{a\\b\\b}$ for all $a,b\in\C$.
\eex

\bdefi
We say two \nbd{E_0}semigroups $\vt$ and $\theta$ on $\cA$ are \hl{cocycle equivalent} if there exist automorphisms $\beta_t$ of $\cA$ such that $\theta_t=\beta_t\circ\vt_t$. If, \it{vice versa}, $\vt$ is an \nbd{E_0}semigroup and $\beta=\bfam{\beta_t}_{t\in\bS}$ is a family of automorphisms such that $\vt^\beta_t:=\beta_t\circ\vt_t$ defines an \nbd{E_0}semigroup $\vt^\beta$, then we say $\beta$ is a \hl{cocycle on $\cA$} with respect to $\vt$.
\edefi

Clearly, if $\beta$ is a cocycle with respect to $\vt$, then $\beta^{-1}=\bfam{\beta_t^{-1}}_{t\in\bS}$ is a cocycle with respect to $\vt^\beta$. So, cocycle equivalence is symmetric. Clearly, it is reflexive and transitive. In other words, cocycle equivalence is an equivalence relation.

The reader might ask, why we used the name \it{cocycle equivalent} instead of the more common \it{cocycle conjugate}. The reason is that in a minute we will define the second term in a different way, which is closer to what is known as \it{cocycle conjugate}.

\lf
Cocycle equivalence is a notion that involves two semigroups of endomorphisms on the \bf{same} algebra. A relation that allows to compare (semigroups of) endomorphisms on \bf{different} algebras is conjugacy. Before we can investigate two semigroups of endomorphisms on different algebras for cocycle equivalence, we must transport one of them to the other algebra via a conjugacy.

\bdefi
Let $\vt$ and $\theta$ denote \nbd{E_0}semigroups on unital \nbd{C^*}algebras $\cA$ and $\cB$, respectively. We say $\vt$ and $\theta$ are \hl{cocycle conjugate} if there exists an isomorphism $\alpha\colon\cA\rightarrow\cB$ such that the \hl{conjugate} \nbd{E_0}semigroup $\vt^\alpha:=\bfam{\alpha\circ\vt_t\circ\alpha^{-1}}_{t\in\bS}$ on $\cB$ and $\theta$ are cocycle equivalent.

If $\alpha$ satisfies additional conditions, then we will indicate these in front of the word \it{conjugate}. (For instance, if $\alpha$ is an inner isomorphism, we will say $\vt$ and $\theta$ are \hl{cocycle inner conjugate}.) If the cocycle satisfies additional conditions, then we will indicate these in front of the word \it{cocycle}. (For instance, if $\beta$ consists of inner automorphisms we will say $\vt$ and $\theta$ are \hl{inner cocycle conjugate}.)
\edefi

Also for two \nbd{E_0}semigroups $\vt$ and $\theta$ on the same unital \nbd{C^*}algebra $\cA$ we may ask, whether they are cocycle conjugate. Of course, cocycle equivalent \nbd{E_0}semigroups are cocycle (inner) conjugate via $\alpha=\id_\cA$. But the converse need not be true.

\bex\label{C2ex}
Let $\cA=\C^2$, the diagonal subalgebra of $M_2$. Define the one-sided shift $\vt\raisebox{.3ex}{\tMatrix{a\\b}}:=\raisebox{.3ex}{\tMatrix{a\\a}}$ and the flip automorphism $\alpha\raisebox{.3ex}{\tMatrix{a\\b}}:=\raisebox{.3ex}{\tMatrix{b\\a}}$. Then $\vt$ and $\theta:=\vt^\alpha=\alpha\circ\vt\circ\alpha^{-1}\colon\raisebox{.3ex}{\tMatrix{a\\b}}\mapsto\raisebox{.3ex}{\tMatrix{b\\b}}$ and, therefore, the whole semigroups $\vt^n$ and $\theta^n$ generated by them are conjugate. \it{A fortiori} these two semigroups are cocycle conjugate via $\alpha$ by the identity cocycle $\beta_n=\id_\cA$. But, no automorphism $\beta$ can recover $\theta$ as $\beta\circ\vt$. So, these semigroups are not cocycle equivalent.
\eex

\bprop\label{Acocprop}
Two \nbd{E_0}semigroups $\vt$ and $\theta$ on $\cA$ are cocycle equivalent if and only if they are cocycle inner conjugate.
\eprop

\proof
The forward implication being clear, suppose $\alpha=u\bullet u^*$ is an inner automorphism (for some unitary $u\in\cA$) and $\beta=\bfam{\beta_t}_{t\in\bS}$ a family of automorphisms $\beta_t$ of $\cA$ such that $\beta_t\circ\vt^\alpha_t=\theta_t$. That is,
\beqn{
\theta_t(a)
~=~
\beta_t(u\vt_t(u^*au)u^*)
~=~
\beta_t\bfam{(u\vt_t(u^*))\vt_t(a)(u\vt_t(u^*))^*}.
}\eeqn
In other words, $\theta_t=\beta'_t\circ\vt_t$ for the automorphism $\beta'_t:=\beta_t\circ\bfam{(u\vt_t(u^*))\bullet(u\vt_t(u^*))^*}$.\qed

\lf
Since $\cA=\sB(H)$ ($H$ some Hilbert space) has only inner automorphisms, the notions of cocycle equivalence and cocycle conjugacy for \nbd{E_0}semigroups $\sB(H)$ coincide. But for $\cA=\sB^a(E)$, of course, this is not so. In fact, Example \ref{C2ex} gives a counterexample via $E:=\cA=\sB^a(E)$.

Among the inner cocycles $\beta$ with respect to $\vt$, a particularly important class consists of those cocycles that are generated as $\beta_t=u_t\bullet u_t^*$ where $u=\bfam{u_t}_{t\in\bS}$ is a family of unitaries in $\cA$ fulfilling
\beqn{
u_0
~=~
\U,
\text{~~~~~~and~~~~~~}
u_{s+t}
~=~
u_s\vt_s(u_t)
}\eeqn
for all $s,t\in\bS$. Such a family is called a \hl{unitary left cocycle in $\cA$} with respect to $\vt$ (or simply a \hl{left cocycle} if the $u_t$ are not necessarily unitary). It is easy to check that every unitary left cocycle implements a cocycle $\beta^u$ via $\beta^u_t:=u_t\bullet u_t^*$. We will say $\beta^u$ is a \hl{unitary cocycle on $\cA$}, and we will denote $\vt^u:=\vt^{\beta^u}$.

\bdefi \label{unicocdef}
Two \nbd{E_0}semigroups are \hl{unitary cocycle conjugate} (\hl{equivalent}) if the conjugacy (the equivalence) can be implemented by a unitary left cocycle.
\edefi

\bex
Suppose two \nbd{E_0}semigroups $\vt$ and $\theta$ on $\cA$ are inner conjugate via a unitary $u\in\cA$. It is easy to check that $u_t:=u\vt_t(u^*)$ is a unitary left cocycle with respect to $\vt$ and that $\theta=\vt^u$. In other words, inner conjugate \nbd{E_0}semigroups on $\cA$ are unitary cocycle equivalent. More generally, if the cocycle $\beta=\beta^v$ in the proof of Proposition \ref{Acocprop} is implemented by a unitary left cocycle $v$ with respect to $\vt^u$, then $\beta'=\beta'^{v'}$ where $v'$ is the unitary left cocycle with respect to $\vt$ defined by $v'_t:=v_tu\vt_t(u^*)$. Indeed, from $v_{s+t}=v_su\vt_s(u^*v_tu)u^*$ on easily verifies that
\beqn{
v'_{s+t}
~=~
v_{s+t}u\vt_{s+t}(u^*)
~=~
v_su\vt_s(u^*v_tu)u^*u\vt_{s+t}(u^*)
~=~
v_su\vt_s(u^*)\vt_s(v_tu\vt_t(u^*))
~=~
v'_s\vt_s(v'_t).
}\eeqn
\eex

\bcor\label{Aucoccor}
Two \nbd{E_0}semigroups $\vt$ and $\theta$ on $\cA$ are unitary cocycle equivalent if and only if they are unitary cocycle inner conjugate.
\ecor

It is easy to check that also unitary cocycle conjugacy or equivalence are equivalence relations. It is unitary cocycle conjugacy that, usually, corresponds to cocycle conjugacy in the literature; see, for instance, Takesaki \cite[Definition X.1.5]{Tak03a}. But be aware that this notion refers rather to the context of groups, not so much to semigroups.

We do not tackle the questions whether every inner cocycle is implemented by a unitary left cocycle, or to what extent the cocycle $\beta$ is non-unique. (See, however, again Remark \ref{nunicocrem}.) We just mention the following easy to prove fact. (Recall that a left cocycle $u$ with respect to $\vt$ is \hl{local} if $u_t\vt_t(a)=\vt_t(a)u_t$ for all $t\in\bS,a\in\cA$. Every local left cocycle is also a right cocycle. Therefore, we usually say just local cocycle.)

\bprop
Two unitary left cocycles $u$ and $v$ implement the same inner cocycle $\beta$ if and only if the elements $v_t^*u_t$ form a local cocycle.
\eprop

\section{$E_0$--Semigroups, product systems, and unitary cocycles}\label{E0cocSEC}

In this section we, finally, explain how the representation theory of $\sB^a(E)$ for a Hilbert \nbd{\cB}mod\-ule $E$ gives rise to the construction of a product system of \nbd{\cB}correspondences from an \nbd{E_0}semi\-group on $\sB^a(E)$. We show that the \nbd{E_0}semi\-groups acting on a fixed $\sB^a(E)$ are classified by their product systems up to unitary cocycle equivalence. We also give a criterion when \nbd{E_0}semigroups acting on possibly different $\sB^a(E)$, varying the (full) Hilbert module $E$ but over a fixed \nbd{C^*}algebra have the same product system. This criterion does, however, not involve cocycle conjugacy and is, therefore, not the criterion we really want.

Let $E$, $F$, and $G$ denote a Hilbert \nbd{\cB}, a \nbd{\cC}, and a \nbd{\cD}module, respectively. Suppose $\vt\colon\sB^a(E)\rightarrow\sB^a(F)$ and $\theta\colon\sB^a(F)\rightarrow\sB^a(G)$ are unital strict homomorphisms. Then the multiplicity correspondences $F_\vt$ and $G_\theta$ (see Theorem \ref{strirepthm}) compose contravariantly as tensor product by the isomorphism
\beq{\label{homcomp}
F_\vt\odot G_\theta
~~~\ni~~~
(x^*\odot_\vt y)\odot(y'^*\odot_\theta z)
~~~\longmapsto~~~
x\odot_{\theta\circ\vt}\theta(yy'^*)z
~~~\in~~~
G_{\theta\circ\vt}.
}\eeq
Moreover, under iterations these isomorphisms compose associatively; see \cite[Theorem 1.14]{MSS06}.

It follows that every equality between compositions of unital strict homomorphisms is reflected by an isomorphism of the corresponding tensor products of the multiplicity correspondences in the reverse order. If $E$ is a Hilbert \nbd{\cB}module and if $\vt$ is a strict \nbd{E_0}semigroup on $\sB^a(E)$, then the semigroup property $\vt_t\circ\vt_s=\vt_{s+t}$ gives rise to isomorphisms
\beqn{
u_{s,t}
\colon
E_s\odot E_t
~\longrightarrow~
E_{s+t}
}\eeqn
of the multiplicity \nbd{\cB}correspondences $E_t:=E_{\vt_t}$, $t>0$. (We shall abbreviate $x^*\odot_{\vt_t}y=:x^*\odot_t y$.) The ``multiplication''
\beqn{
\bfam{(x^*\odot_s x')\,,\,(y^*\odot_t y')}
~\longmapsto~
(x^*\odot_s x')(y^*\odot_t y')
~:=~
u_{s,t}\bfam{(x^*\odot_s x')\odot(y^*\odot_t y')}
~=~
x^*\odot_{s+t}\vt_t(x'y^*)y'
}\eeqn
is associative. If $E$ is full, then everything also extends to $t=0$ with $E_0=E^*\odot E=\cB$, and $u_{0,t}$ and $u_{t,0}$ are just the canonical left and right action, respectively, of $E_0=\cB$. If $E$ is not full, then we put $E_0:=\cB$ by hand and the canonical actions $u_{0,t}$ and $u_{t,0}$ extend uniquely the above identifications from $\cB_E=E^*\odot E=E^*\odot_0 E$ to $\cB=E_0$.

A family $E^\odot=\bfam{E_t}_{t\in\bS}$ of correspondences over $\cB$ with associative identifications $u_{s,t}$ and the conditions on $E_0,u_{0,t},u_{0,t}$ at zero is what has been called a \hl{product system} in Bhat and Skeide \cite{BhSk00}. We call $E^\odot$ constructed as above from $\vt$ the product system \hl{associated} with the strict \nbd{E_0}semigroup $\vt$.

\bemp[Convention.~]\label{fconv}
In the sequel, we restrict our attention to full Hilbert modules.
\eemp

For full Hilbert modules the multiplicity correspondence of a unital strict homomorphism is unique (up to isomorphism) and the condition $E_0=\cB$ is automatic. As far as we are dealing with the connection between \nbd{E_0}semigroups on $\sB^a(E)$ and product systems associated with them, it is natural to restrict to full Hilbert \nbd{\cB}modules, as $\cB$ can always be replaced with $\cB_E$. When we take also into account continuity questions, then $E_0:=\cB=\cB_E$ is forced. (Observe that $\cB_{E_t}=\cB_E$. Therefore, $E^\odot$ will never have continuous sections reaching every point of $E_0=\cB$, unless $\cB=\cB_E$.)

Note that the product system associated with a strict \nbd{E_0}semigroup on $\sB^a(E)$ for a full Hilbert \nbd{\cB}module $E$ must be \hl{full} in the sense that $E_t$ is full for each $t\in\bS$. Note, too, that the \nbd{E_0}semigroup consists of faithful endomorphisms if and only if the associated product system is \hl{faithful} in the sense that all $E_t$ have a faithful left action.

Remember (see again Theorem \ref{strirepthm}) that the multiplicity correspondence $E_t$ of $\vt_t$ is related with $\vt_t$ via a unitary $v_t\colon E\odot E_t\rightarrow E$ such that $\vt_t(a)=v_t(a\odot\id_t)v_t^*$. Moreover, the ``multiplication'' $(x,y_t)\mapsto xy_t:=v_t(x\odot y_t)$ iterates associatively with the product system multiplication, that is, $(xy_s)z_t=x(y_sz_t)$. A family $v=\bfam{v_t}_{t\in\bS}$ fulfilling these properties is what we started calling a \hl{left dilation} of the full product system $E^\odot$ to the full Hilbert module $E$ in \cite{Ske06} (for the Hilbert space case) and in \cite{Ske07,Ske11a}. (For nonfull $E^\odot$ the term \it{left dilation} is not defined.%
\footnote{If $E$ is not necessarily full, then we speak of a \hl{left quasi dilation}. This is an interesting concept, too. But it has no nice relation with \nbd{E_0}semigroups. ($E$ may be very well $\zero$.) In these notes we are interested only in the relation between \nbd{E_0}semigroups and product systems. There is also a relation of product systems with \nbd{E}semigroups, that is, semigroups of not necessarily unital endomorphisms. In that case the $v_t$ need not be unitary but just isometric, and we speak of \hl{semidilations}. In each case, it is forced that $v_0$ is the canonical identification $x\odot b\mapsto xb$.

There is also the concept of \hl{right dilation} \cite{Ske11a} of faithful product systems, which is practically synonymous with \it{faithful nondegenerate representation} of a product system; see \cite{Ske09}. Also this concept, for Hilbert modules, is not directly related to \nbd{E_0}semigroups, while, for von Neumann modules, it parallels Arveson's approach to Arveson systems. We come back to right dilations in Appendix B.2.}%
) If $E^\odot$ is the product system associated with a strict \nbd{E_0}semigroup and the left dilation arises in the prescribed way, then we refer to it as the \hl{standard dilation of $E^\odot$}.

For every left dilation of a product system $E^\odot$ to $E$, by $\vt^v_t(a):=v_t(a\odot\id_t)v_t^*$ we define a strict \nbd{E_0}semigroup $\vt^v$ on $\sB^a(E)$. We say a strict \nbd{E_0}semigroup $\vt$ is \hl{associated} with a full product system $E^\odot$ if it can be obtained as $\vt=\vt^v$ for some left dilation $v$ of $E^\odot$. Of course, every strict \nbd{E_0}semigroup is associated with its associated product system via the standard dilation. But, this need not be the only left dilation that gives back the \nbd{E_0}semigroup. We now investigate the possibilities in the slightly more general situation when two left dilations of two product systems induce conjugate \nbd{E_0}semigroups.

\bdefi
A \hl{morphism} between two product systems $E^\odot$ and $F^\odot$ of \nbd{\cB}correspondences is a family $w^\odot=\bfam{w_t}_{t\in\bS}$ of bilinear adjointable maps $w_t\colon E_t\rightarrow F_t$ such that $w_s(x_s)w_t(y_t)=w_{s+t}(x_sy_t)$ and $w_0=\id_\cB$. An \hl{isomorphism} is a morphism that consists of unitaries. Of course, the inverse of an isomorphism is an isomorphism.
\edefi

The following proposition is very important. It summarizes answers to all questions about product systems associated via left dilations with \nbd{E_0}semigroups, when these \nbd{E_0}semigroups are inner conjugate.

\bprop\label{conisoprop}
Let $E$ and $E'$ be full Hilbert \nbd{\cB}modules. Let $u\in\sB^a(E,E')$ be a unitary and define the inner isomorphism $\alpha:=u\bullet u^*$. Suppose $v,v'$ are left dilations of product systems $E^\odot,E'^\odot$ to $E$ and $E'$, respectively, such that $(\vt^v)^\alpha=\vt^{v'}$. Then there is a unique isomorphism $w^\odot$ from $E^\odot$ to $E'^\odot$ such that\vspace{-.7ex}
\beq{\label{isouni}
u(xy_t)
~=~
(ux)(w_ty_t)
}\eeq
for all $t\in\bS,x\in E,y_t\in E_t$. (That means $uv_t=v'_t(u\odot w_t)$.) In particular:
\begin{enumerate}
\item\label{CIP1}
If $u$ is the identity of $E=E'$, so that $\vt^v=\vt^{v'}$, then $w^\odot$ is the unique isomorphism satisfying
\beqn{
xy_t
~=~
x(w_ty_t).
}\eeqn

\item\label{CIP2}
If, in the situation of \ref{CIP1}, $E'^\odot$ is the product system associated with $\vt^v$ and $v'$ its standard dilation, then $w^\odot$ is the unique isomorphism satisfying
\beqn{
w_t(\AB{x,y}z_t)
~=~
x^*\odot_t(yz_t).
}\eeqn

\item\label{CIP3}
Let $\vt$ be a strict \nbd{E_0}semigroup on $\sB^a(E)$. If $E^\odot$ and $E'^\odot$ are the product systems associated with $\vt$ and $\vt':=\vt^\alpha$, respectively, and if $v$ and $v'$ are their respective standard dilations, then $w^\odot$ is the unique isomorphism determined by
\beqn{
w_t(x^*\odot_ty)
~=~
(ux)^*\odot_t(uy).
}\eeqn
\end{enumerate}
\eprop

\proof
For uniqueness, suppose that $w_t$ and $w'_t$ are maps satisfying $(ux)(w_ty_t)=(ux)(w'_ty_t)$. Then $(ux)(w_ty_t-w'_ty_t)=v'_t(u\odot\id_t)(x\odot(w_ty_t-w'_ty_t))=0$ for all $x\in E$. Since $E$ is full and $\cB$ acts nondegenerately, one easily verifies that, in a tensor product, $x\odot y=0$ for all $x\in E$ implies $y=0$. In other words, $w_ty_t-w'_ty_t=0$ for all $y_t\in E_t$ or $w_t=w'_t$.

By associativity of left dilations, maps $w_t$ fulfilling \eqref{isouni} for all $t\in\bS$ satisfy
\beqn{
(ux)w_{s+t}(y_sz_t)
~=~
u(xy_sz_t)
~=~
(u(xy_s))(w_tz_t)
~=~
(ux)(w_sy_s)(w_tz_t)
}\eeqn
for all $x\in E$. Once more, by uniqueness we get $w_{s+t}(y_sz_t)=(w_sy_s)(w_tz_t)$. Putting $s=0$, they are left-linear. Further, from \eqref{isouni} we conclude $\AB{y_t,\AB{x,x'}y'_t}=\AB{w_ty_t,\AB{x,x'}w_ty'_t}$. Once more, by fullness of $E$ and by nondegeneracy of $\cB$, the maps $w_t$ must be isometries.

It follows that maps $w_t$ satisfying \eqref{isouni}, once they exist, are uniquely determined and form an isometric morphism. It remains to establish mappings $w_t$ with total range that satisfy \eqref{isouni}.

Note that \eqref{isouni} determines the $w_t$ only indirectly, while the properties stated in \ref{CIP2} and \ref{CIP3} in either special case can be used as a direct definition. We just observe that it is enough to prove the cases \ref{CIP2} and \ref{CIP3}, separately to prove also the general statement. (The general situation can be decomposed into an isomorphism of $E^\odot$ and the product system associated with $\vt^v$, an isomorphism between the product systems associated with $\vt$ and with $\vt^\alpha$, and an isomorphism between the product system associated with $\vt^{v'}$ and $E'^\odot$.) \ref{CIP1}, instead, is simply the restriction of the general statement to the special case. So, it only remains to prove \ref{CIP2} and \ref{CIP3}.

To prove \ref{CIP2}, we mention the identities $\AB{x^*\odot_ty,x'^*\odot_ty'}=\AB{y,\vt^v_t(xx'^*)y'}$ and $\vt^v_t(a)(xy_t)=(ax)y_t$, which follow directly from the definitions. Observe that the stated $w_t$ are isometric (and, therefore, well-defined) and surjective. Since
\beqn{
z(w_t(\AB{x,x'}y_t))
~=~
z(x^*\odot_t(x'y_t))
~=~
\vt^v_t(zx^*)(x'y_t)
~=~
(z\AB{x,x'})y_t
~=~
z(\AB{x,x'}y_t),
}\eeqn
the $w_t$ fulfill \eqref{isouni}.

Similarly, to prove \ref{CIP3}, we mention, in addition to the preceding relations, that $\vt^{v'}_t((ux)(uy)^*)=\vt^{v'}_t(u(xy^*)u^*)=u\vt^v_t(xy^*)u^*$. Observe that the stated $w_t$ are isometric (and, therefore, well-defined) and surjective. Since
\beqn{
(ux)(w_t(y^*\odot_tz))
~=~
(ux)((uy)^*\odot_t(uz))
~=~
\vt^{v'}_t((ux)(uy)^*)uz
~=~
u\vt^v_t(xy^*)z
~=~
u(x(y^*\odot_tz)),
}\eeqn
the $w_t$ fulfill \eqref{isouni}.\qed

\bdefi
In the situation of Proposition \ref{conisoprop} we say the pairs $(v,E^\odot)$ and $(v',E'^\odot)$ are \hl{conjugate}, and in the particular situation of Number \ref{CIP1} we say they are \hl{equivalent}.
\edefi

Roughly speaking, two left dilations of two product systems are conjugate (equivalent) if they induce inner conjugate (the same) \nbd{E_0}semigroup(s). In either case, the two product systems are necessarily isomorphic. The isomorphism is uniquely determined by \eqref{isouni}, once the unitary between the dilation spaces is fixed.

We see that associating different product systems with the same \nbd{E_0}semigroup means establishing a unique isomorphism between the product systems that behaves well with respect to the left dilations providing the association. This remains even true if the two \nbd{E_0}semigroups live on different but isomorphic $E$. But what happens if we have two ways $v$ and $v'$ to associate the \it{same} product system $E^\odot$ with the \it{same} \nbd{E_0}semigroup $\vt$ on $\sB^a(E)$?

In general, an endomorphism $w^\odot$ of $E^\odot$ induces a family $u_t:=v_t(\id_E\odot w_t)v_t^*$ of elements in $\sB^a(E)$ that form a local cocycle for $\vt=\vt^v$. (Exercise.) If $w_t$ is the automorphism that fulfills \eqref{isouni}, that is, that fulfills $v_t(x\odot y_t)=v'_t(x\odot w_ty_t)$, then we find $v'_t=v_t(\id_E\odot w_t^*)=u_t^*v_t$. In other words, $v$ and $v'$ are related by the local cocycle $u^*$ for $\vt=\vt^v$ and $\vt^{u^*}=\vt^{v'}=\vt^v=\vt$.

We have seen that two product systems associated with the same \nbd{E_0}semigroup or with inner conjugate \nbd{E_0}semigroups are isomorphic in an essentially unique way, and  we have seen the relation between two ways of associating the same product system with the same \nbd{E_0}semigroup. The next natural question is which \nbd{E_0}semigroups acting on a fixed $\sB^a(E)$ have the same or (equivalently, by the preceding discussion) isomorphic product systems. In other words, how are the \nbd{E_0}semigroups that act on a fixed $\sB^a(E)$ classified by their product systems. This is a generalization of the result \cite[Corollary of Definition 3.20]{Arv89} for Hilbert spaces and of \cite[Theorem 2.4]{Ske02} for Hilbert modules with a unit vector.

\bthm\label{ucethm}
Let $\vt$ and $\vt'$ be two strict \nbd{E_0}semigroups on $\sB^a(E)$ ($E$ a full Hilbert \nbd{\cB}module). Then their associated product systems $E^\odot$ and $E'^\odot$ are isomorphic if and only if $\vt$ and $\vt'$ are unitary cocycle equivalent.
\ethm

\proof
Denote by $v,v'$ the standard dilations of $E^\odot,E'^\odot$.

Suppose $w^\odot$ is a morphism from $E^\odot$ to $E'^\odot$. One checks that $u_t:=v'_t(\id_E\odot w_t)v_t^*$ defines a left cocycle with respect to $\vt$. The cocycle $u$ is unitary if and only if $w^\odot$ is an isomorphism. We find
\beqn{
u_t\vt_t(a)u_t^*
~=~
v'_t(\id_E\odot w_t)v_t^*v_t(a\odot\id_t)v_t^*v_t(\id_E\odot w_t^*)v'^*_t
~=~
v'_t(a\odot\id'_t)v'^*_t
~=~
\vt'_t(a),
}\eeqn
so that $\vt$ and $\vt'$ are cocycle equivalent.

Conversely, suppose $u$ is a unitary left cocycle such that $\vt^u=\vt'$. By $u_t\vt_t(a)=\vt'_t(a)u_t$, we see that $u_t$ is an isomorphism between the \nbd{\sB^a(E)}\nbd{\cB}correspondences $_{\vt_t}E$ and $_{\vt'_t}E$. It follows that $w_t:=\id_{E^*}\odot u_t$ defines a bilinear unitary from $E_t=E^*\odot_tE$ to $E'_t=E^*\odot'_tE$. We find
\bmun{
w_s(x^*\odot_sx')w_t(y^*\odot_ty')
~=~
(x^*\odot'_su_sx')(y^*\odot'_tu_ty')
~=~
x^*\odot'_{s+t}\vt'_t(u_sx'y^*)u_ty'
~=~
x^*\odot'_{s+t}u_t\vt_t(u_sx'y^*)y'
\\
~=~
x^*\odot'_{s+t}u_{s+t}\vt_t(x'y^*)y'
~=~
w_{s+t}(x^*\odot_{s+t}\vt_t(x'y^*)y')
~=~
w_{s+t}((x^*\odot_sx')(y^*\odot_ty')),
}\emun
so that the $w_t$ form a morphism.\qed

\bcor
Let $E$ and $E'$ be isomorphic full Hilbert \nbd{\cB}modules and suppose $\vt$ and $\vt'$ are two strict \nbd{E_0}semigroups on $\sB^a(E)$ and $\sB^a(E')$, respectively. Then their associated product systems $E^\odot$ and $E'^\odot$ are isomorphic if and only if $\vt$ and $\vt'$ are unitary cocycle inner conjugate.
\ecor

\proof
Fix a unitary $u\in\sB^a(E,E')$ and define the isomorphism $\alpha:=u\bullet u^*$. Then by Proposition \ref{conisoprop}\eqref{CIP3} $\vt$ and $\vt^\alpha$ have the same product system. The statement now follows by applying the theorem to $\vt^\alpha$ and $\vt'$.\qed

\lf
Before analyzing in the following sections the relation between product systems of \nbd{E_0}semi\-groups acting on not necessarily inner conjugate $\sB^a(E)$s in terms of more general notions of cocycle conjugacy, we close this section by giving a general result that does not involve cocycles. The result, well known for (separable) Hilbert spaces, provides a necessary and sufficient criterion for that all members of a family of \nbd{E_0}semigroups have isomorphic product systems. 

\bthm\label{Bclassthm}
Let $E^i$ $(i\in I)$ be a family of full Hilbert \nbd{\cB}modules and suppose that for each $i\in I$ we have a strict \nbd{E_0}semigroup $\vt^i$ on $\sB^a(E^i)$. Denote $E:=\bigoplus_{i\in I}E^i$.

Then the $\vt^i$ have mutually isomorphic product systems if and only if there exists a strict \nbd{E_0}semi\-group $\vt$ on $\sB^a(E)$ such that $\vt\upharpoonright\sB^a(E^i)=\vt^i$ for all $i\in I$. In the affirmative case, the product system of $\vt$ is in the same isomorphism class as those of the $\vt^i$.
\ethm

\proof
Denote by $p^i\in\sB^a(E)$ the projection onto $E^i$. We observe that an \nbd{E_0}semigroup $\vt$ on $\sB^a(E)$ leaves all $\sB^a(E^i)$ invariant if and only if $\vt_t(p^i)=p^i$ for all $i\in I,t\in\bS$.

Suppose all $\vt^i$ have isomorphic product systems. By Proposition \ref{conisoprop} we may fix one product system $E^\odot$ in this isomorphism class, and left dilations $v^i$ of $E^\odot$ to $E^i$ such that $\vt^i=\vt^{v^i}$ for all $i\in I$. Then $v_t:=\bigoplus_{i\in I}v^i_t$ defines a left dilation of $E^\odot$ to $E$. Clearly, the \nbd{E_0}semigroup $\vt:=\vt^v$ leaves all $\sB^a(E^i)$ invariant, and the restriction of $\vt$ to $\sB^a(E^i)$ is $\vt^i$.

On the contrary, suppose that $\vt$ is a strict \nbd{E_0}semigroup on $\sB^a(E)$ that leaves each $\sB^a(E^i)$ invariant. Suppose $E^\odot$ is a product system and $v$ is a left dilation of $E^\odot$ to $E$ such that $\vt^v=\vt$. Since $\vt_t(p^i)=p^i$, it follows that $v^i_t:=p^iv_t\upharpoonright(E^i\odot E_t)$ defines a left dilation of $E^\odot$ to $E^i$. Clearly, $\vt^{v^i}$ is just the restriction $\vt^i$ of $\vt$ to $\sB^a(E^i)$. Therefore, again by Proposition \ref{conisoprop} the product system of $\vt^i$ is isomorphic to $E^\odot$.\qed

\brem
The problem dealt with in the preceding proof is somewhat similar to showing that a functor between two categories of Hilbert modules is uniquely determined by what it does to a single full object. (On a single object, Theorem \ref{strirepthm} tells us that the functor is given by tensoring with a multiplicity correspondence. This is crucial for the proof of Blecher's \it{Eilenberg-Watts theorem} for Hilbert modules \cite{Ble97} from \cite[Section 2]{MSS06}.) In fact, here we are concerned with a semigroup of endofunctors (for each $t\in\bS$ induced by tensoring with the multiplicity correspondence $E_t$) leaving all objects fixed (that is, acting only on the morphisms) of the minicategory that has only the (full) objects $E,E^i~(i\in I)$. However, thanks to the simple structure (only full objects which are fixed by the functor) the direct proof we gave here is considerably simpler than reducing the statement to \cite[Section 2]{MSS06}.
\erem

\section{Conjugate $E_0$--Semigroups and Morita equivalent product systems}\label{E0MeSEC}

In the preceding section we showed that \nbd{E_0}semigroups on the same $\sB^a(E)$ have the same product system (up to isomorphism) if and only if they are unitary cocycle equivalent. This remains true if we replace unitary cocycle equivalence with unitary cocycle inner conjugacy, even between different $\sB^a(E)$ and $\sB^a(E')$, as long as $\sB^a(E)$ and $\sB^a(E')$ are inner conjugate, that is, as long as $E$ and $E'$ are isomorphic.

In this section we deal with the question what happens with \nbd{E_0}semigroups on two conjugate $\sB^a(E)$ and $\sB^a(F)$ under a strict conjugacy provided by an arbitrary strict isomorphism $\alpha\colon\sB^a(E)\rightarrow\sB^a(F)$ where $E$ is a Hilbert \nbd{\cB}module and $F$ is a Hilbert \nbd{\cC}module. Following Convention \ref{fconv}, we shall assume that $E$ and $F$ are full.

By Corollary \ref{Morcor}, there is a Morita equivalence $M$ from $\cB$ to $\cC$ such that $F\cong E\odot M$ and $E\cong F\odot M^*$, so that $E$ and $F$ are Morita equivalent. Moreover, $\alpha$ is the homomorphism implemented by the isomorphism $F\cong E\odot M$ and $\alpha^{-1}$ is the homomorphism implemented by the isomorphism $E\cong F\odot M^*$. Here, we wish to be more specific than making statements just up to isomorphism. We fix $M:=E^*\odot_\alpha F$ and $N:=F^*\odot_{\alpha^{-1}}E$ with identifications according to Theorem \ref{strirepthm}. For the isomorphisms in \eqref{Medefeq} we choose \eqref{homcomp}, that is,
\baln{
N\odot M
~\ni~
(y^*\odot_{\alpha^{-1}}x')\odot(x^*\odot_\alpha y')
&
~\longmapsto~
y^*\odot_{\id_F}\alpha(x'x^*)y'
~=~
\AB{y,\alpha(x'x^*)y'}
~\in~
\cC,
\\
M\odot N
~\ni~
(x^*\odot_\alpha y')\odot(y^*\odot_{\alpha^{-1}}x')
&
~\longmapsto~
x^*\odot_{\id_E}\alpha^{-1}(y'y^*)x'
~=~
\AB{x,\alpha^{-1}(y'y^*)x'}
~\in~
\cB.
}\ealn
By \cite[Theorem 1.14]{MSS06}, the identifications according to \eqref{homcomp} compose associatively. That is, we are in the situation required in Convention \ref{Massconv}. We easily check that for every $m:=x^*\odot_\alpha y\in M$ the element $n:=y^*\odot_{\alpha^{-1}}x\in N$ allows to recover $m^*\in M^*$ as $m'\mapsto n\odot m'\in N\odot M=\cC$. (We leave it as an instructive exercise to verify that the map $n\mapsto m^*$ is the map $u$ used in the proof of Theorem \ref{Meinvthm}.)

Morita equivalence of correspondences has been defined by Muhly and Solel \cite{MuSo00}. Recall the following version for product systems from \cite{Ske09}.

\bdefi\label{Mepsdef}
Let $E^\odot$ be a product system of \nbd{\cB}correspondences and let $M$ be a Morita equivalence from $\cB$ to $\cC$. Then the \hl{\nbd{M}transformed product system} of $E^\odot$ is the product system $M^*\odot E^\odot\odot M$ with $(M^*\odot E^\odot\odot M)_t:=M^*\odot E_t\odot M$ and identifications
\beqn{
M^*\odot E_{s+t}\odot M
~\cong~
M^*\odot E_s\odot E_t\odot M
~=~
(M^*\odot E_s\odot M)\odot(M^*\odot E_t\odot M).
}\eeqn
Clearly, $E^\odot\mapsto M^*\odot E^\odot\odot M$ and $w^\odot\mapsto\id_{M^*}\odot w^\odot\odot\id_M$ define an equivalence between the category of product systems of \nbd{\cB}correspondences and the category of product systems of \nbd{\cC}correspondences.

Two product systems $E^\odot$ and $F^\odot$ are \hl{Morita equivalent} if there exists a Morita equivalence $M$ such that $M^*\odot E^\odot\odot M$ and $F^\odot$ are isomorphic. Clearly, Morita equivalence of product systems is an equivalence relation.
\edefi

Putting all these identifications together and taking into account, once more, the associativity result \cite[Theorem 1.14]{MSS06} for the identifications according to \eqref{homcomp}, we immediately read off the following result.

\bprop
Let $E$ be a full Hilbert \nbd{\cB}module, let $\vt$ be a strict \nbd{E_0}semigroup on $\sB^a(E)$, and denote by $E^\odot$ the product system associated with $\vt$. Suppose $F$ is a full Hilbert \nbd{\cC}module with a strict isomorphism $\alpha\colon\sB^a(E)\rightarrow\sB^a(F)$, and denote by $M$ the associated Morita \nbd{\cB}\nbd{\cC}e\-quiv\-a\-lence (as discussed before).

Then the product system associated with $\vt^\alpha$ is isomorphic to $M^*\odot E^\odot\odot M$. In particular, the product systems of $\vt$ and of $\vt^\alpha$ are Morita equivalent.
\eprop

\proof
All we have to do is to write down a left dilation of $M^*\odot E^\odot\odot M$ to $F$ that gives $\vt^\alpha$. The chain $F\odot(M^*\odot E_t\odot M)\cong E\odot E_t\odot M\cong E\odot M\cong F$ of isomorphisms (where $F\odot M^*\cong E$ is the unique isomorphism from the proof of Corollary \ref{Morcor}, and where $E\odot E_t\cong E$ via the left dilation of $E^\odot$) does the job.\qed

\lf
As a simple corollary, we obtain the main classification result for \nbd{E_0}semigroup acting on strictly conjugate operator algebras.

\bthm\label{conjthm}
Let $\vt$ and $\theta$ be strict \nbd{E_0}semigroups on $\sB^a(E)$ and $\sB^a(F)$, respectively. Then $\vt$ and $\theta$ are unitary cocycle strictly conjugate if and only if there exists a strict isomorphism $\alpha\colon\sB^a(E)\rightarrow\sB^a(F)$ and their associated product systems are Morita equivalent via the same Morita equivalence inducing $\alpha$.
\ethm

\bex\label{nonidrem}
Note that for $F=E$ the notion of unitary cocycle strict conjugacy is strictly wider than the notion of unitary cocycle equivalence. We may suspect this, because the Picard group of $\cB$, in general, consists of more than multiplier inner automorphisms. But $M$ being a nontrivial Morita equivalence over $\cB$ such that $F=E\odot M$ does not yet guarantee that $M^*\odot E^\odot\odot M$ is not isomorphic to $E^\odot$. (For example, take the trivial product system $\cB^\odot$. But also the time ordered product systems $\DG^\odot(F)$ do not change under a transformation $M^*\odot\bullet\odot M$, whenever $M^*\odot F\odot M\cong F$.)

But we may obtain a concrete example in the following way. Let $F$ be a correspondence over $\cB$ and $M$ a Morita equivalence over $\cB$ such that $M^*\odot F\odot M\ncong F$. (Example \ref{nuvMex} helps. Indeed, we choose $M=E$ from that example, and $F=\C$, the \nbd{\C}component of $\cB=\C\oplus M_2$. We easily check that $M^*\odot F\odot M=M_2\ncong F$.) In that case, also the time ordered product systems $\DG^\odot(F)$ and $\DG^\odot(M^*\odot F\odot M)=M^*\odot\DG^\odot(F)\odot M$ (see \cite{BhSk00,LiSk01}) cannot be isomorphic, because the index $F$ of the time ordered product system $\DG^\odot(F)$ is an isomorphism invariant; see \cite{Ske06d}. Now $\DG^\odot(F)$ is the product system of the CCR-flow on $\sB^a(\DG(F))$ and $\DG^\odot(M^*\odot F\odot M)$ is the product system of the CCR-flow on $\sB^a(\DG(M^*\odot F\odot M))$. If $F$ and $M$ (like in the example) are countably generated, then $\DG(F)$ and $\DG(M^*\odot F\odot M)$ are stably isomorphic. So, the respective amplifications of the CCR-flows are unitary cocycle strictly conjugate. But they are not unitary cocycle inner conjugate, because the have non-isomorphic product systems. (Note that the discussion about stable conjugacy and amplification of \nbd{E_0}semigroups anticipates some arguments from Section \ref{sccSEC}; see there for details.)
\eex

We briefly specialize to the case, when the conjugacy of $\sB^a(E)$ and $\sB^a(F)$ can be chosen ternary. In that case, $M={_\vp}\cC$ where $\vp\colon\cB\rightarrow\cC$ is an isomorphism. One easily verifies that $M^*={_{\vp^{-1}}}\cB$ and that $M^*\odot E_t\odot M$ can be identified with $E_t$ via $b\odot x_t\odot c\mapsto bx_t\vp^{-1}(c)$ where, however, the inner product is $\AB{x_t,y_t}_\cC:=\vp(\AB{x_t,y_t})$ and the left action is $cx_t:=\vp^{-1}(c)x_t$.

We call \hl{ternary equivalent} two product systems that are Morita equivalent via a Morita equivalence $M$ that induces a ternary isomorphism. Theorem \ref{conjthm} remains true replacing `strictly conjugate' with `ternary conjugate' everywhere. Also Example \ref{nonidrem} remains valid in either direction:

\brem
The notion of unitary cocycle ternary conjugacy lies strictly in between unitary cocycle strict conjugacy and (where it applies) unitary strict cocycle equivalence. This follows from Example \ref{nonidrem}, from existence of non inner ternary isomorphisms, and from the observation that either composition of an isomorphism $\alpha$ with a ternary (an inner) isomorphism is ternary (inner) if and only if $\alpha$ is ternary (inner).
\erem

\brem\label{nunicocrem}
We think the potential of the translation of equations between homomorphisms into equations between multiplicity correspondences, as discussed in the beginning of Section \ref{E0cocSEC}, is by far not yet exhausted. It would be an interesting exercise to do the computations of Section \ref{E0cocSEC} in these terms. We did not do it in that way, because we do not gain simpler identifications, but rather a considerable complication concerning abstraction. A question where it appears unavoidable to proceed in that way is, what happens if we pass from unitary cocycles to arbitrary cocycles implementing the equivalence of \nbd{E_0}semigroups on the same $\sB^a(E)$. Already for ternary unitary cocycles we do not know the answer. (The main problem is that it is completely unclear what $\vt_s(u_t)$ might be for a ternary unitary $u_t$; see \cite[Remark 3.8]{Ske06c} and \cite[Section 4]{AbSk07}.) As for the present notes, we do not need an answer to this question. So, we do not tackle the problem here.
\erem


\section{Stable unitary cocycle (inner) conjugacy of $E_0$--semigroups}\label{sccSEC}

In the two preceding sections we established the main results about classification of \nbd{E_0}semi\-groups in terms of product systems in the situation where the involved \nbd{E_0}semigroups act on the same $\sB^a(E)$ or at least on strictly (or ternary) isomorphic operator algebras. (The only exception is Theorem \ref{Bclassthm}, which is, however, not in terms of cocycle conjugacy.) On the same $\sB^a(E)$ we found classification of \nbd{E_0}semigroups up to unitary cocycle equivalence by product systems up to isomorphism. For conjugate $\sB^a(E)$ and $\sB^a(F)$ we found classification of \nbd{E_0}semigroups up to unitary cocycle strict (ternary) conjugacy by product systems up to Morita (ternary) equivalence. We also showed that on the intersection of their domains, in general, the notions are all different. Only for Hilbert spaces the differences disappear.

But the question when two \nbd{E_0}semigroups have isomorphic or Morita equivalent product systems, has a meaning also if the \nbd{E_0}semigroups act on operator algebras of \it{a priori} unrelated Hilbert modules $E$ and $F$. In the present section we use the results from Section \ref{stabMeSEC} (and Section \ref{terSEC}) combining them with Sections \ref{E0cocSEC} and \ref{E0MeSEC} to answer this question under the (reasonable) countability conditions of Section \ref{stabMeSEC}. (The von Neumann case, without such countability conditions, will be discussed in Section \ref{vNalgSEC}.)

The problem is that before we can apply the results of Sections \ref{E0cocSEC} and \ref{E0MeSEC} in order to compare \nbd{E_0}semigroups on $\sB^a(E)$ and $\sB^a(F)$ in terms of their product systems, we must be able to compare $\sB^a(E)$ and $\sB^a(F)$ in terms of strict conjugacy. If $\sB^a(E)$ and $\sB^a(F)$ are not strictly conjugate then, maybe, $\sB^a(E\otimes\K)$ and $\sB^a(F\otimes\K)$ are. As soon as this is the case, we can apply Sections \ref{E0cocSEC} and \ref{E0MeSEC} to the amplified \nbd{E_0}semigroups. The following proposition, which is a simple corollary of Theorem \ref{Bclassthm}, taking into account that $E\otimes\K=\bigoplus_{n\in\N}E$, shows that passing to the stable versions does not change the product systems.

\bprop\label{ampPSprop}
Suppose $\vt$ is a strict \nbd{E_0}semigroup on $\sB^a(E)$ for some full Hilbert \nbd{\cB}module $E$. Then $\vt$ and its amplification $\vt^\K:=\vt\otimes\id_\sK$ to $\sB^a(E\otimes\K)$ have the same product system.
\eprop

\bdefi
Two strict \nbd{E_0}semigroup are \hl{stably cocycle conjugate} (\hl{equivalent}) if their amplifications are cocycle conjugate (equivalent). We use all supplements (like unitary, inner, ternary, etc.) in the same way as before.
\edefi


The following theorem merely collects most of the results of Sections \ref{stabMeSEC}, \ref{terSEC}, \ref{E0cocSEC} and \ref{E0MeSEC}.

\bthm\label{sclthm}
\begin{enumerate}
\item\label{s1}
If $\vt$ is a strict \nbd{E_0}semigroup on $\sB^a(E)$ for a full countably generated Hilbert module $E$ over a \nbd{\sigma}unital \nbd{C^*}algebra $\cB$, then the amplification $\vt^\K$ is inner conjugate to an \nbd{E_0}semigroup $\vt^\cB$ on $\sB^a(\K_\cB)$.

\item\label{s2}
Let $\vt$ and $\vt'$ be strict \nbd{E_0}semigroups on $\sB^a(E)$ and $\sB^a(E')$, respectively, where $E$ and $E'$ are full countably generated Hilbert modules over a \nbd{\sigma}unital \nbd{C^*}algebra $\cB$. Then the following conditions are equivalent:
\begin{enumerate}
\item\label{suc1}
$\vt$ and $\vt'$ are stably unitary cocycle inner conjugate.

\item\label{suc2}
$\vt^\cB$ and $\vt'^\cB$ are unitary cocycle equivalent.

\item\label{suc3}
$\vt$ and $\vt'$  have isomorphic product systems.
\end{enumerate}

\item\label{s3}
Let $\vt$ and $\theta$ be strict \nbd{E_0}semigroups on $\sB^a(E)$ and $\sB^a(F)$, respectively, where $E$ and $F$ are full countably generated Hilbert modules over \nbd{\sigma}unital \nbd{C^*}algebras $\cB$ and $\cC$, respectively. Then $\vt$ and $\theta$ are stably unitary cocycle strictly (ternary) conjugate if and only if they have Morita (ternary) equivalent product systems.
\end{enumerate}
\ethm

\proof
\ref{s1}. By Proposition \ref{Lanprop}, $E\otimes\K\cong\K_\cB$. Choose an isomorphism $u$. Then, conjugation of $\vt^\K$ with $\alpha:=u\bullet u^*$ gives the conjugate semigroup on $\sB^a(\K_\cB)$.

\ref{s2}. By definition, \ref{suc1} is equivalent to that $\vt^\K$ and $\vt'^\K$ are unitary cocycle inner conjugate. By Part \ref{s1}, this is the same as \ref{suc2}, and by Theorem \ref{ucethm} and Proposition \ref{ampPSprop} this is equivalent to \ref{suc3}.

\ref{s3}. By Theorem \ref{sMethm} (Theorem \ref{ter=isothm}), either condition means that $\cB$ and $\cC$ must be Morita equivalent (isomorphic). (Otherwise, none of the conditions can be satisfied.) So, if one of the conditions is satisfied, then there exists a Morita equivalence $M$ from $\cB$ to $\cC$. Since $\sK(M)=\cB$ is \nbd{\sigma}unital, by \cite[Proposition 6.7]{Lan95} $M$ is countably generated, and with $M$ and $E$ also the Hilbert \nbd{\cC}module $E\odot M$ is countably generated. We are now ready to apply Part \ref{s2} replacing $\vt$ with the conjugate \nbd{E_0}semigroup on $\sB^a(E\odot M)$. The specialization to the ternary case (where $M={_\vp}\cC$ for some isomorphism $\vp\colon\cB\rightarrow\cC$) is obvious.\qed

\section{About continuity}\label{contSEC}

So far, we have answered completely the question, to what extent strict \nbd{E_0}semigroups acting on operator algebras on countably generated Hilbert modules over \nbd{\sigma}unital \nbd{C^*}algebras are classified by their product systems (up to isomorphism, up to ternary equivalence, or up to Morita equivalence of the product systems). The answer is: Up to a suitable notion of stable unitary cocycle conjugacy. (The variation is just in the adjective preceding the word \it{conjugacy}.) For a complete treatment, there remains the problem to indicate which product systems can arise as product systems of \nbd{E_0}semigroups.

In this section we recall the known results about existence of \nbd{E_0}semigroups for product systems. We will see that in order that the constructed \nbd{E_0}semigroups live on spaces that are compatible with the countability assumptions (originating in Section \ref{stabMeSEC}), we can no longer avoid to introduce technical constraints both on the \nbd{E_0}semigroup side (strong continuity) and on the product system side (continuity and countability hypothesis).

First of all, recall that the product system of an \nbd{E_0}semigroup is always full. Recall, too, that in these notes we deal only with the case of \it{unital} \nbd{C^*}algebras. The case of not necessarily unital \nbd{C^*}algebras (with similar results, though technically more involved) is done in Skeide \cite{Ske09a}. For $\bS=\N_0$, by \cite[Theorem 7.6]{Ske09} we obtain an \nbd{E_0}semigroup that acts on the operators on a Hilbert module $E$. It is easy to check that $E$ is countably generated if and only if $E_1$ is countably generated as a (right) Hilbert module. As sketched only very briefly in \cite{Ske07}, in the case $\bS=\R_+$ but without continuity conditions, the algebraic part of the construction in the Hilbert space case in \cite{Ske06} generalizes (making also use of the result of \cite{Ske09}) easily to Hilbert modules. But the \nbd{E_0}semigroup obtained acts on the operators on a Hilbert module which is definitely not countably generated. Without continuity conditions, there is no construction known that would lead to a countably generated Hilbert module.

This negative statement ends the discussion of the continuous time case $\bS=\R_+$ in the purely algebraic situation.

\lf
Speaking about Hilbert modules, there remains the case $\bS=\R_+$ \it{with} continuity conditions. (The case of von Neumann modules will be discussed in Section \ref{vNcontSEC}.) To deal with that case, we have to repeat to some extent what these conditions are, and how the results from \cite{Ske07} allow to prove, as a new result, that suitable countability conditions on the \it{continuous} product system are preserved under the construction of an \nbd{E_0}semigroup.

\lf
Recall that an \nbd{E_0}semigroup $\vt$ on $\sB^a(E)$ is \hl{strongly continuous} if $t\mapsto\vt_t(a)x$ is continuous for all $a\in\sB^a(E)$ and $x\in E$. Obviously, the amplification $\vt^\K$ is strongly continuous if and only $\vt$ is strongly continuous. A family $u=\bfam{u_t}_{t\in\R_+}$ of elements $u_t\in\sB^a(E)$ (that is, in particular, a left cocycle) is \hl{strongly continuous} if $t\mapsto u_tx$ is continuous for all $x\in E$.

Following the definitions in \cite{Ske03b,Ske07}, a \hl{continuous product system} is a product system $E^\odot=\bfam{E_t}_{t\in\R_+}$ with a family $i_t\colon E_t\rightarrow\wh{E}$ of isometric embeddings of the \nbd{\cB}cor\-re\-spond\-ences $E_t$ into a common Hilbert \nbd{\cB}module $\wh{E}$ as right modules (there is no left action on $\wh{E}$), fulfilling the following conditions: Denote by
\beqn{
CS_i(E^\odot)
~=~
\BCB{\,\bfam{x_t}_{t\in\R_+}\colon x_t\in E_t,~t\mapsto i_tx_t\text{~is continuous~}}
}\eeqn
the set of \hl{continuous sections} of $E^\odot$ (with respect to the embeddings $i_t$). Then, firstly,
\beqn{
\BCB{\,x_s\colon\bfam{x_t}_{t\in\R_+}\in CS_i(E^\odot)\,}
~=~
E_s
}\eeqn
for all $s\in\R_+$ (that is, $E^\odot$ has sufficiently many continuous sections), and, secondly,
\beqn{
(s,t)
~\longmapsto~
i_{s+t}(x_sy_t)
}\eeqn
is continuous for all $\bfam{x_t}_{t\in\R_+},\bfam{y_t}_{t\in\R_+}\in CS_i(E^\odot)$ (that is, the `product' of continuous sections is continuous). A morphism between continuous product systems is \hl{continuous} if it sends continuous sections to continuous sections. An \hl{isomorphism} of continuous product systems is a continuous isomorphism. In \cite[Theorem 2.2]{Ske09a} we showed that a continuous isomorphism has continuous inverse. So, an isomorphism provides a bijection between the sets of continuous sections.

If $E$ is a full Hilbert module over a unital \nbd{C^*}algebra $\cB$, and if $\vt$ is a strongly continuous strict \nbd{E_0}semigroup on $\sB^a(E)$, then (generalizing on Skeide \cite{Ske03b}) we define a continuous structure on the associated product system $E^\odot$ in the following way:

Pass to the strongly continuous strict \nbd{E_0}semigroup $\theta:=\vt^\K$ on $\sB^a(F)$ for $F:=E\otimes\K$ and associated product system $F^\odot$. We know from Proposition \ref{ampPSprop} that $\vt$ and $\theta$ have the same product system. But  we wish to have more specific identifications. Obviously,
\beqn{
F_t
~=~
(E\otimes\K)^*\odot({_{\vt_t\otimes\sid_\K}}(E\otimes\K))
~=~
(E\otimes\K)^*\odot(({_{\vt_t}}E)\otimes\K)
~=~
(E^*\odot_tE)\otimes(\K^*\odot\K)
~=~
E_t\otimes\C
~=~
E_t,
}\eeqn
everywhere the canonical identifications and the natural action of $\sK$ on $\K$. (We leave it as an intriguing exercise for the reader to show that the concrete prescription
\beqn{
(x\otimes f)^*\odot_t(y\otimes g)
~\longmapsto~
(x^*\odot_ty)\AB{f,g}
}\eeqn
suggested by the preceding identifications, define bilinear unitaries $F_t\rightarrow E_t$ that form an isomorphism $F^\odot\rightarrow E^\odot$. The inverse is determined by $x^*\odot_ty\mapsto(x\otimes\vk)^*\odot_t(y\otimes\vk)$ where $\vk$ can be any unit vector in $\K$.)

Even if $E$ has no unit vector, by \cite[Lemma 3.2]{Ske09} $F$ has one, $\zeta$ say. That is, we are ready for the construction of the product system of $\theta$ following the first construction in \cite{Ske02} (imitating Bhat's construction \cite{Bha96} in the Hilbert space case) based on existence of a unit vector: For every $t\in\R_+$ define the Hilbert \nbd{\cB}submodule $\breve{F}_t:=\theta_t(\zeta\zeta^*)F$ of $F$. Turn it into a correspondence over $\cB$ by setting $b\breve{y}_t=\theta_t(\zeta b\zeta^*)\breve{y}_t$. Then $\breve{F}^\odot=\bfam{\breve{F}_t}_{t\in\R_+}$ is a product system via $\breve{y}_s\odot\breve{z}_t\mapsto\theta_t(\breve{y}_t\zeta^*)\breve{z}_t$ and $y\odot\breve{z}_t\mapsto\theta_t(y\zeta^*)\breve{z}_t$ defines a left dilation of $\breve{F}^\odot$ to $F$ giving back $\theta$. By Proposition \ref{conisoprop}, $F^\odot$ and $\breve{F}^\odot$ have to be isomorphic, and it is easy to verify directly that $y^*\odot_t z\mapsto\theta_t(\zeta y^*)z$ defines an isomorphism. (See the old version \cite[p.5]{Ske04p} for details.) We now define the isometric embedding
\beqn{
j_t
\colon
y^*\odot_tz
~\longmapsto~
\theta_t(\zeta y^*)z
~\in~
F
}\eeqn
of $F_t$ into $F$. It is easy to prove that this equips $F^\odot$ with a continuous structure. (See \cite{Ske03b,Ske07} for details. For instance, for every $y_t\in F_t\subset F$ the continuous section $\bfam{\theta_s(\zeta\zeta^*)y_t}_{s\in\R_+}$ meets $y_t$ for $s=t$.) By composing the isomorphism $E_t\rightarrow F_t$ with $j_t$ we define isometric embeddings $i_t\colon E_t\rightarrow F$, turning $E^\odot$ into a continuous product system isomorphic to $F^\odot$. It is $E^\odot$ with this continuous structure we have in mind if we speak about the continuous product system \hl{associated} with $\vt$.

It is noteworthy that the continuous structure does not depend on the choice of the unit vector $\zeta\in F$. In fact, if $\zeta'\in F$ is another unit vector, then the isomorphism $\theta_t(\zeta'\zeta^*)$ from the product system $\breve{F}^\odot$ constructed from $\zeta$ to the product system $\breve{F}'^\odot$ constructed from $\zeta'$ (see \cite[Proposition 2.3]{Ske02} for details), clearly, sends continuous sections to continuous sections, and so does its inverse $\theta_t(\zeta\zeta'^*)$. Even if $E$ has already a unit vector, $\xi$ say, and we started the construction from $\vt$ and that $\xi$, obtaining embeddings $E_t\rightarrow\vt_t(\xi\xi^*)\subset E$, the continuous structure would be the same. In fact, we may choose $\zeta=\xi\otimes\vk$ for a unit vector $\vk\in\K$. If, then, we identify $E$ with the subspace $E\otimes\vk$ of $F$, then the embeddings $E_t\rightarrow\vt_t(\xi\xi^*)E\subset E\rightarrow E\otimes\vk\subset F$ and $E_t\rightarrow F_t\rightarrow\theta(\zeta\zeta^*)F\subset F$ are the same.

In conclusion, the continuous structure of the product system associated with $\vt$ is determined uniquely by the preceding construction and is isomorphic to the continuous structure defined in \cite{Ske03b} in presence of a unit vector in $E$.

We show that the classification in Theorem \ref{ucethm} of \nbd{E_0}semigroups of a fixed $\sB^a(E)$ behaves well with respect to continuity. (Recall that the strong topology on unitaries coincides with the \nbd{*}strong topology. So, if the family $u$ of unitaries is strongly continuous, then so is the family of inverses. This makes cocycle equivalence by a strongly continuous cocycle an equivalence.)

\bthm\label{ucecthm}
Let $\vt$ and $\vt'$ be two strongly continuous strict \nbd{E_0}semigroups on $\sB^a(E)$ ($E$ a full Hilbert \nbd{\cB}module). Then their associated continuous product systems $E^\odot$ and $E'^\odot$ are isomorphic if and only if $\vt$ and $\vt'$ are unitary cocycle equivalent via a strongly continuous unitary left cocycle.
\ethm

\proof
If $\vt$ and $\vt'$ are not unitary cocycle equivalent, then by Theorem \ref{ucethm} $E^\odot$ and $E'^\odot$ are not even algebraically isomorphic. So, let us assume that $\vt$ and $\vt'$ are unitary cocycle equivalent and denote by $u$ and $w^\odot$ the unitary left cocycle and the isomorphism, respectively, related by the \it{formulae} in the proof of Theorem \ref{ucethm}. Then $v:=u\otimes\id_\K$ is a unitary left cocycle in $\sB^a(F)$ providing cocycle equivalence of $\theta$ and $\theta'$. Clearly, $v$ is strongly continuous if and only if $u$ is. Moreover, the isomorphism related to $v$ is the same $w^\odot$.

We shall show that $u\otimes\id_\K$ is strongly continuous if and only if $w^\odot$ is continuous. To that goal we switch to the product systems $\breve{F}^\odot$ and $\breve{F}'^\odot$, where continuity of sections is defined. So let $\breve{y}_t=\theta_t(\zeta\zeta^*)\breve{y}_t$ be a section of $\breve{F}^\odot$. Note that $v_t\breve{y}_t=v_t\theta_t(\zeta\zeta^*)\breve{y}_t=\theta'_t(\zeta\zeta^*)v_t\breve{y}_t\in\breve{F}'_t$, so that $v_t\breve{y}_t$ is a section of $\breve{F}'^\odot$. One easily verifies that this is precisely what $w^\odot$ does in the picture $\breve{F}^\odot\rightarrow\breve{F}'^\odot$.

So, if $v$ is strongly continuous, then with the section $\breve{y}$ also the section $v\breve{y}$ is continuous (and \it{vice versa}). That is, $w^\odot$ is continuous.

For the other direction, suppose $w^\odot$ is continuous. Choose $y\in F$ and a continuous section $\breve{z}$ of $\breve{F}^\odot$. Observe that with $\breve{z}$ also the function $t\mapsto\theta_t(y\zeta^*)\breve{z}_t$ is continuous. Likewise, this holds for $\theta'$ and continuous sections of $\breve{F}'^\odot$. So,
\bmun{
\norm{(v_t-v_s)\theta_t(y\zeta^*)\breve{z}_t}
~\le~
\norm{v_t\theta_t(y\zeta^*)\breve{z}_t-v_s\theta_s(y\zeta^*)\breve{z}_s}+\snorm{v_s\bfam{\theta_s(y\zeta^*)\breve{z}_s-\theta_t(y\zeta^*)\breve{z}_t}}
\\
~=~
\snorm{\theta'_t(y\zeta^*)v_t\breve{z}_t-\theta'_s(y\zeta^*)v_s\breve{z}_s}+\snorm{v_s\bfam{\theta_s(y\zeta^*)\breve{z}_s-\theta_t(y\zeta^*)\breve{z}_t}}
}\emun
is small for $s$ close to $t$, because $v_t\breve{z}_t$ is a continuous section of $\breve{F}'^\odot$. Since the set of all $\theta_t(y\zeta^*)\breve{z}_t$ is total in $F$, it follows that $v$ is strongly continuous.\qed

\lf
The basic result of \cite{Ske07} asserts that every full continuous product system of correspondences over a unital \nbd{C^*}algebra is (isomorphic to) the continuous product system associated with a strongly continuous strict \nbd{E_0}semigroup on some $\sB^a(E)$. (See, however, Remark \ref{gaprem}.) Theorem \ref{ucecthm} tells us that all strongly continuous strict \nbd{E_0}semigroups on that $\sB^a(E)$ are determined by the isomorphism class of their continuous product system up to continuous unitary cocycle equivalence. If, in order to complete the classification of strongly continuous strict \nbd{E_0}semigroups by continuous product systems, we wish to apply Theorem \ref{sclthm} (in particular, the equivalence of \ref{suc1} and \ref{suc3}), then we have to analyze to what extent we can guarantee that the \nbd{E_0}semigroup constructed in \cite{Ske07} lives on a $\sB^a(E)$ where $E$ satisfies the necessary countability hypotheses. Since $\cB$ is unital, that is, in particular, \nbd{\sigma}unital, we only have to worry about whether $E$ is countably generated.

Recall that starting from a strongly continuous \nbd{E_0}semigroup $\vt$ on $\sB^a(E)$, the module $\wh{E}$ is simply $F=E\otimes\K$. It seems, therefore, reasonable to require for the opposite direction that $\wh{E}$ is countably generated. But, in general, that would not even guarantee that the submodules $i_tE_t$ of $\wh{E}$ are countably generated. (Submodules of countably generated Hilbert modules need not be countably generated. Think of the \nbd{\sB(H)}submodule $\sK(H)$ of the singly generated \nbd{\sB(H)}module $\sB(H)$ for nonseparable $H$.) In addition, having a look at the construction in \cite{Ske07}, the question whether the constructed $E$ is countably generated reduces to the question whether the direct integral $\int_0^1E_\alpha\,d\alpha$ (defined in \cite{Ske07}) is countably generated. It is a submodule of the countably generated $L^2(\SB{0,1},\wh{E})$. However, once more submodules of countably generated Hilbert modules need not be countably generated. 

The problem disappears if we take into account that the submodules $\theta_t(\zeta\zeta^*)F$ of $F$ are the range of a projection. Indeed, if $E$ is countably generated, then so is $F$. Take a countable generating subset $S$ of $F$. Then the countable set of sections $\BCB{\bfam{\theta_t(\zeta\zeta^*)y}_{t\in\R_+}\colon y\in S}$ shows that $E^\odot$ is countably generated in the following sense. (A proof goes like that of \cite[Theorem 2.2]{Ske09a}.)
 
\bdefi
A continuous product system $E^\odot$ is \hl{countably generated} if it admits a countable subset of $CS_i(E^\odot)$ that is total in the locally uniform topology of $CS_i(E^\odot)$.
\edefi

\bthm\label{existcthm}
If $E^\odot$ is a countably generated continuous product system of correspondences over a unital \nbd{C^*}algebra $\cB$ then there exist a countably generated full Hilbert \nbd{\cB}module $E$ and a strongly continuous strict \nbd{E_0}semigroup on $\sB^a(E)$ such that $E^\odot$ is isomorphic to the continuous product systems associated with $\vt$.
\ethm

\proof
The continuous sections restricted to $\SB{0,1}$ take an inner product $\AB{x,y}=\int_0^1\AB{x_\alpha,y_\alpha}\,d\alpha$ and $\int_0^1E_\alpha\,d\alpha$ is defined as the norm completion. A countable set of sections generating the restriction to $\SB{0,1}$ in the uniform topology, is \it{a fortiori} generating for the $L^2$ topology.\qed

\lf
The classification theorem immediately follows:


\bthm\label{C*contclassthm}
Let $\cB$ be a unital \nbd{C^*}algebra. Then there is a one-to-one correspondence between equivalence classes (up to stable unitary cocycle inner conjugacy with strongly continuous unitary cocycles) of strongly continuous strict \nbd{E_0}semigroups acting on the operators of countably generated full Hilbert \nbd{\cB}mod\-ules and isomorphism classes of countably generated continuous product systems of correspondences over $\cB$.
\ethm

We dispense with stating the obvious variants for unitary cocycle strict or ternary conjugacy following from Theorem \ref{sclthm}\eqref{s3}. However:

\brem\label{MecPSrem}
The only problem is to make $M^*\odot E^\odot\odot M$ inherit a continuous structure from $E^\odot$, and to show that $M\odot(M^*\odot E^\odot\odot M)\odot M^*$ gives back the original structure. The simplest way, is to start with a strongly continuous \nbd{E_0}semigroup on some $\sB^a(E)$ that has product system $E^\odot$, and then to equip $M^*\odot E^\odot\odot M$ with the continuous structure that emerges from the conjugate \nbd{E_0}semigroup on $\sB^a(E\odot M)$. This assures that the structures are compatible, and that iterating we get back the original one.

A direct way to equip $M^*\odot E^\odot\odot M$ with a continuous structure, \it{without} making reference to \nbd{E_0}semigroups, is the following. Observe that, by \cite[Lemma 3.2]{Ske09}, there exists $n\in\N$ such that $M^n$ has a unit vector, $\mu$ say. Then $m^*\odot x_t\odot m'\mapsto\mu\odot m^*\odot x_t\odot m'$ defines an isometry $M^*\odot E_t\odot M\rightarrow M^n\odot M^*\odot E_t\odot M=(E_t\odot M)^n$. Combining this with $i_t\odot\id_M$ on each component, we get isometries $k_t\colon M^*\odot E_t\odot M\rightarrow(\wh{E}\odot M)^n$. It is not difficult to see that the $k_t$ turn $M^*\odot E^\odot\odot M$ into a continuous product system. (Approximate everything locally uniformly by sections that are finite linear combinations of sections $\bfam{m^*\odot y_t\odot m'}_{t\in\R_+}$ with $y\in CS_i(E^\odot)$ and $m,m'\in M$.) It also is not difficult to see that the double iteration of this procedure gives a new continuous structure on $E^\odot$ where all section from $CS_i(E^\odot)$ are continuous. By \cite[Theorem 2.2]{Ske09a}, the new and the old continuous structure coincide.

 Either way also works in the von Neumann case. (Just that in the second way, the finite cardinality $n\in\N$ must be allowed to be an arbitrary cardinal number.)
\erem

This ends the classification of \nbd{E_0}semigroups by product systems in the case of Hilbert \nbd{C^*}modules (under the manifest countability assumptions, of course). In the next section we apply the results to characterize the Markov semigroups that admit a special sort of dilation as the \it{spatial} ones. After that, the rest of the paper is devoted to discuss all the results for von Neumann algebras and modules. Apart from the absence of countability assumptions, in particular, the result about dilations is much more powerful, because there are much more interesting spatial Markov semigroups on von Neumann algebras.

\brem\label{gaprem}
We mention that the proof of \cite[Proposition 4.9]{Ske07}, which asserts that the \nbd{E_0}semigroup constructed in \cite{Ske07} for a continuous product system $E^\odot$, gives back the same continuous structure on $E^\odot$ we started with, has a gap. \cite[Theorem 2.2]{Ske09a} does not only fix that gap, but provides a considerably more general statement on the bundle structure of a continuous product system.
\erem

\newpage

\section{Hudson-Parthasarathy dilations of spatial  Markov semigroups}\label{HuPaSEC}

\it{Markov semigroups} are models for irreversible evolutions both of classical and of quantum systems. \it{Dilation} means to understand an irreversible evolution as \it{projection} from a reversible evolution of a \it{big} system onto the \it{small} subsystem via a conditional expectation. \it{Noises} are models for \it{big} systems in which the \it{small} system is \it{unperturbed}, that is, dilations of the \it{trivial} evolution of the \it{small} system or, yet in another way, a \it{big} physical system in \it{interaction picture} with the interaction \it{switched off}. Finding a dilation of the Markov semigroup that is a \it{cocycle perturbation} of a noise, means establishing a so-called \it{quantum Feynman-Kac formula}. If the perturbation is by a unitary left cocycle, then we speak of a \it{Hudson-Parthasarathy dilation}.

\it{Noises} come along with filtrations of subalgebras that are \it{conditionally independent} in the conditional expectation for some notion of \it{quantum independence}. That is why they are called \it{noises}. They also are direct quantum generalizations of \it{noises} in the sense of Tsirelson \cite{Tsi98p1,Tsi03p1}. In Skeide \cite{Ske04} we showed that every noise on $\sB^a(E)$ comes along with a \it{conditionally monotone independent} filtration. In Skeide \cite{Ske09r2,Ske08p2}, we show that every such \it{monotone noise} may be ``blown up'' to a \it{free noise} by making a relation with \it{free product systems}. See the survey Skeide \cite{Ske09r1} for more details.

\it{Spatiality} is a property that exists on the level of \nbd{E_0}semigroups, on the level of product systems, and on the level of CP-semigroups. If spatiality is present on one level, then it is present also on the other levels. Noises are spatial and, therefore, so are their product systems and all Markov semigroups that can be dilated to noises. Once more, see \cite{Ske09r1} for more details.

The scope of this section is to show by the means developed in the preceding sections that a Markov semigroup admits a Hudson-Parthasarathy dilation if and only if it is \it{spatial}. The key point is that, starting from the spatial Markov semigroup, we will construct two \nbd{E_0}semigroups having the same product system. One is another sort of dilation, a so-called \it{weak} dilation, while the other is a noise. So far, it was unclear how to compare these two \nbd{E_0}semigroups. But now, with the results obtained in the preceding sections, we know that (under countability conditions) their amplifications will act on the operator algebras of isomorphic Hilbert modules, so that there is a unitary left cocycle sending one amplification to the other. The only thing is to adjust the identification of the modules in such a way that they behave nicely in terms of the dilation. In \cite{Ske08p3} we performed that program for Markov semigroups on $\sB(H)$, where all the necessary classification results had already been known long before.

We start by explaining the terms that we used in the preceding introduction.

\lf
A \hl{CP-semigroup} is a semigroup $T=\bfam{T_t}_{t\in\bS}$ of completely positive maps $T_t$ on a \nbd{C^*}algebra $\cB$. In the sequel, we fix a unital \nbd{C^*}algebra $\cB$. A \hl{Markov semigroup} is a CP-semigroup $T$ where all $T_t$ are unital.

Suppose $(E,\vt,\xi)$ is a triple consisting of a Hilbert \nbd{\cB}module $E$, a strict \nbd{E_0}semigroup $\vt$ on $\sB^a(E)$, and a unit vector $\xi\in E$. Then by \cite[Proposition 3.1]{Ske02}, the family of maps $T_t\colon b\mapsto\AB{\xi,\vt_t(\xi b\xi^*)\xi}$ defines a CP-semigroup on $\cB$ (which is unital automatically) if and only if the projection $\xi\xi^*$ is \hl{increasing} for $\vt$, that is if an only if $\vt_t(\xi\xi^*)\ge\xi\xi^*$ for all $t\in\bS$. In this case, $(E,\vt,\xi)$ is a \hl{weak dilation} of $T$ in the sense of \cite{BhSk00}, that is, with the embedding $i\colon b\mapsto\xi b\xi^*$ and the vector expectation $\ep\colon a\mapsto\AB{\xi,a\xi}$ the diagram
\beqn{
\xymatrix{
\cB	\ar[d]_{i=\xi\bullet\xi^*}	\ar[r]^{T_t}	&\cB
\\
\sB^a(E)				\ar[r]_{\vt_t}	&\sB^a(E)	\ar[u]_{\ep=\AB{\xi,\bullet\xi}}
}
}\eeqn
commutes for all $t\in\bS$. A weak dilation is \hl{primary} if $\vt_t(\xi\xi^*)$ converges strongly to $\id_E$ for $t\to\infty$. If the diagram commutes with $i$ replaced by an arbitrary embedding, then we say just $(E,\vt,\xi,i)$ is a \hl{dilation}. A dilation $(E,\vt,\xi,i)$ is \hl{unital} if $i$ is unital. It is \hl{reversible} if $\vt$ consists of automorphisms. Note that whatever the dilation is, putting $t=0$ it follows that $i$ is injective and that $i\circ\ep$ is a conditional expectation onto $i(\cB)$. In the case of a unital dilation, this means that $i$ turns $E$ into a faithful correspondence over $\cB$. If we wish to think of $E$ as a correspondence in that way, we will identify $\cB$ as a unital subalgebra of $\sB^a(E)$ and write, slightly abusing notation, $(E,\vt,\xi,\id_\cB)$ for the unital dilation.

By a \hl{noise} over a unital \nbd{C^*}algebra $\cB$ we understand a triple $(E,\sS,\om)$ consisting of a (necessarily faithful) correspondence $E$ over $\cB$, an \nbd{E_0}semigroup $\sS$ on $\sB^a(E)$ (referred to as \hl{time shift}), and a unit vector $\om\in E$ (referred to as \hl{vacuum}), such that:
\begin{enumerate}
\item\label{N1}
$\sS$ leaves $\cB\subset\sB^a(E)$ \hl{pointwise invariant}, that is, $\sS_t(b)=b$ for all $t\in\R_+,b\in\cB$.

\item\label{N2}
$(E,\sS,\om,\id_\cB)$ is a unital dilation. (That is, with embedding $b=(x\mapsto bx)\in\sB^a(E)$.)

\item\label{N3}
$(E,\sS,\om)$ is a weak dilation. (That is, with embedding $b\mapsto\om b\om^*\in\sB^a(E)$.)
\end{enumerate}
By \ref{N2}, $\AB{\om,b\om}=b$ for all $b\in\cB$. Calculating the norm of $b\om-\om b$, it follows that $b\om=\om b$. By \ref{N1}, $(E,\sS,\om,\id_\cB)$ is a unital dilation of the trivial semigroup on $\cB$. By \ref{N3}, the projection $p:=\om\om^*$ is increasing. From
\beq{\label{basomcom}
\om b\om^*
~=~
\om\om^*b\om\om^*
}\eeqn
it follows that $(E,\sS,\om)$ is a weak dilation of the trivial semigroup.

\brem
This definition of noise is more or less from Skeide \cite{Ske06d}. In the scalar case (that is, $\cB=\C$) it corresponds to noises in the sense of Tsirelson \cite{Tsi98p1,Tsi03p1}. A reversible noise is closely related to a \it{Bernoulli shift} in the sense of Hellmich, Köstler and Kümmerer \cite{HKK04p}.
\erem

A noise is \hl{inner} and \hl{strongly continuous} and \hl{strict} if the time shift $\sS$ is inner and strongly continuous and strict, respectively. We use similar conventions for other properties of $\sS$, or of the weak dilation $(E,\sS,\om)$. For instance, a \hl{reversible noise} is a noise where $\sS$ consists of automorphisms. But, since there are, in general, noninner automorphism semigroups, a reversible noise need not be inner. An inner noise is \hl{vacuum preserving} if the implementing unitary semigroup $u$ can be chosen such that $u_t\om=\om\,(=u_t^*\om)$ for all $t\in\R_+$.

\bprop
An inner noise fulfills $\sS_t(\om\om^*)=\om\om^*$ for all $t\in\bS$. Moreover, the elements $\AB{\om,u_t\om}$ form a unitary semigroup in the center of $\cB$ such that the unitary semigroup $u'_t:=u_t\AB{u_t\om,\om}$ implements an inner noise that is vacuum preserving.
\eprop

\proof
For that $\sS_t(\om\om^*)=(u_t\om)(u_t\om)^*\ge\om\om^*$, it is necessary and sufficient that $\AB{u_t\om,\om}$ is an isometry. For that $\sS$ leaves $\cB$ invariant, it is necessary and sufficient, that all $u_t$ commute with all $b\in\cB$. It follows that also $u_t\om$ commutes with all $b\in\cB$. Therefore, $\AB{u_t\om,\om}$ is in the center of $\cB$. But an isometry in a commutative \nbd{*}algebra is a unitary. From this, also the inverse inequality $(u_t\om)(u_t\om)^*\le\om\om^*$ follows, so that $(u_t\om)(u_t\om)^*=\om\om^*$.

Observe that $u_t\om=(u_t\om)(u_t\om)^*(u_t\om)=\om\om^*(u_t\om)=\om\AB{\om,u_t\om}$. Applying $(u_{s+t}\om)(u_{s+t}\om)^*=(u_t\om)(u_t\om)^*$ to $u_{s+t}\om$, we find
\beqn{
u_{s+t}\om
~=~
(u_t\om)\AB{u_t\om,u_{s+t}\om}
~=~
(u_t\om)\AB{\om,u_s\om}
~=~
\om\AB{\om,u_t\om}\AB{\om,u_s\om}.
}\eeqn
Multiplying with $\om^*$ from the left, we see that the $\AB{\om,u_t\om}$ form a unitary semigroup in the center of $\cB$. The rest is obvious.\qed

\bob
The time shift $\sS_t$ differs from the modified time shift $\sS'_t:=u'_t\bullet{u'_t}^*$ by conjugation with the unitary semigroup $\AB{u_t\om,\om}$ in $\cB\subset\sB^a(E)$. By \cite[Theorem 4.2.18]{Ske01} (for instance), the center of $\cB$ is isomorphic to the center of $\sB^a(E)$ (acting on $E$ by multiplication from the \bf{right}), but it need not coincide with the center of $\sB^a(E)\supset\cB\supset C(\cB)$ (that is, acting from the \bf{left}). So $\sS'$ is, in general, different from $\sS$. But, since $\sS$ leaves $\cB\subset\sB^a(E)$ invariant, the unitaries $\AB{u_t\om,\om}$ form a (generally, nonlocal!) cocycle for $\sS$.
\eob

\bex
If $\cB$ has trivial center, for instance if $\cB=\sB(G)$ for some Hilbert space $G$, then the cocycle $\AB{u_t\om,\om}$ is local and does not change $\sS$. But, suppose $\cB\subset\sB(G)$ is a von Neumann algebra with nontrivial center $C(\cB)$. Put $E=G\,\ol{\otimes}^s\cB\subset\sB(G,G\otimes G)$ (exterior tensor product of von Neumann modules; see \cite[Section 4.3]{Ske01}). Then $\sB^a(E)=\sB(G)\,\ol{\otimes}^s\cB\subset\sB(G\otimes G)$ and the center of $\sB^a(E)$ is $\id_G\otimes C(\cB)$. We turn $E$ into a correspondence over $\cB$ by letting act $\cB$ on the factor $G$ of $E$. Clearly, conjugation with the left action of a unitary semigroup in $C(\cB)$ defines an automorphism semigroup leaving invariant the left action of $\cB$, but not $\sB(G)\otimes\id_G\subset\sB^a(E)$.
\eex

\bdefi
Let $T$ be a Markov semigroup on a unital \nbd{C^*}algebra $\cB$. A \hl{Hudson-Par\-tha\-sa\-rathy dilation} of $T$ is a noise $(E,\sS,\om)$ together with a unitary left cocycle $u$ with respect to $\sS$, such that $(E,\sS^u,\om,\id_\cB)$ becomes an (automatically unital) dilation of $T$. We shall often write $(E,\sS^u,\om)$ for a Hudson-Par\-tha\-sa\-rathy dilation.

A Hudson-Parthasarathy dilation is \hl{inner}, \hl{vacuum preserving}, and so forth, if the underlying noise is inner, vacuum preserving, and so forth. We will say the Hudson-Parthasarathy dilation is \hl{strongly continuous} if both the time shift $\sS$ and the cocycle $u$ are strongly continuous.

A Hudson-Par\-tha\-sa\-rathy dilation of $T$ is \hl{weak} if $(E,\sS^u,\om)$ is also a weak dilation (by \eqref{basomcom} necessarily of the same Markov semigroup $T$).
\edefi

Note that a Hudson-Par\-tha\-sa\-rathy dilation cannot be inner and weak at the same time. But we will see in Theorem \ref{inextthm} that every weak Hudson-Par\-tha\-sa\-rathy dilation arises as the restriction from an inner one.

\brem
The name Hudson-Parthasarathy dilation refers to the seminal work of Hudson and Parthasarathy \cite{HuPa84}. Perturbations of Markov semigroups by cocycles have been introduced by Accardi \cite{Acc78} under the name of \it{quantum Feynman-Kac formula}. Hudson and Parthasarathy \cite{HuPa84} constructed, for the first time, Hudson-Par\-tha\-sa\-rathy dilations for uniformly continuous Markov semigroups on $\sB(H)$ with a \it{Lindblad generator} of finite degree of freedom. The construction is with the help of their quantum stochastic calculus developed precisely for that purpose. Quantum stochastic calculus has been generalized to allow to find dilations of Markov semigroups with arbitrary Lindblad generator (Hudson and Parthasarathy \cite{HuPa84a}), unbounded versions (Chebotarev and Fagnola \cite{CheFa98}), and arbitrary von Neumann algebras (Goswami and Sinha \cite{GoSi99}). While the cited works all deal with $\sB(H)$ and more general von Neumann algebras, the quantum stochastic calculus in Skeide \cite{Ske00} deals completely within the \nbd{C^*}framework (and generalizes easily to von Neumann algebras).
\erem

We said that our results about classification of \nbd{E_0}semigroups up to stable cocycle conjugacy will allow to establish existence of the unitary cocycle of the Hudson-Parthasarathy dilation. As these results depend on continuity conditions, we switch immediately to sufficiently continuous Markov semigroups. For semigroups on unital \nbd{C^*}algebras, apart from the uniform topology, there is only the strong topology.%
\footnote{ \label{sgFN}
Weak continuity implies strong continuity; see Engel and Nagel \cite[Theorem I.1.6]{EnNa06}. Hille and Phillips \cite{HiPhi57} is still \it{the} source, in particular, for information about measurable semigroups and continuity at $0$. But the book \cite{EnNa06} is very concise and contains whatever we need easily detectable within the first 37 pages.
}
A semigroup $T$ of bounded linear maps on $\cB$ is \hl{strongly continuous} if $t\mapsto T_t(b)$ is continuous for all $b\in\cB$. We will see in a minute that in the \nbd{C^*}case such Markov semigroups that admit a Hudson-Par\-tha\-sa\-rathy dilations have bounded generators. (The von Neumann case is more interesting; see Section \ref{vNHPSEC}. It is even simpler in a sense, as it does not have some pathologies of the \nbd{C^*}case. Knowing the \nbd{C^*}case will help appreciating better the von Neumann case.) To understand this, we investigate better the product systems of the involved dilations.

\lf
Suppose $(E,\vt,\xi)$ is a strongly continuous weak dilation of an (automatically strongly continuous) Markov semigroup $T$. Then the projection $\xi\xi^*$ is increasing. If we construct the product system of $\vt$ with the unit vector construction (as described for $\theta$ on $\sB^a(E\otimes\K)$ on the pages preceding Theorem \ref{ucecthm}), then the $E_t=\vt_t(\xi\xi^*)E\ni\xi$ become an increasing family of subsets of $E$, all containing $\xi$. It is easy to check that the family $\xi^\odot=\bfam{\xi_t}_{t\in\bS}$ with $\xi_t:=\xi\in E_t$ form a \hl{unit}, that is, $\xi_s\xi_t=\xi_{s+t}$ and $\xi_0=\U$.%
\footnote{There is no unit defined for nonunital $\cB$. The condition $\xi_0=\U$ reflects that all our semigroups are actually monoids. In contexts with continuity, dropping the condition at $t=0$ would mean to speak about semigroups that are continuous only for $t>0$. It is well known that continuity at $t=0$ is often not automatic.}%
~The unit $\xi^\odot$ is even \hl{continuous} in that it is among the continuous sections of $E^\odot$. (After all, under the embedding into $E$ the section $\xi^\odot$ is constant.) Obviously, we recover $T$ from the unit $\xi^\odot$ as $T_t=\AB{\xi_t,\bullet\xi_t}$.

Now suppose, further, that $\vt=\sS^u_t$ is the perturbation by a strongly continuous cocycle of a strongly continuous noise $(E,\sS,\om=\xi)$. (In other words, suppose that $T$ admits a strongly continuous weak Hudson-Parthasarathy dilation.) Then the continuous product system of $\sS$ is also (isomorphic to) $E^\odot$. Since the noise is a weak dilation of the trivial CP-semigroup, its product system contains a continuous unit $\om^\odot=\bfam{\om_t}_{t\in\bS}$ such that $b=\AB{\om_t,b\om_t}$ for all $b\in\cB,t\in\bS$. One easily concludes that $b\om_t=\om_tb$, that is, the unit is \hl{central}. Moreover, the unit $\om^\odot$ is (like $\xi^\odot$) \hl{unital}, in the sense that all $\om_t$ are unit vectors. In Skeide \cite{Ske06d} we introduced \hl{spatial} product systems as pairs $(E^\odot,\om^\odot)$ consisting of a product system $E^\odot$ and a central unital \hl{reference unit} $\om^\odot$. We agree here to say a spatial product system is \hl{continuous} if $E^\odot$ is a continuous product system and if the reference unit $\om^\odot$ is among its continuous sections.

We just proved:

\bprop\label{spembprop}
If a Markov semigroup $T$ admits a strongly continuous weak Hudson-Par\-tha\-sa\-rathy dilation, then there is a continuous spatial product system $(E^\odot,\om^\odot)$ with a continuous unit $\xi^\odot$ such that $T_t=\AB{\xi_t,\bullet\xi_t}$.
\eprop

The statement that for every  CP-semigroup $T$ on a unital \nbd{C*}algebra there exists a product system $E^\odot$ with a unit $\xi^\odot$ such that $T_t=\AB{\xi_t,\bullet\xi_t}$, is not new. In fact, by a GNS-type construction, Bhat and Skeide \cite[Section 4]{BhSk00} construct a product system ${E^T}^\odot$ with a suitable unit $\xi^\odot$, the \hl{GNS-system} of $T$. The GNS-system is minimal in the sense that there is no proper subsystem containing the unit $\xi^\odot$, and the pair $({E^T}^\odot,\xi^\odot)$ is determined by these properties up to unit preserving isomorphism. The point about Proposition \ref{spembprop} is that the GNS-system of that Markov semigroup embeds continuously into a continuous spatial product system. After recalling the definition of a spatial strongly continuous Markov semigroup and a theorem from Bhat, Liebscher and Skeide \cite{BLS10}, this is equivalent to that $T$ is spatial.

\bdefi
(From \cite{BLS10} modeled after Arveson \cite{Arv97a}.) A \hl{unit} for a strongly continuous CP-semigroup $T$ on a unital \nbd{C^*}algebra $\cB$ is a continuous semigroup $c=\bfam{c_t}_{t\in\bS}$ of elements in $\cB$ such that $T_t$ \hl{dominates} the CP-map $b\mapsto c_t^*bc_t$ for all $t\in\bS$ (that is, for all $t\in\bS$ the map $T_t-c_t^*\bullet c_t$ is a CP-map). $T$ is \hl{spatial} if it admits units.
\edefi

CP-semigroups of the form $b\mapsto c_t^*bc_t$ are also called \hl{elementary} CP-semigroups. Continuity requirements for an elementary CP-semigroup refer to continuity of the implementing semigroup $c=\bfam{c_t}_{t\in\R_+}$.

\bitemp[Theorem \cite{BLS10}.~]\label{spthm}
Let $T$ be a strongly continuous CP-semigroup on a unital \nbd{C^*}al\-gebra. Then the following conditions are equivalent:
\begin{enumerate}
\item\label{spCP1}
$T$ is spatial.

\item\label{spCP2}
The GNS-system of $T$ embeds into a continuous spatial product system in such a way that the unit $\xi^\odot$ giving back $T$ is among the continuous sections.

\item\label{spCP3}
The generator $\cL$ of $T$ has \hl{Christensen-Evans form}, that is, $\cL(b)=\cL_0(b)+b\beta+\beta^*b$ for a CP-map $\cL_0$ on $\cB$ and an element $\beta$ of $\cB$.

\end{enumerate}
\eitemp

\brem
By \ref{spCP3} it follows, in particular, that a spatial CP-semigroup is uniformly continuous. It also follows by \cite[Theorem 6.3]{Ske06d} that the spatial product system into which the GNS-system embeds can be chosen to be a product system of time ordered Fock modules. Simply take the system generated by the two units $\om^\odot$ and $\xi^\odot$. But, the example in \cite{BLS10} shows that the GNS-system alone (that is, the subsystem generated by $\xi^\odot$) need not be spatial.
\erem

We sketch, very briefly, the proof from \cite{BLS10}. To show that a CP-semigroup fulfilling \ref{spCP2} is spatial, observe that $c_t:=\AB{\om_t,\xi_t}$ defines a semigroup in $\cB$. From $T_t-c_t^*\bullet c_t=\AB{\xi_t,q\bullet\xi_t}$ for the bilinear projection $q:=\id_t-\om_t\om_t^*$ in $\sB^a(E_t)$ we see that $T_t$ dominates $c_t^*\bullet c_t$. Moreover, since $b\mapsto\AB{\om_t,b\xi_t}=bc_t$ is strongly continuous, $c_t=\U c_t=\AB{\om_t,\U\xi_t}$ is norm continuous. So, $T$ is spatial. It is also easy to check that as soon as $\cL$ is bounded, the generator of $c$ provides a suitable $\beta$ as required for the Christensen-Evans form in \ref{spCP3}; see \cite[Lemma 5.1.1]{BBLS04}. Not so obvious is to see that $\cL$ is, indeed, bounded. (See \cite{BLS10} for the reduction to \cite[Theorem 7.7]{Ske03b}.) For the opposite direction, observe that if $T$ is a spatial CP-semigroup with unit $c$, say, then the maps
\beqn{
\fMatrix{b_{11}&b_{12}\\b_{21}&b_{22}}
~\longmapsto~
\fMatrix{b_{11}&b_{12}c_t\\c_t^*b_{21}&T_t(b_{22})}
}\eeqn
define a CP-semigroup on $M_2(\cB)$. Interpreting this in terms of so-called \it{CPD-semigroups}, by \cite[Theorem 4.3.5]{BBLS04}, there is a product system $E^\odot$ generated by a strongly continuous set of two units, one $\xi^\odot$ such that $\AB{\xi_t,\bullet\xi_t}=T_t$, the other a central unital unit $\om^\odot$ such that $\AB{\om_t,\xi_t}=c_t$. By Theorem \ref{cdilthm} in Appendix A.1, this product system is continuous with $\xi^\odot$ and $\om^\odot$ among the continuous sections. We refer to $(E^\odot,\xi^\odot,\om^\odot)$ as the \hl{spatial} continuous product system of $T$ \hl{associated} with the unit $c$.

Returning to our problem, Proposition \ref{spembprop} together with Theorem \ref{spthm} tells us that we must seek among the spatial Markov semigroups if we wish that they admit a weak Hudson-Par\-tha\-sa\-rathy dilation. We now wish to show that every spatial Markov semigroup admits such a dilation.

By Theorem \ref{spthm}, we may start with a continuous product system $E^\odot$ that has at least two unital units among its continuous sections. One is $\xi^\odot$ and generates $T$ as $T_t=\AB{\xi_t,\bullet\xi_t}$, the other the central unital reference unit $\om^\odot$. Already for Arveson systems it was known since \cite[Appendix]{Arv89} that a unital unit allows easily to construct an \nbd{E_0}semigroup. The construction for Hilbert modules is from \cite[Section 5]{BhSk00}: Take a product system $E^\odot$ and a unital unit $\zeta^\odot$. Embed $E_t$ into $E_{s+t}$ as $\zeta_sE_t$. The family of these embeddings forms an inductive system with inductive limit $E^\zeta$. The factorization $E_s\odot E_t=E_{s+t}$, under the limit, survives as $E^\zeta\odot E_t=E^\zeta$. In other words, we obtain a left dilation of $E^\odot$ to $E^\zeta$, inducing a strict \nbd{E_0}semigroup $\vt^\zeta$ on $\sB^a(E^\zeta)$. Moreover, $E^\zeta$ contains a unit vector $\zeta$ (the image of the elements $\zeta_t\in E_t\subset E^\zeta$) that factors as $\zeta=\zeta\zeta_t$ under the left dilation, and $(E^\zeta,\vt^\zeta,\zeta)$ is a weak dilation of the Markov semigroup $T^\zeta$ defined by $T^\zeta_t=\AB{\zeta_t,\bullet\zeta_t}$; see \cite{BhSk00,Ske02} for details.

\bob \label{primeob}
The weak dilation obtained via that inductive limit is primary and every primary weak dilation coincides with the inductive limit of its product system over its unit. (This is seen most easily in the picture $E_t=\vt_t(\zeta\zeta^*)E\subset E$.) This one-to-one correspondence survives in the continuous version in the next theorem. And, as explained below (after the theorem), if we restrict to unital central units, it survives to (continuous or not) primary noises.
\eob

\bitemp[Theorem {\cite[Theorem 7.5]{Ske03b}}.~] \label{indcontthm}
If $E^\odot$ is a continuous product system and if the unital unit $\zeta^\odot$ is among its continuous sections, then $\vt^\zeta$ is strongly continuous and the continuous structure on $E^\odot$ derived from $\vt^\zeta$ coincides with the original one. In particular, the continuous structure does not depend on the choice of $\zeta^\odot$.%
\footnote{\label{xiE0FN}In the proof of \cite[Theorem 7.5]{Ske03b} (with the unit denoted $\xi^\odot$ instead of $\zeta^\odot$), in proving that $\vt^\xi$ is strongly continuous, we were negligent regarding left continuity. However, in the case of the proof of \cite[Theorem 7.5]{Ske03b}, the omission is marginal. (To complete that proof, simply test strong left continuity at $t$, by checking on the total set of vectors $xy_t$ for some continuous section $y$, taking also into account that $y_t\approx y_{t-\ve}\xi_\ve\approx\xi_\ve y_{t-\ve}$.) The rest of the proof is okay.}
\eitemp

Constructing $E^\om$ and $\sS:=\vt^\om$ from $\om^\odot$, we obtain a weak dilation $(E^\om,\sS,\om)$ of the trivial semigroup. Since $\om^\odot$ is central, the left action of $\cB$ on $E_t$ survives the inductive limit ($b\om_sx_t=\om_sbx_t$). So, $E^\om$ with that left action becomes a correspondence over $\cB$, and the unit vector $\om$ fulfills $b\om=\om b$. Once more, by \eqref{basomcom} and since $\om\om^*$ is increasing for $\sS$, we see that $(E^\om,\sS,\om)$ is a strongly continuous noise. Moreover, $(E^\xi,\vt:=\vt^\xi,\xi)$ provides us with a strongly continuous weak dilation of $T$, sharing the product system with that noise.

The strategy is, like in Theorem \ref{sclthm} and its continuous version Theorem \ref{C*contclassthm}, to amplify the two \nbd{E_0}semigroup, appealing to that the modules $E^\om$ and $E^\xi$ are stably isomorphic, so that there will be a unitary cocycle. But for that, we must be sure that both modules are countably generated. Also, if we wish that the dilations are related somehow, then we have to make sure that also the amplified semigroups can be turned into a noise and a weak dilation of $T$, respectively, that are related in the sense of Hudson-Parthasarathy dilation.

We first look at the hypothesis to be countably generated.

\bprop
Let $\cB$ be a separable \nbd{C^*}algebra and let $E^\odot$ be a continuous product system of \nbd{\cB}correspondences. Suppose, further, that there is a countable set $S\subset CS_i(E^\odot)$ of units that generate $E^\odot$.

Then all $E_t$ are separable. Moreover, if $\zeta^\odot\in CS_i(E^\odot)$ is a unital unit, then also the inductive limit $E^\zeta$ is separable.
\eprop

\proof
One just has to observe that for each $t>0$ the set
\beq{\label{Etgenset}
\BCB{b_n\xi^n_{t_n}\ldots b_1\xi^1_{t_1}b_0\colon n\in\N,b_i\in\cB,{\xi^i}^\odot\in S,t_i>0,t_1+\ldots+t_n=t}
}\eeq
is total in $E_t$. (Every product subsystem of $E^\odot$ containing the units in $S$, must contain these elements. On the other hand, it is easy to check that the closed linear spans form a product subsystem; see \cite[Proposition 4.2.6]{BBLS04}. Since $E^\odot$ is generated by $S$, the subsystem must coincide with $E^\odot$.) Now, by continuity, the standard argument applies that in \eqref{Etgenset} the points $t_i$ can be restricted to the rational numbers and, of course, the elements $b_i$ can be restricted to a countable total subset of $\cB$, without changing totality of the set. This subset is, then, a countable union of countable sets and, therefore, a countable subset of \eqref{Etgenset}. So, $E_t$ is separable.

Moreover, the inductive limit of the $E_t$ over $t\in\R_+$ is increasing. It coincides, therefore, with the inductive limit of $E_n$ over $n\in\N_0$. So, also $E^\zeta$, as countable inductive limit over separable spaces, is separable.\qed

\lf
If $\cB$ is separable, we see that both $E^\om$ and $E^\xi$ are countably generated. So, in principle, we could now apply Theorem \ref{C*contclassthm}. But now we really have to worry about how to choose the identifications of the amplified modules $E^\om\otimes\K$ and $E^\xi\otimes\K$ in order that they behave nice with respect to the dilations carried by the original modules.

Let us start by observing that the inductive limit $E^\zeta$ obtained from a continuous unital unit $\zeta^\odot$ has the unit vector $\zeta$. In other words, the submodule $\zeta\cB$ is a direct summand of $E^\zeta$. So, from the two parts of the proof that the amplification $E^\zeta\otimes\K$ is isomorphic to $\K_\cB$, namely, Parts \ref{Lp1} and \ref{Lp2} of Proposition \ref{Lanprop}, we need only the second part. However, instead of applying Proposition \ref{Lanprop}\eqref{Lp2} directly to $E^\zeta\otimes\K$,  we, first, take away a piece. Like in the discussion in Section \ref{contSEC}, we choose a unit vector $\vk\in\K$. Then $\zeta_\vk:=\zeta\otimes\vk$ is a unit vector in $E^\zeta\otimes\K$. Moreover, $(E^\zeta\otimes\K,(\vt^\zeta)^\K,\zeta_\vk)$ remains a weak dilation of $T^\zeta$, sharing all the properties of the weak dilation $(E^\zeta,\vt^\zeta,\zeta)$ (except for that it is no longer primary). In particular, we know that the product system is the same. Now fix an isometry $v$ from $\CB{\vk}^\perp\subset\K$ onto $\K$. Then
\bmun{
E^\zeta\otimes\K
~=~
(\zeta\cB\otimes\K)\oplus(\CB{\zeta}^\perp\otimes\K)
~=~
(\zeta\cB\otimes\vk)\oplus(\zeta\cB\otimes\CB{\vk}^\perp)\oplus(\CB{\zeta}^\perp\otimes\K)
\\
~\cong~
(\zeta\cB\otimes\vk)\oplus(\zeta\cB\otimes\K)\oplus(\CB{\zeta}^\perp\otimes\K)
~=~
\zeta_\vk\cB\oplus(E^\zeta\otimes\K),
}\emun
where from the first to the second line we applied the isomorphism $\id_{\zeta\cB}\otimes v$ to the middle summand. Applying Proposition \ref{Lanprop}\eqref{Lp2} to the second summand of the last term, we obtain
\beqn{
E^\zeta\otimes\K
~\cong~
\zeta_\vk\cB\oplus\K_\cB.
}\eeqn
If we now do the same for $E^\om$ and $E^\xi$ we obtain
\beqn{
E^\om\otimes\K
~\cong~
\om_\vk\cB\oplus\K_\cB
~\cong~
\xi_\vk\cB\oplus\K_\cB
~\cong~
E^\xi\otimes\K.
}\eeqn
It is crucial to observe that this isomorphism identifies the distinguished unit vectors $\om_\vk$ and $\xi_\vk$. More precisely, we just have shown that there exists a unitary $u\colon E^\xi\otimes\K\rightarrow E^\om\otimes\K$ such that $u\xi_\vk=\om_\vk$.

By Theorem \ref{ucethm} there exists a strongly continuous unitary left cocycle $u_t$ with respect to $\sS^\K$  that fulfills
\beqn{
u\vt^\K_t(u^*au)u^*=u_t\sS^\K_t(a)u_t^*
}\eeqn
for all $a\in\sB^a(E^\om\otimes\K)$. We find
\bmun{
T_t(b)
~=~
\xi_\vk^*\vt^\K_t(\xi_\vk b\xi_\vk^*)\xi_\vk
~=~
\xi_\vk^*u^*u\vt^\K_t(u^*u\xi_\vk b\xi_\vk^*u^*u)u^*u\xi_\vk
\\
~=~
\om_\vk^*u\vt^\K_t(u^*\om_\vk b\om_\vk^*u)u^*\om_\vk
~=~
\om_\vk^*u_t\sS^\K_t(\om_\vk b\om_\vk^*)u_t^*\om_\vk,
}\emun
so that $u_t\sS^\K_t(\bullet)u_t^*$ with the unit vector $\om_\vk$ is a weak dilation of $T$. In particular, the projection $\om_\vk\om_\vk^*$ must be increasing, that is, $u_t\sS^\K_t(\om_\vk\om_\vk^*)u_t^*\om_\vk\om_\vk^*=\om_\vk\om_\vk^*$ or $u_t\sS^\K_t(\om_\vk\om_\vk^*)u_t^*\om_\vk=\om_\vk$. Now, recall that also $\om_\vk$ fulfills \eqref{basomcom}. It follows that
\beqn{
T_t(b)
~=~
\om_\vk^*u_t\sS^\K_t(\om_\vk\om_\vk^*)u_t^*u_t\sS^\K_t(b)u_t^*u_t\sS^\K_t(\om_\vk\om_\vk^*)u_t^*\om_\vk
~=~
\om_\vk^*u_t\sS^\K_t(b)u_t^*\om_\vk.
}\eeqn
In other words, the cocycle perturbation of the noise $(E^\om\otimes\K,\sS^\K,\om_\vk)$ by the cocycle $u_t$ is a weak Hudson-Parthasarathy dilation of $T$.

We collect what we have proved so far in the following characterization of Markov semigroups admitting weak Hudson-Parthasarathy dilations.

\bthm\label{C*HPthm}
Let $\cB$ be a separable unital \nbd{C^*}algebra and let $T$ be a Markov semigroup on $\cB$. Then $T$ admits a strongly continuous strict weak Hudson-Parthasarathy dilation if and only if $T$ is spatial.
\ethm 

By our construction, the correspondence on which the noise acts can be chosen to be isomorphic to $\K_\cB$ as right module. This does not at all mean that $\K_\cB$ would carry the canonical left action of $\cB$ that acts on each summand $\cB$ in $\K_\cB$ just by multiplication from the left. Also, due to the amplification procedure, the weak dilation of $T$ coming shipped with the weak Hudson-Parthasarathy dilation, in our construction will never be the minimal one. (There is a similarity of these facts to what happens in Goswami and Sinha \cite{GoSi99}.  There a Hudson-Parthasarathy dilation is constructed on $\cB\otimes\K$ where $\K$ is identified with a symmetric Fock space with the help of a quantum stochastic calculus. We mention, however, that the left action there is not even unital. Our construction improves this aspect.)

Apart from the mentioned problems with minimality, we can even say the following: There exist spatial Markov semigroups whose minimal weak dilation does not arise from a weak Hudson-Parthasarathy dilation. In fact, whenever the GNS-system of the spatial CP-semigroup is nonspatial, then a weak dilation obtained from a weak Hudson-Parthasarathy dilation is not minimal, because the product system of the dilation is spatial and, therefore, too big. An example is the counter example studied in \cite{BLS10}.

We close this long section on spatial CP-semigroups with the following result on inner Hudson-Parthasarathy dilations.

\bthm\label{inextthm}
~For every strongly continuous strict weak Hudson-Parthasarathy dilation $(E,\sS^u,\om)$ there exists a strongly continuous inner vacuum preserving Hudson-Parthasarathy dilation $(E',\sS'^{u'},\om')$ of the same Markov semigroup that ``contains'' $(E,\sS^u,\om)$ in the following sense:
\begin{enumerate}
\item
There is a strict unital representation of $\sB^a(E)$ on $E'$ that allows to identify $\sB^a(E)$ as a unital subalgebra of $\sB^a(E')$ (including the identification of $\cB\subset\sB^a(E)$ with $\cB\subset\sB^a(E')$).

\item
$\sS'$ leaves $\sB^a(E)\subset\sB^a(E')$ \hl{globally invariant} (that is, $\sS'_t(\sB^a(E))\subset\sB^a(E)$).

\item
$u'_t=u_t\in\sB^a(E)\subset\sB^a(E')$.
\end{enumerate}
\ethm

\noindent
This is the result of \cite{Ske07a} applied to the noise $(E,\sS,\om)$ ornamented by the embedding of the cocycle $u\subset\sB^a(E)$ into the bigger $\sB^a(E')\supset\sB^a(E)$. The algebraic properties are checked easily in the construction of \cite{Ske07a}. Continuity, a matter completely neglected in \cite{Ske07a}, follows very similarly as many other proofs of continuity like, for instance, continuity of $(E^\zeta,\vt^\zeta,\zeta)$. We do not give any detail.

\clearpage

\section{Von Neumann case: Algebraic classification}\label{vNalgSEC}

{\parskip0.5ex plus 0.5ex minus 0.5ex
For the balance of these notes (Appendix A.1 being the only exception) we discuss the analogues for von Neumann or \nbd{W^*}al\-ge\-bras (respectively, modules and correspondences) of the statements we obtained so far for \nbd{C^*}algebras and modules.

The algebraic part of the classification gets even simpler. This is mainly for two reasons. Firstly, a Morita equivalence $M$ from $\cA$ to $\cB$ relates the \nbd{C^*}algebra $\cA$ to the compact operators $\sK(M)$ on $M$. This obscures somehow that the representation theory of $\sB^a(E)$ is actually an operation of Morita equivalence, in that the statements that have interpretations in terms of Morita equivalence must be extended from $\sK(E)$ to $\sB^a(E)$ via strictness. In the von Neumann case this obstruction disappears, and the representations theory becomes pure Morita equivalence. Secondly, the stable versions of isomorphism results, in the \nbd{C^*}case, depend on Kasparov's stabilization theorem and, therefore, on countability assumptions on the modules. (The hypothesis of \nbd{\sigma}unitality of $\cB$, on the other hand, does not play a role, because $\cB$ is even assumed unital.) Also this obstruction disappears in the von Neumann case. (A small price to be paid is that now the dimension of the space with which a module must be stabilized depends on the module. But we are more than happy to pay that price, because we are payed off by getting a general theory without countability limitations.)

On the other hand, the weaker topologies of von Neumann objects require more work. Also this work is, however, payed off by much wider applicability of the results. (For instance, the results about Hudson-Parthasarathy dilations now apply to the large class of spatial Markov semigroups on a von Neumann algebra, which need no longer be uniformly continuous as in the \nbd{C^*}case.) For the first time, we give a concise definition of \it{strongly continuous product system}. (Let us state that we consider this definition only a ``working definition'', because it does not behave sufficiently ``nice'' with respect to the \it{commutant} of von Neumann correspondences. But it is enough for our purposes here. For instance, for strongly full and faithful product systems, like spatial ones, we obtain a full duality between product systems and their commutants; see Theorem \ref{scpsdualthm}.) Parts of the results for continuous product systems generalize more or less directly to strongly continuous versions. Other parts do not.

In the present section we repeat what we need to know about von Neumann modules. We specify the versions for von Neumann modules of the results of Sections \ref{MeSEC} -- \ref{sccSEC} (and, actually, a part of Section \ref{contSEC}) that go through without any further complication. This is the algebraic part of the classification. Actually, that part simplifies for von Neumann modules. In Section \ref{vNcontSEC} we deal with the analogues for the strong topology of the continuity results in Section \ref{contSEC}. In particular, we give our ``working definition'' of strongly continuous product systems. Here, most proofs go similar to the continuous \nbd{C^*}case or, at least, the necessary modifications are more or less obvious. An exception is the proof of existence of an \nbd{E_0}semigroup for every strongly full strongly continuous product system; a main result of these notes. Although, the strategy, in principle, is the same as for the \nbd{C^*}case, the technical differences are so substantial that we decided to discuss the proof in Appendix B. (The technical preparation allows us to prove in Appendix B.2 existence of faithful nondegenerate representations for a faithful strongly continuous product system. Arveson's proof of existence on an \nbd{E_0}semigroup for every Arveson system is done by establishing a faithful nondegenerate representation for every Arveson system. Our result allows the same for von Neumann correspondences; as explained in \cite[Section 9]{Ske06d}, this is \it{dual} to our treatment via the \it{commutant}, and it does not apply to the \nbd{C^*}case. In Appendix B, instead of \it{faithful nondegenerate representation}, we use the more recent terminology \it{right dilation} that underlines the \it{duality} with \it{left dilation} via \it{commutant}.) In Section \ref{vNHPSEC}, finally, we prove the results of Section \ref{HuPaSEC} about Hudson-Parthasarathy dilations for spatial Markov semigroups in the case of von Neumann algebras. In that section the differences become most substantial. It is necessary to fix in Appendix A.2 a gap in the proof of \cite[Theorem 12.1]{BhSk00} on dilations of Markov semigroups on von Neumann algebras. But, Appendix A.2 does more. Apart from presenting a more general version of \cite[Theorem 12.1]{BhSk00}, a result from \cite{MaSha10} is reproved as a corollary. We also prove the fundamental result that, in the von Neumann case, the GNS-system of a spatial Markov semigroup is spatial. This means a considerable simplification of the \nbd{C^*}case, where we only have embedding into a spatial product system. Appendix A.2 contains, thus, considerable parts of a beginning theory of \it{strong type I} product systems of von Neumann correspondences (that is, strongly continuous product systems of von Neumann correspondences that are generated by their strongly continuous units).

\lf\
The heart of Morita equivalence of \nbd{C^*}algebras is Example \ref{fundMeex} together with Corollary \ref{KE*cor}: What the inner products of elements of $E$ generate in norm coincides with the compact operators on $E^*$. All the rest is writing down suitable isomorphisms of certain \nbd{C^*}algebras with $\sK(E)$ or with $\cB_E$. For utilizing the relation between $\sK(E)$ and $\sB^a(E)$ in the representation theory, we had to work. In particular, we had to require that the representations are strict.

For von Neumann or \nbd{W^*}algebras and modules, once accepted the premise that all reasonable mappings between them be normal (or \nbd{\sigma}weak), everything is simpler. The range ideal $\cB_E$ of a von Neumann (or \nbd{W^*}) module over a von Neumann (or \nbd{W^*}) algebra $\cB$ will be replaced by its strong (or \nbd{\sigma}weak) closure $\ol{\cB_E}^s$. It coincides with the von Neumann (or \nbd{W^*}) algebra $\sB^a(E^*)$. The list of results or proofs involving Morita equivalence, where the proofs of the von Neumann version runs considerably more smoothly than that of the \nbd{C^*}version (or where the \nbd{C^*}version even fails), is still getting longer. We resist the temptation to give such a list and refer the reader to \cite{Rie74a,Ske09,Ske06b}.

In order to avoid the notorious distinction between von Neumann and \nbd{W^*}modules, we have to make a decision. Although they form equivalent categories, von Neumann modules \cite{Ske00b,Ske06b} are technically simpler. In fact, many proofs  of results about \nbd{W^*}modules run best, after transforming the modules into von Neumann modules by choosing a faithful normal representation of the involved \nbd{W^*}algebras. Other proofs do not even possess intrinsic \nbd{W^*}versions. However, the notion of \nbd{W^*}modules is more wide-spread. Our decision is: We shall formulate for \nbd{W^*}algebras, \nbd{W^*}modules, and \nbd{W^*}correspondences; but we will always tacitly assume that a \nbd{W^*}algebra $\cB$ is acting concretely as a von Neumann algebra $\cB=\cB''\subset\sB(G)$ on a Hilbert space $G$. (Note that $\cB=\cB''$ this implies that the identification map is normal and that $\U_\cB=\id_G$.) It is well known (see Rieffel \cite{Rie74a} or Murphy \cite{Mur97} or Skeide \cite{Ske00b}) that, once a (pre-)\nbd{C^*}algebra $\cB$ is identified as a concrete operator algebra $\cB\subset\sB(G)$, every (pre-)Hilbert \nbd{\cB}module $E$ is identified as a subset of $\sB(G,H)$ where the Hilbert space $H$ is $E\odot G$ and where $x\in E$ acts as $x\colon g\mapsto x\odot g$. Moreover, one directly checks that if $E\subset\sB(G,H')$ is another identification fulfilling (i) $x\circ b=xb$, (ii) $x^*y=\AB{x,y}$, and (iii) $\cls EG=H'$, then $H'$ and $H$ are unitarily equivalent identifying $xg\in H'$ with $x\odot g\in H$. Now, a pre-Hilbert module $E$ over a von Neumann algebra $\cB\subset\sB(G)$ is a \hl{von Neumann \nbd{\cB}module} if $E$ is strongly closed in $\sB(G,H)$; \cite[Definition 4.4+Proposition 4.5]{Ske00b}. By \cite[Theorem 5.5]{Ske00b}, $E$ is a von Neumann \nbd{\cB}module if and only if it is self-dual, that is, if and only if it is a \hl{\nbd{W^*}module}. (Note that self-duality is a property of $E$ that does not depend on the representation $\cB$.) Moreover, since the strong closure of $E$, $\ol{E}^s$, is a von Neumann module and since it contains $E$ as a strongly dense subset, $\ol{E}^s$ is the unique minimal \hl{self-dual extension} (Paschke \cite{Pas73}) of $E$ to a \nbd{W^*}module. (It coincides with $\sB^r(E,\cB)$, but the problem in \cite{Pas73} is how to compute inner products of elements from $\sB^r(E,\cB)$.) Together with the identification of $E\subset\sB(G,H)$ there comes the identification $\sB^a(E)=\sB^a(E)\odot\id_G\in\sB(H)$. If $E$ is a von Neumann module, then $\sB^a(E)$ is a \nbd{W^*}algebra, always acting as von Neumann algebra on $H=E\odot\id_G$; see \cite[Proposition 4.5]{Ske00b} (\cite[Proposition 3.10]{Pas73}). A \hl{\nbd{W^*}correspondence} from $\cA$ to $\cB$ is a \nbd{W^*}module over $\cB$ with a nondegenerate and suitably normal left action of $\cA$. The \hl{tensor product} of \nbd{W^*}correspondences $E$ and $F$, denoted by $E\sodots F$, is the self-dual extension of $E\odot F$). Every subset $S$ of a \nbd{W^*}module $E$ generates a \nbd{W^*}submodule of $E$ that coincides with the orthogonal bicomplement of $S$ in $E$. An (algebraic) submodule of $E$ is strongly total if its bicomplement is $E$. It is dense in $E$ in the natural \nbd{\sigma}weak topology of $E$, respectively, in the strong topology when considered as von Neumann module. The strong topology on a \nbd{W^*}module, like the strong topology on the \nbd{W^*}algebra itself, depends on the choice of the representation; but our results do not. The \nbd{\sigma}weak topology is the relative \nbd{\sigma}weak topology from the linking von Neumann algebra and, therefore, unique. A \nbd{W^*}module $E$ over $\cB$ is \hl{strongly full} if $\ol{\cB_E}^s=\cB$.
}

\lf
Let us start with the notions and results of Section \ref{MeSEC}.

\newpage
\bdefi
A \nbd{W^*}correspondence $M$ from $\cA$ to $\cB$ is a \hl{Morita \nbd{W^*}equivalence} if there exists a \nbd{W^*}correspondence $N$ from $\cB$ to $\cA$ such that
\baln{
N\sodots M
&
~\cong~
\cB,
&
M\sodots N
&
~\cong~
\cA,
}\ealn
as \nbd{W^*}correspondences over $\cB$ and over $\cA$, respectively. Also here, we call $N$ an \hl{inverse} of $M$ under tensor product.

Following Rieffel \cite{Rie74a}, two \nbd{W^*}al\-ge\-bras $\cA$ and $\cB$ are \hl{Morita equivalent} if there exists a Morita \nbd{W^*}equivalence from $\cA$ to $\cB$.
\edefi

Also here, a Morita equivalence is necessarily faithful and strongly full. Since $\sK(E)$ is strictly dense dense in $\sB^a(E)$ (a bounded approximate unit for $\sK(E)$, by definition, converges \nbd{*}strongly, and, therefore, strictly, to the identity), $\sK(E)$ is, \it{a fortiori}, dense in the \nbd{*}strong topology of $\rho(\sB^a(E))\subset\sB(K)$ for whatever (normal or not) representation $\rho$ on a Hilbert space $K$. Consequently, $E\sodots E^*=\ol{\sK(E)}^s=\sB^a(E)$. So, a strongly full \nbd{W^*}module over $\cB$ is a Morita \nbd{W^*}equivalence from $\sB^a(E)$ to $\cB$. (Clearly, the finite-dimensional spaces in Example \ref{nuvMex}, are also \nbd{W^*}objects.)

Considerations like these show that the results that follow in Section \ref{MeSEC} (including their proofs) remain true until Theorem \ref{strirepthm} if we replace everywhere:
\begin{itemize}
\item
\nbd{C^*}Algebras with \nbd{W^*}algebras;

\item
Hilbert modules with \nbd{W^*}modules;

\item
correspondences with \nbd{W^*}correspondences;

\item
tensor products of correspondences with the tensor products of \nbd{W^*}correspondences;

\item
range ideals with their strong closures;

\item
full with strongly full;

\item
$\sK(E)$ with $\sB^a(E)$;

\item
strict maps with normal maps.
\end{itemize}
For instance, the considerations we did above are the direct explanation why Corollary \ref{KE*cor} remains true under these changes. In particular, we have the compatibility result for tensor products in Proposition \ref{assdiagprop} and Convention \ref{Massconv} ($M\sodots(N\sodots M)=M=(M\sodots N)\sodots M$), we have the characterization of Morita equivalences by Theorem \ref{Meinvthm} (the left action of a Morita \nbd{W^*}equivalence from $\cA$ to $\cB$ defines an isomorphism from $\cA$ onto $\sB^a(M)$) and Convention \ref{N=M*conv} ($N=M^*$), and we have Corollary \ref{Krepcor} ($F\cong\sB^a(E)\sodots F\cong(E\sodots E^*)\sodots F\cong E\sodots(E^*\sodots F)$ as \nbd{W^*}correspondences from $\sB^a(E)$ to $\cC$). All these results are essentially algebraic; just that now the extensions are done in a weaker topology. To see that this goes without any problems, it is key that any normal representation of $\cB$ extends to a normal representation of the linking von Neumann algebra of a \nbd{W^*}module (this can be deduced, for instance, from the remark following Muhly and Solel \cite[Lemma 2.16]{MuSo02}, but see also \cite[Lemma 3.3.2]{Ske01}, whose proof illustrates applicable techniques), the fact that normal representations of von Neumann algebras are strongly continuous (this follows, for instance, from suitable applications of \cite[Lemma II.2.5 and Theorem II.2.6]{Tak02}), and the fact that the operations when viewed as operator multiplication in the von Neumann linking algebra, share all the ``good'' properties of strong continuity of operator multiplication (see, for instance, \cite[Proposition 4.6]{Ske00b}).

While in deducing Theorem \ref{strirepthm} from Corollary \ref{Krepcor}, we had to discuss an extension from $\sK(E)$ to $\sB^a(E)$ by strictness, now the \nbd{W^*}version of Corollary \ref{Krepcor} gives us directly  the \nbd{W^*}version of Theorem \ref{strirepthm} (normal nondegenerate representations of $\sB^a(E)$ on a \nbd{W^*}module $F$ are amplifications of the identity representation with the \nbd{W^*}correspondence $E^*\sodots F$; see \cite[Theorem 1.16]{MSS06}) by simply fixing the isomorphism $u$.

Isomorphisms between von Neumann algebras are normal, automatically. Therefore, in the modified version of Corollary \ref{Morcor} we may leave out the word `normal' (which, according to our rules, has substituted the word `strict'). The same is true for Corollary \ref{striMocor}, once we stated the following:

\bdefi
A \nbd{W^*}module $E$ over $\cB$ and a \nbd{W^*}module $F$ over $\cC$ are \hl{Morita equivalent} if there is a Morita \nbd{W^*}equivalence $M$ from $\cB$ to $\cC$ such that $E\sodots M\cong F$ (or $F\sodots M^*\cong E$).
\edefi

So, we get the \nbd{W^*}version of Corollary \ref{striMocor}, which is worth a theorem:

\bthm
Strictly full \nbd{W^*}modules have isomorphic operator algebras if and only if they are Morita equivalent.
\ethm

Corollary \ref{bistrcor} remains true, independently, in its original formulation. Apart from that the stated isomorphism has no choice but being normal, the corollary states a criterion for when it is also strict (which may happen or not).

\lf
Section \ref{stabMeSEC} has to be overworked considerably. The stabilization results in the von Neumann context are more general, but they depend on the choice of a sufficiently big cardinal number $\en$. Their proofs are completely different and less sophisticated than their \nbd{C^*}counterparts. We take them mainly from \cite{Ske09}, where it is also pointed out that these facts resemble statements from the representation theory of von Neumann algebras.

Given a cardinal number $\en$, by $\C^\en$ we denote the canonical Hilbert space of dimension $\en$. We denote $M_\en^s:=\sB(\C^\en)$. By $\ol{E^\en}^s:=E\sbars{\otimes}\C^\en$ we mean the von Neumann or \nbd{W^*}version of the exterior tensor product; see \cite[Section 4.3]{Ske01} for details and for the facts we are using in the sequel. We have $\sB^a(E\sbars{\otimes}\C^\en)=\sB^a(E)\sbars{\otimes}M_\en^s$ (tensor product of von Neumann algebras). For an infinite cardinal number $\en$, we say \nbd{W^*}algebras $\cA$ and $\cB$ are \hl{\nbd{\en}stably isomorphic} if $\cA\sbars{\otimes}M_\en^s$ and $\cB\sbars{\otimes}M_\en^s$ are isomorphic. We say $\cA$ and $\cB$ are \hl{stably isomorphic} if they are \nbd{\en}stably isomorphic for some infinite cardinal number $\en$. Since $\ol{\cB^\en}^s$ is a Morita \nbd{W^*}equivalence from $\cB\sbars{\otimes}M_\en^s$ to $\cB$ it follows that \nbd{W^*}algebras $\cA$ and $\cB$ are Morita equivalent if they are stably isomorphic. By \cite[Corollary 9.4]{Ske09}, also the converse is true.

Proposition \ref{Lanprop} gets the following shape.

\bprop\label{vNLanprop}
Suppose $E$ is a strongly full \nbd{W^*}module over $\cB$. Then there exists a cardinal number $\en$ such that:
\begin{enumerate}
\item\label{vNLp1}
$\ol{E^\en}^s$ has a direct summand $\cB$.

\item\label{vNLp2}
$\ol{E^\en}^s\cong\ol{\cB^\en}^s$.
\end{enumerate}
\eprop
Part \ref{vNLp1} is \cite[Lemma 4.2]{Ske09}. Part \ref{vNLp2} is stated and proved in front of \cite[Corollary 4.3]{Ske09}. Of course, \ref{vNLp2} implies \ref{vNLp1}. But, like in the proof of Proposition \ref{Lanprop}, Part \ref{vNLp2} is proved using Part \ref{vNLp1}. It may be noted that $\en$ cannot always be chosen to be the smallest cardinality of a subset that generates $E$ as a \nbd{W^*}module:

\bex
Let $H$ be a nonseparable Hilbert space and choose a nonzero vector $h\in H$. Then the strongly full \nbd{W^*}module $H^*$ over $\cB=\sB(H)$ is generated by the single element $h^*$. But no cardinality $\en$ strictly smaller than $\dim H$ makes $\ol{{H^*}^\en}^s$ isomorphic to $\ol{\cB^\en}^s$. In fact, $\ol{{H^*}^\en}^s=\sB(H,\C^\en)$ does not contain a single copy of $\cB$, because it contains only operators of ``rank'' not greater than $\en$. In particular, it does not contain any unit vector.
\eex

\bdefi
For an infinite cardinal number $\en$, two \nbd{W^*}modules $E$ and $F$ are \hl{\nbd{\en}stably Morita equivalent} if $\ol{E^\en}^s$ and $\ol{F^\en}^s$ are Morita equivalent. They are \hl{stably Morita equivalent} (as \nbd{W^*}modules) if they are \nbd{\en}stably Morita equivalent for some infinite cardinal number $\en$.
\edefi

\bob
Suppose $\ol{E^\ek}^s$ and $\ol{F^\el}^s$ are Morita equivalent for some arbitrary (also finite) cardinal numbers $\ek$ and $\el$. Then $E$ and $F$ are \nbd{\en}stably Morita equivalent for every infinite cardinal number $\en\ge\max(\ek,\el)$. (Simply choose isomorphisms $\C^\ek\otimes\C^\en\cong\C^\en\cong\C^\el\otimes\C^\en$.) The same is true, of course, for stable isomorphism of \nbd{W^*}algebras.
\eob

The analogue of Theorem \ref{sMethm} reads as follows.

\bthm\label{vNsMethm}
Let $E$ and $F$ denote strongly full \nbd{W^*}modules over \nbd{W^*}algebras $\cB$ and $\cC$, respectively. Then the following conditions are equivalent:
\begin{enumerate}
\item\label{vNsM1}
$E$ and $F$ are stably Morita equivalent.

\item\label{vNsM2}
$\sB^a(E)$ and $\sB^a(F)$ are Morita equivalent.

\item\label{vNsM3}
$\cB$ and $\cC$ are Morita equivalent.

\item\label{vNsM4}
$\cB$ and $\cC$ are stably isomorphic.

\end{enumerate}
\ethm

\proof
Since $\sB^a(E)$ and $\cB$ are Morita equivalent (similarly, for $\sB^a(F)$ and $\cC$) and since Morita equivalence is an equivalence relation, \ref{vNsM2} and \ref{vNsM3} are equivalent. Equivalence of \ref{vNsM3} and \ref{vNsM4} is \cite[Corollary 9.4]{Ske09}. Of course, \ref{vNsM1} implies \ref{vNsM3}; and if \ref{vNsM3} holds, then by Proposition \ref{vNLanprop}\eqref{vNLp2} also \ref{vNsM1} holds, so that also \ref{vNsM1} and \ref{vNsM3} are equivalent.
\qed

\lf
To save space we do not spend much time on ternary isomorphisms, because everything is quite obvious. We mention only one thing, which facilitates to understand why everything is obvious. A ternary homomorphism between \nbd{W^*}modules extends to a normal homomorphism between their linking von Neumann algebras if and only if it is \nbd{\sigma}weak. For that this happens, it is already sufficient that the restriction of the extension to the corner $\ol{\cB_E}^s$ or to the corner $\sB^a(E)$ is normal. (See again  the remark following \cite[Lemma 2.16]{MuSo02}.) With this observation, everything in Section \ref{terSEC} goes through for the obvious modifications. In particular, the \nbd{W^*}version of Theorem \ref{ter=isothm} asserts that strongly full \nbd{W^*}modules are stably ternary isomorphic if and only if they are modules over isomorphic \nbd{W^*}algebras.

Section \ref{ccSEC}, of course, remains unchanged, as it is completely on algebras without any modules or topologies.

With the same global substitutions as for Section \ref{MeSEC}, also Section \ref{E0cocSEC} remains essentially unchanged. Only in Theorem \ref{Bclassthm} we have to replace the direct sum with the \nbd{W^*}module direct sum. The same is true for Section \ref{E0MeSEC} and \it{cum grano salis} also for Section \ref{sccSEC}.%
\footnote{\label{repMEfn}We dispense with giving a formal \nbd{W^*}version of Definition \ref{Mepsdef} as the changes belong to our list of changes. But we would like to mention that Morita equivalence of correspondences (as introduced by Muhly and Solel \cite{MuSo00}) in the \nbd{W^*}case has a particularly nice interpretation in terms of our representation theory when applied to a unital normal endomorphism $\vt$ of $\sB^a(E)$. In fact, $_\vt\sB^a(E)$ is a \nbd{W^*}correspondence over $\sB^a(E)$ with left action via $\vt$. If $E$ is strongly full, then $E$ is a Morita \nbd{W^*}equivalence from $\sB^a(E)$ to $\cB$. The multiplicity correspondence of the endomorphism $\vt$ is nothing but $E_\vt:=E^*\sodots E=E^*\sodots{_\vt\sB^a(E)}\sodots E$, the ``conjugate'' of $_\vt\sB^a(E)$ with the Morita equivalence $E^*$. For a normal \nbd{E_0}semigroup $\vt$ on $\sB^a(E)$, we see that the product system of $\vt$ is simply conjugate to the one-dimensional product system of $\vt$ (see the methodological introduction) via the Morita equivalence $E^*$. In the proof of \cite[Theorem 5.12]{Ske09} and its corollary we used the statement that a strongly full product system of \nbd{W^*}correspondences is the product system of an \nbd{E_0}semigroup if and only if it arises in the described way by conjugation from a one-dimensional product system, that is, if and only if it is Morita equivalent to a one-dimensional product system.}%
~\it{Cum grano salis} for Section \ref{sccSEC} means that stably, of course, has to be replaced with the version where stably means \nbd{\en}stably for some infinite cardinal number $\en$. We only reformulate the main results of Sections \ref{E0MeSEC} and \ref{sccSEC}.

\bthm
Let $\vt$ and $\theta$ be normal \nbd{E_0}semigroups on $\sB^a(E)$ and $\sB^a(F)$, respectively, and suppose that $\sB^a(E)$ and $\sB^a(F)$ are isomorphic. Then $\vt$ and $\theta$ are unitary cocycle conjugate if and only if their associated product systems are Morita equivalent via the same Morita \nbd{W^*}equivalence inducing the isomorphism of $\sB^a(E)$ and $\sB^a(F)$.
\ethm

If $\en$ is a cardinal number and if $\vt$ is an \nbd{E_0}semigroup on $\sB^a(E)$, denote by $\vt^\en$ the \hl{amplification} of $\vt$ to $\sB^a(\ol{E^\en}^s)$.

\bthm\label{vNsclthm}
\begin{enumerate}
\item
Let $\vt$ be a normal \nbd{E_0}semigroup on $\sB^a(E)$ for a strongly full \nbd{W^*}mod\-ule $E$ over $\cB$. Then there exists a cardinal number $\en$ such that the amplification $\vt^\en$ is inner conjugate to an \nbd{E_0}semigroup $\vt^\cB$ on $\sB^a(\ol{\cB^\en}^s)$.

\item
Let $\vt$ and $\vt'$ be normal \nbd{E_0}semigroups on $\sB^a(E)$ and $\sB^a(E')$, respectively, where $E$ and $E'$ are strongly full \nbd{W^*}modules over $\cB$. Then the following conditions are equivalent:
\begin{enumerate}
\item
$\vt$ and $\vt'$ are stably unitary cocycle inner conjugate.

\item
There exists a cardinal number $\en$ such that $\vt^\cB$ and $\vt'^\cB$ are unitary cocycle equivalent.

\item
$\vt$ and $\vt'$ have isomorphic product systems.
\end{enumerate}

\item
Let $\vt$ and $\theta$ be normal \nbd{E_0}semigroups on $\sB^a(E)$ and $\sB^a(F)$, respectively, where $E$ and $F$ are strongly full \nbd{W^*}modules over $\cB$ and $\cC$, respectively. Then $\vt$ and $\theta$ are stably unitary cocycle (ternary) conjugate if and only if they have Morita (ternary) equivalent product systems.
\end{enumerate}
\ethm

\noindent
We explain briefly in which sense a part of the results of Section \ref{contSEC} are available already now in the \nbd{W^*}context. Essentially, we mean all results that are algebraic without continuity conditions. The reason why we can allow this, is that Proposition \ref{vNLanprop}, the analogue Proposition \ref{Lanprop}, now, does no longer depend on countability conditions. (Remember: The main reason, why in Section \ref{contSEC} we had to restrict to the continuous case, was precisely to guarantee these countability conditions.)

One of the main results of \cite{Ske09} asserts that every \it{discrete} product system $\bfam{E_n}_{n\in\N_0}$ of strongly full \nbd{W^*}correspondences is the product system of a discrete \nbd{E_0}semigroup. This completed the classification for the case of the discrete semigroup $\bS=\N_0$. But, for the continuous case $\bS=\R_+$ it also means that for every strongly full product system $E^\odot=\bfam{E_t}_{t\in\R_+}$ of \nbd{W^*}correspondences we can find a left dilation of the discrete subsystem $\bfam{E_t}_{t\in\N_0}$ to a strongly full \nbd{W^*}module $\breve{E}$. Such a left dilation of the discrete subsystem is precisely the main input for the construction in \cite{Ske06} of an \nbd{E_0}semigroup for every Arveson system. We mentioned already in \cite{Ske06} that the construction works without any problem if all the direct integrals are with respect to the counting measure. So, if we define $\int_a^bE_\alpha\,d\alpha:=\ol{\bigoplus_{\alpha\in\RO{a,b}}E_\alpha}^s$ and if we put $E:=\breve{E}\sodots\!\int_0^1E_\alpha\,d\alpha$, then the following formula from \cite{Ske06}
\bal{
\notag
E\sodots\! E_t
~=~
\breve{E}\sodots\!\family{\int_0^1E_\alpha\,d\alpha}\sodots\! E_t
&
~=~
\breve{E}\sodots\!\!\int_t^{1+t}E_\alpha\,d\alpha
\\[2ex]
\notag
&
~\cong~
\family{\breve{E}\sodots\! E_n\sodots\!\!\int_{t-n}^1E_\alpha\,d\alpha}
\oplus
\family{\breve{E}\sodots\! E_{n+1}\sodots\!\!\int_0^{t-n}E_\alpha\,d\alpha}
\\[2ex]\label{idea}
&
~\cong~
\family{\breve{E}\sodots\!\!\int_{t-n}^1E_\alpha\,d\alpha}
\oplus
\family{\breve{E}\sodots\!\!\int_0^{t-n}E_\alpha\,d\alpha}
~=~
E
}\eal
($n\in\N_0$ such that $t-n\in\RO{0,1}$) suggests an isomorphism $v_t\colon E\sodots E_t\rightarrow E$ for every $t\in\R_+$. By \cite[Proposition 3.1]{Ske06}, these maps $v_t$ form a left dilation of $E^\odot$ to $E$. Just that, by construction, the induced \nbd{E_0}semigroup is definitely not continuous with time in any reasonable topology. Nevertheless, we can formulate the classification theorem for the classification of algebraic normal \nbd{E_0}semigroups by algebraic product systems---a result that has no analogue in the \nbd{C^*}theory.

\bthm\label{vNalgclthm}
\begin{enumerate}
\item
Let $\cB$ denote a \nbd{W^*}algebra. Then there is a one-to-one correspondence between equivalence classes (up to stable unitary cocycle inner conjugacy) of normal \nbd{E_0}semigroups acting on the algebras of operators on strongly full \nbd{W^*}modules over $\cB$ and isomorphism classes of strongly full product systems of \nbd{W^*}correspondences over $\cB$.

\item
There is a one-to-one correspondence between equivalence classes (up to stable unitary cocycle conjugacy) of normal \nbd{E_0}semigroups acting on the algebras of operators on strongly full \nbd{W^*}modules and Morita equivalence classes of strongly full product systems of \nbd{W^*}correspondences.
\end{enumerate}
\ethm

%


\section{Von Neumann case: Topological classification}\label{vNcontSEC}

We now come to strongly continuous \nbd{E_0}semigroups in the \nbd{W^*}case%
\footnote{ \label{scFN}
Once for all, when we speak about \hl{strongly continuous} semigroups on a \bf{\nbd{W^*}algebra} $\cB$, what we have in mind is the point-strong topology of a subalgebra $\cB$ of some $\sB(G)$: A semigroup $T$ on $\cB$ is \hl{strongly continuous} if $t\mapsto T_t(b)g$ is continuous for all $b\in\cB,g\in G$. The strong topology depends on the identification $\cB\subset\sB(G)$, continuity results do not. Usually, continuity for semigroups on \nbd{W^*}algebras is formulated in terms of the weak$^*$ topology induced by the pre-dual. But, by result of Markiewicz and Shalit \cite{MaSha10}, this implies strong continuity for every faithful normal representation of $\cB$. We reprove this result in Corollary \ref{MScor}. (For \nbd{E_0}semigroups it is elementary; see Part (iv) of the proof of Theorem \ref{scdilthm}.)

By Elliot's result in \cite{Ell00}, a semigroup of normal CP-maps on a \nbd{W^*}algebra that is strongly (or weakly) continuous in the \nbd{C^*}sense, is uniformly continuous. So it is pointless to speak about them.
}%
, a property that has to be reflected by a property of the associated product system of \nbd{W^*}correspondences. This brings us to the problem that we have to give a concise definition of strongly continuous product system. Already in Skeide \cite{Ske03b} we indicated briefly how this can be done, following the procedure in the \nbd{C^*}case. This is what we will do here in order to be able to work. But we do not hide that we think this definition should be considered as a preliminary working definition. The reason is as follows: In the von Neumann way to see things, von Neumann correspondences come shipped with a commutant generalizing the commutant of a von Neumann algebra; see \cite{Ske03c,Ske06b}. The same is true for whole product systems; see \cite{Ske03b,Ske09,MuSo07}. In \cite{MuSo07}, Muhly and Solel introduced a weakly measurable version of product systems and showed (under separability assumptions, and using their independent way \cite{MuSo04} to look at the commutant) that also the commutant system has a measurable structure by reducing it to Effros' analogue result \cite{Eff65} for fields of von Neumann algebras. It is not difficult to see that the definition we will use here, is manifestly asymmetric under commutant. Our scope in \cite{Ske07p}, among others, will be to provide a notion of strongly continuous product system that is compatible with the commutant. Therefore, we would like to consider the definition used here as preliminary. (See, however, Theorem \ref{scpsdualthm}.)

Recall that the \nbd{W^*}algebra $\cB$ is assumed to act as a von Neumann algebra on a Hilbert space $G$. Following the suggestion in \cite{Ske03b}, we define as follows:

\bdefi\label{scPSdefi}
Suppose $E^\odot$ is a product system of \nbd{W^*}correspondences over $\cB$. Suppose further that $i_t$ is a family of isometric embeddings of $E_t$ into a fixed \nbd{W^*}module $\wh{E}$ over $\cB$, and denote by
\beqn{
CS_i^s(E^\odot)
~:=~
\BCB{\bfam{x_t}_{t\in\R_+}\colon x_t\in E_t, t\mapsto i_tx_t\odot g\in\wh{H}:=\wh{E}\odot G\text{~is continuous for all~}g\in G}
}\eeqn
the set of \hl{strongly continuous sections} (with respect to the embedding $i$). We say $E^\odot$ is a \hl{strongly continuous} product system if
\beqn{
\BCB{x_s\colon\bfam{x_t}_{t\in\R_+}\in CS_i^s(E^\odot)}
~=~
E_s
}\eeqn
for all $s\in\R_+$, and if the function
\beqn{
(s,t)
~\longmapsto~
i_{s+t}(x_sy_t)\odot g\in\wh{H}
}\eeqn
is continuous for all $\bfam{x_t}_{t\in\R_+},\bfam{y_t}_{t\in\R_+}\in CS_i^s(E^\odot)$ and for all $g\in G$.

A morphism between strongly continuous product systems is \hl{continuous} if it sends strong\-ly continuous sections to strongly continuous sections. By Theorem \ref{scisothm}, a continuous isomorphism has a continuous inverse.
\edefi

\brem
As said in the beginning of this section, we consider this a working definition. We do not know if the condition that each $x_s\in E_s$ is the value $x_s=y_s$ of a strongly continuous section $y$ can be weakened to strong totality of such $y_s$ in $E_s$. (In the \nbd{C^*}version this is so.) But the definition is justified by two facts: Firstly, every product system of a strongly continuous \nbd{E_0}semigroup admits such a strongly continuous structure. (This will be explained immediately.) Secondly, every (strongly full) strongly continuous product system can be obtained in that way from an \nbd{E_0}semigroup. (The proof of the latter fact is postponed to Appendix B.1.)
\erem

If $E$ is a \nbd{W^*}module over $\cB$, then we turn $\sB^a(E)$ into a von Neumann algebra by embedding it faithfully as $\sB^a(E)\odot\id_G$ into $\sB(H)$, where $H:=E\odot G$. Like in the \nbd{C^*}case, if we have a normal \nbd{E_0}semigroup $\vt$ acting on $\sB^a(E)$, then it is strongly continuous (with respect to the strong topology of $\sB(H)$) if and only if each amplification $\vt^\en$ to $\sB^a(\ol{E^\en}^s)$ is strongly continuous (with respect to the strong topology of $\sB(H^\en)$). If $E$ is strongly full, then by \cite[Lemma 4.2]{Ske09} $F:=\ol{E^\en}^s$ has a unit vector $\zeta$ as soon as the cardinal number $\en$ is big enough. Like in the \nbd{C^*}case, we may use that unit vector to construct embeddings $i_t\colon E_t\rightarrow F$. It is easy to show (similar to the \nbd{C^*}case) that these embeddings equip $E^\odot$ with a strongly continuous structure, and that this strongly continuous structure does not depend neither on the choice of cardinal number $\en$ nor on the choice of $\zeta$. In particular, if $E$ has already a unit vector $\xi$, then $\en=1$ is among the admissible cardinal numbers and the strongly continuous structure derived from that $\xi$ coincides with all others.

Once more, if $u$ is a unitary cocycle for $\vt$, then it is strongly continuous (in $\sB(H)$) if and only if the automorphism of the associated product system $E^\odot$ corresponding to $u$ by Theorem \ref{ucethm}, is strongly continuous. (The proof is exactly like that of Theorem \ref{ucecthm} for the \nbd{C^*}case, except that we have to tensor everywhere in the estimates with an element $g\in G$.) We find:

\bthm\label{sucecthm}
Let $E$ be a strongly full \nbd{W^*}module and suppose $\vt$ and $\vt'$ are two strongly continuous normal \nbd{E_0}semigroups on $\sB^a(E)$. Then the following are equivalent:
\begin{enumerate}
\item
$\vt$ and $\vt'$ are unitary cocycle equivalent via a strongly continuous cocycle.

\item
The strongly continuous product systems associated with $\vt$ and $\vt'$ are isomorphic.
\end{enumerate}
\ethm

\noindent
This is the classification of strongly continuous normal \nbd{E_0}semigroups acting all on the same $\sB^a(E)$. Of course, also Theorem \ref{vNsclthm} remains true if we simply add everywhere strongly continuous, since (as pointed out before) amplification is compatible with strong continuity. To save space, we do not repeat it.

What is missing to obtain the strongly continuous analogue also for Theorem \ref{vNalgclthm}, is the following existence result:

\bthm\label{vNexithm}
Every strongly full strongly continuous product system of \nbd{W^*}correspondences is isomorphic to the strongly continuous product system associated with a normal strongly continuous \nbd{E_0}semigroup acting on the algebra of all adjointable operators of a strongly full \nbd{W^*}module.
\ethm

We prove this theorem in Appendix B.1. It immediately follows that the following strongly continuous version of Theorem \ref{vNalgclthm} holds, too, in all of its parts. (Two strongly continuous product systems $E^\odot$ and $F^\odot$ are in the same strongly continuous Morita equivalence class if the isomorphism between $F^\odot$ and the product system $M^*\sodots E^\odot\sodots M$ Morita equivalent to $E^\odot$ can be chosen continuous. See also Remark \ref{MecPSrem} for how the strongly continuous structure of $M^*\sodots E^\odot\sodots M$ can be defined.)

\bthm\label{vNalgsclthm}
\begin{enumerate}
\item
Let $\cB$ denote a \nbd{W^*}algebra. Then there is a one-to-one correspondence between equivalence classes (up to stable unitary cocycle inner conjugacy by a strongly continuous cocylce) of normal strongly continuous \nbd{E_0}semigroups acting on the algebras of operators on strongly full \nbd{W^*}modules over $\cB$ and continuous isomorphism classes of strongly full strongly continuous product systems of \nbd{W^*}correspondences over $\cB$.

\item
There is a one-to-one correspondence between equivalence classes (up to stable unitary cocycle conjugacy by a strongly continuous cocycle) of normal strongly continuous \nbd{E_0}semigroups acting the algebras of operators on strongly full \nbd{W^*}modules and strongly continuous Morita equivalence classes of strongly full product systems of \nbd{W^*}cor\-re\-spond\-ences.

\end{enumerate}
\ethm


\section{Von Neumann case: Spatial Markov semigroups}\label{vNHPSEC}

The discussion of spatial Markov semigroups on a \nbd{W^*}algebra $\cB$ and their Hudson-Par\-tha\-sa\-rathy dilations, apart from the weaker topologies, is very similar to the \nbd{C^*}case. (The weaker topologies sometimes  require different proofs, which we discuss in Appendix A.2.) We even have the simplification that spatial Markov semigroups, here, turn out to have a spatial product system; see Theorem \ref{vNspthm}. (In the \nbd{C^*}case, we had only embedding into a spatial product system.) Thanks to the weaker topology, the results are applicable to a much wider (thus, more interesting) class of Markov semigroups. (In fact, in the case $\cB=\sB(G)$ we do not know examples of nonspatial Markov semigroups, except for Markov semigroups that arise as a tensor product of a Markov semigroup with a nonspatial \nbd{E_0}semigroup. Only recently, in Skeide \cite[Corollary]{Ske12} we could prove that Floricel \cite{Flo08} provides us with a proper type III Markov semigroup, provided we find a type III Arveson system that does not factor into a tensor product. Existence of such an Arveson sytem is likely, but unproved yet. On the other hand, we know from Fagnola, Liebscher, and Skeide \cite{FLS05p,Ske05b} that the Brownian and Ornstein-Uhlenbeck semigroups acting on the commutative von Neumann algebra $\cB=L^\infty(\R)\subset\sB(L^2(\R))$ have nonspatial product systems. In fact, their product systems do not contain a single nonzero element that would commute with the algebra.)

We would like to mention that the discussion of the case $\sB(G)$ in Skeide \cite{Ske08p3}, actually, was inspired by the preparation of Section \ref{HuPaSEC} and the present section. But, while in \cite{Ske08p3} we used mainly well-known results about spatial \nbd{E_0}semigroups and spatial Arveson systems (that is, formulated with measurability conditions rather than continuity conditions), here we present a completely new treatment adapted to our notions of strong continuity.

\lf
While in Section \ref{HuPaSEC} we could build on the results on \it{continuous} units in \it{continuous} product systems from \cite{Ske03b}, in this section we have to develop the \it{strongly continuous} analogues. On the one hand, there are results, like the following analogue of Theorem \ref{indcontthm} (notation is introduced in front of Observation \ref{primeob}), that can be proved simply by tensoring the vectors in $\wh{E}$ occurring in the estimates in \cite{Ske03b}, with a fixed vector $g$ in the representation space $G$ of $\cB$.

\bthm\label{vNscsecthm}
Let $E^\odot$ be a strongly continuous (with respect to embeddings $i_t\colon E_t\rightarrow\wh{E}$, say) product system of \nbd{W^*}correspondences over $\cB$ and suppose that $\zeta^\odot\in CS_i^s(E^\odot)$ is unital unit among the strongly continuous sections.

Then the normal \nbd{E_0}semigroup $\vt$ on $\sB^a(\ol{E^\zeta}^s)$ is strongly continuous, and the strongly continuous structure induced on $E^\odot$ by $\vt^\zeta$ via the unit vector $\zeta$ coincides with the original one. In particular, the induced strongly continuous structure does not depend on the choice of the strongly continuous unital unit $\zeta^\odot$.
\ethm

\proof
Except for the modifications stated in front of the theorem, the proof goes exactly like the corresponding proof of \cite[Theorem 7.5]{Ske03b}, including the \it{add-on} mentioned in Footnote \ref{xiE0FN} to Theorem  \ref{indcontthm}.\qed

\lf
On the other hand, there are results like the construction of a strongly continuous weak dilation from a product system $E^\odot$ (so far, without a strongly continuous structure) and a unital strongly continuous unit $\xi^\odot$ that \hl{strongly generates} $E^\odot$ (that is, the smallest product subsystem of \nbd{W^*}cor\-re\-spond\-ences containing $\xi^\odot$ is $E^\odot$ itself).  Like in Section \ref{HuPaSEC}, from product system and strongly generating unit we construct the triple $(\ol{E^\xi}^s,\vt^\xi,\xi)$ which is a weak dilation, the so-called unique \hl{minimal dilation}, of the strongly continuous normal Markov semigroup $T^\xi:=\AB{\xi_t,\bullet\xi_t}$ on $\cB$. By \cite[Theorem 12.1]{BhSk00}, every strongly continuous normal Markov semigroup $T$ on a \nbd{W^*}algebra arises in that way from its \hl{GNS-system} $E^\odot$ and a strongly continuous unit $\xi^\odot$. But \cite[Theorem 12.1]{BhSk00} also claims that the minimal dilation is strongly continuous. Since the proof in \cite{BhSk00} had a gap, we give a complete proof of a slightly more general statement, Theorem \ref{scdilthm}, in Appendix A.2. The strongly continuous \nbd{E_0}semigroup of the minimal dilation of $T$ equips the GNS-system, $E^\odot$, with a strongly continuous structure and, of course, the unit $\xi^\odot$ (being the constant element $\xi\in E^\xi$) is among the strongly continuous sections.

If a (necessarily strongly continuous and normal) Markov semigroup $T$ admits a strongly continuous normal weak Hudson-Parthasarathy dilation $(E,\sS^u,\om)$ (by \hl{strongly continuous} we mean that also the cocycle $u$ in $\sB^a(E)$ is strongly continuous in $\sB(H)\supset\sB^a(E)$), then the product system of the dilation contains two unital units among its strongly continuous sections: One is the unit $\xi^\odot$ that gives back $T$ as $T^\xi$, because $\sS^u$ is a weak dilation, and the other is the unital reference unit $\om^\odot$ of the noise $\sS$. (As in Section \ref{HuPaSEC}, it suffices just to observe that the strongly continuous product systems of $\sS^u$ and of $\sS$ are continuously isomorphic. The former contains $\xi^\odot$ and the latter contains $\om^\odot$ as strongly continuous sections. So a section that is strongly continuous for one product system has an image in the other that is also strongly continuous.) We find the analogue of Proposition \ref{spembprop}: There is a \hl{spatial} strongly continuous product system (that is, a strongly continuous product system that is spatial with a strongly continuous reference unit $\om^\odot$) with a strongly continuous unit $\xi^\odot$ such that $T=T^\xi$.

If we define the semigroup $c_t:=\AB{\om_t,\xi_t}$, then $T$ is spatial in the sense of the following definition:

\bdefi
A \hl{unit} for a strongly continuous normal CP-semigroup $T$ on a \nbd{W^*}algebra $\cB$ is a strongly continuous semigroup $c=\bfam{c_t}_{t\in\bS}$ of elements in $\cB$ such that $T_t$ dominates the CP-map $b\mapsto c_t^*bc_t$ for all $t\in\bS$. We say $T$ is \hl{spatial} if it admits units.
\edefi

For the backwards direction we are done as soon as we are able to find for every spatial Markov semigroup a strongly continuous spatial product system (with reference unit $\om^\odot$, say) and a strongly continuous unit $\xi^\odot$ such that $T=T^\xi$. This is done in Appendix  A.2 in a way that is much more satisfactory than the \nbd{C^*}case. In fact, Theorem \ref{vNspthm} asserts that a Markov semigroup on a \nbd{W^*}algebra is spatial if and only if its GNS-system is spatial (including all requirements about strong continuity).

Once we have these ingredients, the construction of a Hudson-Parthasarathy dilation goes exactly as in Section \ref{HuPaSEC} on the page preceding Theorem \ref{C*HPthm}. Just that now there are no countability assumptions. The price to be paid is that now we have to choose our amplifications \it{big enough} when establishing stable cocycle conjugacy---a small price, of course. We do not give more details on these steps because they, really, are completely analogous to Section \ref{HuPaSEC}, and all necessary compatibility results regarding the strong topologies have been mentioned.

\bthm\label{vNHPthm}
Let $\cB$ be a \nbd{W^*}algebra and let $T$ be a strongly continuous normal Markov semigroup on $\cB$. Then $T$ admits a strongly continuous normal weak Hudson-Parthasarathy dilation if and only if $T$ is spatial. Like in Theorem \ref{inextthm}, this dilation may be obtained as the restriction of a strongly continuous inner normal Hudson-Parthasarathy dilation.
\ethm 

\brem
It is standard to show that the construction preserves countability assumptions: If $\cB$ has separable pre-dual, then both the minimal dilation and the Hudson-Parthasarathy dilation of $T$ act on a $\sB^a(E)$ with separable pre-dual.
\erem

\bob
In cases where the space of the minimal dilation $E^\xi$ and the space of the noise $E^\om$ are isomorphic from the beginning (for instance, when $\cB=\sB(G)$ with separable $G$), amplification is not necessary. In this case, we see that the Hudson-Parthasarathy dilation interpreted as weak dilation, is the minimal dilation.
\eob


\section*{Appendix A: Strong type I product systems}
\addcontentsline{toc}{section}{Appendix A: Strong type I product systems}
	\stepcounter{section}
        \renewcommand{\thesection}{A}

A \hl{strong type I system} is a product system $E^\odot$ that is generated by a unital strongly continuous set of units $S$. \hl{Unital} means that $S$ contains at least one unital unit $\xi^\odot$. The topology in which the units generate the product system depends on whether we are speaking about unital \nbd{C^*}algebras and modules or about \nbd{W^*}algebras and modules.%
\footnote{ \label{genFN}
In general, by the product subsystem (of pre-correspondences, of correspondences, of \nbd{W^*}correspondences) of a product system (of pre-correspondences, of correspondences, of \nbd{W^*}correspondences) $E^\odot$ \hl{generated} by a certain set of elements $x_t\in E_t,t\in\R_+$ we mean, as usual, the smallest product subsystem (of pre-correspondences, of correspondences, of \nbd{W^*}correspondences) containing all these elements. (By ``dual linear algebra'', for vector suspaces $V_i\subset V$, $W_i\subset W$, we have $(V_1\otimes W_1)\cap(V_2\otimes W_2)=(V_1\cap V_2)\otimes(W_1\cap W_2)$. This allows to conclude that, in each case, a minimal product subsystem exists. If the generating elements come from a set $S$ of units, then the subsystems generated by $S$ my be identified explicitly by formulae like \eqref{unigen}.)
}
Also strongly continuous means two different things: In the \nbd{C^*}case we mean that a family of maps $T_t\colon\cA\rightarrow\cB$ is strongly continuous if $t\mapsto T_t(a)$ is continuous in $\cB$ for all $a\in\cA$. The occurring product systems will be continuous product systems. In the \nbd{W^*}case we refer to the strong operator topology of $\cB\subset\sB(G)$, so, $t\mapsto T_t(a)g$ is continuous in $G$ for all $a\in\cA$ and $g\in G$. (See also Footnote \ref{scFN} in the beginning of Section \ref{vNcontSEC}.) The occurring product systems will be strongly continuous product systems. 

\lf
This appendix serves two purposes.

The first purpose is to show that every strong type I system comes along with a unique continuous structure making $S$ a set of (strongly) continuous sections. (Theorem \ref{cdilthm} in the \nbd{C^*}case and Theorem \ref{scdilthm} in the \nbd{W^*}case.) The idea is to use the unital unit $\xi^\odot$ to construct an \nbd{E_0}semigroup $\vt$ on $\sB^a(E)$ that is strongly continuous and whose associated product system is $E^\odot$. This equips $E^\odot$ with a (strongly) continuous structure. (Note that \cite[Theorem 7.5]{Ske03b} (Theorem \ref{vNscsecthm}) starts from a product system that has already a (strongly) continuous structure and from the assumption that $\xi^\odot$ is among the (strongly) continuous sections.) By \cite[Theorem 7.5]{Ske03b} (Theorem \ref{vNscsecthm}), such a structure, if it exists, is unique. In particular, once we show that $S$ is a subset of the (strongly) continuous sections of this structure, we know that the structure does not depend on the choice of the unital unit $\xi^\odot\in S$.

We thank Orr Shalit for having pointed out that the proofs of strong continuity of $\vt$ in \cite[Theorems 10.2, 12.1]{BhSk00} contain a gap. In fact, in either case only right continuity is shown; but in either case the situation is not among those where right continuity would imply also left continuity. (Strong or weak operator topology of both $\sB^a(E)$ (\nbd{C^*}case) and $\sB(H)\supset\sB^a(E)$ (\nbd{W^*}case) are too weak to achieve this immediately from standard theorems! This error is surprisingly frequent in literature and led Markiewicz and Shalit to prove in \cite{MaSha10}, once for all, that weak operator continuity of a normal CP-semigroup on a von Neumann algebra implies strong operator continuity.) In this appendix we fix this gap in a systematic way for the most general situation. We also reprove the result of \cite{MaSha10} in a different way as Corollary \ref{MScor}.

The second purpose is to provide, in the \nbd{W^*}case, some (new) basic classification results for strong type I systems. (Corollary \ref{wstIcor}, Corollary \ref{elemtrivcor}, and Theorem \ref{vNspthm}; also some more technical results about strongly continuous sets of units, which have nice interpretations in terms of CPD-semigroups. \it{CPD-semigroup} is central a notion for the classification of \it{type I} product systems from \cite{BBLS04} which we avoid to discuss here for reasons of space, but which would be part of a systematic study of strong type I systems.) This part is limited to the \nbd{W^*}case, because we know from Bhat, Liebscher, and Skeide \cite{BLS10} that these results fail in the \nbd{C^*}case.

Our reference for results about one-parameter semigroups is Engel and Nagel \cite{EnNa06}; see Footnote \ref{sgFN}.

\subsection*{A.1: The $C^*$--case}

In the discussion culminating in Theorem \ref{indcontthm}, we have described a construction that, starting from a product system $E^\odot$ and a unital unit $\xi^\odot$, results in a weak dilation $(E^\xi,\vt^\xi,\xi)$ of the Markov semigroup $T^\xi$ defined by $T^\xi_t=\AB{\xi_t,\bullet\xi_t}$. By \cite[Theorem 7.5]{Ske03b}, if $E^\odot$ is continuous and if $\xi^\odot$ is among its continuous sections, then $\vt^\xi$ is strongly continuous and the continuous structure on $E^\odot$ derived from $\vt^\xi$ coincides with the original one. In particular, the continuous structure induced by $\vt^\xi$ does not depend on the choice of $\xi^\odot$. As pointed out in the introduction to this appendix, we will start with a product system that is generated (see Footnote \ref{genFN}) by a unital strongly continuous set of units $S\ni\xi^\odot$ and that has \bf{not yet} a continuous structure. Instead, it will be our scope to show that under these assumptions $\vt^\xi$ is strongly continuous. This $\vt^\xi$ can, then, be used to induce a continuous structure on $E^\odot$.

\bthm\label{cdilthm}
Suppose $S$ is a set of units for the product system $E^\odot$ that generates $E^\odot$, and suppose $\xi^\odot\in S$ is unital. If for all $\zeta^\odot,\zeta'^\odot\in S$ and $b\in\cB$ the map $t\mapsto\AB{\zeta_t,b\zeta'_t}$ is continuous, then the \nbd{E_0}semigroup $\vt^\xi$ is strongly continuous. Moreover, the continuous structure induced on $E^\odot$ by $\vt^\xi$ is the unique one making all $\zeta^\odot\in S$ continuous sections.
\ethm

\proof
We explained already that if $\vt^\xi$ is strongly continuous, then by \cite[Theorem 7.5]{Ske03b} the continuous structure derived from it is the unique one, making $\xi^\odot$ a continuous section. So it remains to show that $\vt^\xi$ is strongly continuous, and that $S$ is a subset of the set of continuous sections of this structure.

It is easy to see that the product subsystem ${E^S}^\odot$ generated by $S$ is formed by the spaces
\beq{ \label{unigen}
E^S_t
~:=~
\cls\BCB{~b_n\xi^n_{t_n}\ldots b_1\xi^1_{t_1}b_0\colon \,n\in\N,\,b_i\in\cB,\,{\xi^i}^\odot\!\!\in S,\,t_i>0,\,t_1+\ldots+t_n=t~}
}\eeq
for $t>0$ and, of course, $E^S_0=E_0=\cB$. Since $S$ is assumed generating, we have $E^S_t=E_t$. The inductive limit $E^\xi$ is, therefore, spanned by elements of the form
\beqn{
x
~=~
\xi b_n\xi^n_{t_n}\ldots b_1\xi^1_{t_1}b_0.
}\eeqn
(No condition on $t_1+\ldots+t_n$, here!) As in \cite{BhSk00} one easily shows that the semigroup of time-shift operators $\s_t\colon y\mapsto y\xi_t$ on $E^\xi\ni y$ is strongly continuous at $t=0$. (Since all $\xi_t$ are contractions, it is sufficient to show this on the total set of elements of the form $x$. Inserting the the concrete form of $x$ and calculating $\AB{x\xi_\ve-x,x\xi_\ve-x}$, one sees that this goes to $0$ for $\ve\to0$. Since that part of the proof in \cite{BhSk00} is okay, and because we discuss in all detail the strong version in the \nbd{W^*}case in Appendix A.2, we do not give more details, here.) It follows that the semigroup $\s_t$ is strongly right-continuous everywhere. (By \cite[Proposition I.1.3]{EnNa06}, it is even strongly continuous, but we do not need this fact.) Like in \cite{BhSk00}, we compute
\beqn{
\vt_t(a)x-ax
~=~
\vt_t(a)x-\vt_t(a)x\xi_t+\vt_t(a)x\xi_t-ax
~=~
\vt_t(a)(x-x\xi_t)+((ax)\xi_t-(ax)),
}\eeqn
and conclude that $\vt_t(a)$ is strongly right-continuous. (So far, this is the old proof.) Additionally, we observe that nothing changes in the proof if we replace the semigroup $\s_t$ with with the semigroup $\s^\zeta_t\colon x\mapsto x\zeta_t$ for any other unit $\zeta^\odot$ in $S$.

To add left continuity, fix $t>0$ and fix an element of the form $x$. Without changing the value of $x$, we may assume that: Firstly, $t_1+\ldots+t_n>t$. (Otherwise, choose $t_{n+1}$ sufficiently big, and $b_{n+1}=\U$.) Secondly, there is $m<n$ such that $t_1+\ldots+t_m=t$. (Otherwise, split that $\xi^m_{t_m}$ where $t_1+\ldots+t_{m-1}<t<t_1+\ldots+t_m$ into a product $\xi^m_{t_m-s}\xi^m_s$ such that $t_1+\ldots+t_{m-1}+s=t$ and insert a $\U$ in between.) Choose $a\in\sB^a(E^\xi)$ and $0<\ve<\min\CB{t_1,\ldots,t_n}$ (so that $\ve\le t$, too). Then
\bmun{
\vt_t(a)x
~=~
(a(\xi b_n\xi^n_{t_n}\ldots b_{m+1}\xi^{m+1}_{t_{m+1}}b_m))\xi^m_{t_m}\ldots b_1\xi^1_{t_1}b_0
\\
~=~
\bfam{a(\xi b_n\xi^n_{t_n}\ldots b_{m+1}\xi^{m+1}_{t_{m+1}}b_m)}\xi^m_\ve\xi^m_{t_m-\ve}\ldots b_1\xi^1_{t_1}b_0
}\emun
and
\beqn{
\vt_{t-\ve}(a)x
~=~
\bfam{a(\xi b_n\xi^n_{t_n}\ldots b_{m+1}\xi^{m+1}_{t_{m+1}}b_m\xi_\ve)}\xi^m_{t_m-\ve}\ldots b_1\xi^1_{t_1}b_0.
}\eeqn
We observe that $\xi^m_{t_m-\ve}\ldots b_1\xi^1_{t_1}b_0$ is bounded uniformly in $\ve$, say, by $M$.
If we abbreviate $y:=\xi b_n\xi^n_{t_n}\ldots b_{m+1}\xi^{m+1}_{t_{m+1}}b_m$, we find
\bmun{
\norm{\vt_{t-\ve}(a)x-\vt_t(a)x}
\\[1ex]
~\le~
M
\norm{a(\xi b_n\xi^n_{t_n}\ldots b_{m+1}\xi^{m+1}_{t_{m+1}}b_m\xi^m_\ve)
-
\bfam{a(\xi b_n\xi^n_{t_n}\ldots b_{m+1}\xi^{m+1}_{t_{m+1}}b_m)}\xi^m_\ve}
~=~
M\snorm{a(\s^{\xi^m}_\ve y)
-
\s^{\xi^m}_\ve(ay)},
}\emun
which goes to $0$ for $\ve\to0$.

We conclude that $\vt^\xi$ is strongly right-continuous and strongly left-continuous and, therefore, strongly continuous.

To show that $\zeta^\odot\in S$ is among the continuous sections, we have to compute the norm of $\xi\zeta_t-\xi\zeta_s$. For $t\ge s$ this is  $\xi\zeta_{t-s}\zeta_s-\xi\zeta_s=(\s^\zeta_{t-s}\xi-\xi)\zeta_s$. Since $\zeta_s$ is bounded uniformly on any compact interval, continuity of $\zeta^\odot$ follows.\qed

\subsection*{A.2: The $W^*$--case}

The preceding part of the appendix fixes only a gap in the proof of  \cite[Theorem 10.2]{BhSk00}, and also proves the case where $S$ is bigger than $\CB{\xi^\odot}$. (The product system of \cite[Theorem 10.2]{BhSk00} is generated by a single strongly continuous unit.) The present part, apart from fixing an analogue gap in the proof of \cite[Theorem 12.1]{BhSk00}, may also be considered as a start-up for the theory of strongly continuous product systems of \nbd{W^*}correspondences that are generated by their strongly continuous units (\it{strong type I}). It is not exhaustive, and derives only those results that we need for Section \ref{vNHPSEC}.

The following lemma is the generalization of Accardi and Mohari \cite[Lemma 3.2]{AcMo96} to \nbd{W^*}algebras $\cB$ with not necessarily separable pre-dual. (The word ``net'' in their proof, should be replaced by ``sequence''. And even then, it seems that a restriction to bounded subsets is still necessary. A similar result, which deals only with convergence but not with continuity, is \cite[Lemma 4.1(b)]{MuSo02}.) By the \nbd{\sigma}weak topology on $\R_+\times\cB$ we mean the product topology of the usual topology on $\R_+$ and the \nbd{\sigma}weak topology on $\cB$.

\blem\label{s*contlem}
Let $T$ be a \nbd{\sigma}weakly continuous one-parameter semigroup on a \nbd{W^*}algebra $\cB$, and fix an arbitrary bounded subset $B$ of $\cB$. Then the map $(t,b)\mapsto T_t(b)$ is a continuous map $\R_+\times B\rightarrow\cB$ for the (relative) \nbd{\sigma}weak topologies on either side.
\elem

\proof
$T$ being \nbd{\sigma}weakly continuous, means that the pre-dual semigroup $T_*$ on the pre-dual $\cB_*$ of $\cB$ is weakly, hence, by \cite[Theorem I.1.6]{EnNa06}, strongly continuous. Therefore, by \cite[Proposition I.1.4]{EnNa06}, $T$, like $T_*$, is bounded by a family of constants $\bfam{Me^{t\gamma}}_{t\in\R_+}$ for suitable positive numbers $M,\gamma$. We shall assume that $\gamma=0$ passing, if necessary, to the semigroup $T$ rescaled by $e^{-t\gamma}$. Denote by $(L_*,\sD(L_*))$ the generator of $T_*$, and choose an element $\vp\in\sD(L_*)$. Then, by \cite[Equation (1.6) of Lemma II.1.3]{EnNa06},
\beqn{
(T_*)_{t'}(\vp)
~=~
(T_*)_t(\vp)+\int_t^{t'} (T_*)_s(L_*(\vp))\,ds
}\eeqn
(in norm of $\cB_*$). Let $\bfam{(t_\lambda,b_\lambda)}_{\lambda\in\Lambda}$ be a net converging \nbd{\sigma}weakly in $\R_+\times B$ to $(t,b)$, that is, $t_\lambda\to t$ and $\vp(b_\lambda)\to\vp(b)$ for every $\vp\in\cB_*$. For $\vp\in\sD(L_*)$ we find
\bmun{
\abs{\vp(T_{t_\lambda}(b_\lambda)-T_t(b))}
~=~
\abs{\SB{(T_*)_{t_\lambda}(\vp)}(b_\lambda)-\SB{(T_*)_t(\vp)}(b)}
\\
~=~
\Babs{\SB{(T_*)_t(\vp)}(b_\lambda-b)+\int_t^{t_\lambda}\SB{(T_*)_s(L_*(\vp))}(b_\lambda)\,ds}
\\
~\le~
\babs{\SB{(T_*)_t(\vp)}(b_\lambda-b)}+M\abs{t_\lambda-t}\,\babs{\SB{L_*(\vp)}(b_\lambda)}.
}\emun
The first summand converges to $0$. The second summand converges to $M\cdot 0\cdot\SB{L_*(\vp)}(b)=0$, too. Now, since an arbitrary $\vp\in\cB_*$ is the norm limit of elements in $\sD(L_*)$ and since the $b_\lambda\in B$ are bounded uniformly in $\lambda$ , it follows that $\sabs{\vp(T_{t_\lambda}(b_\lambda)-T_t(b))}\to0$ for all $\vp\in\cB_*$.\qed

\lf
Since the functionals $\AB{g_1,\bullet g_2}$ for $g_1,g_2$ from a total subset of $G$, form a total subset of the pre-dual of the von Neumann algebra $\cB\subset\sB(G)$, for checking \nbd{\sigma}weak continuity of locally bounded semigroups (like \nbd{E_0}semigroup and other contraction semigroups) it is enough to check weak continuity with functionals from that total subset. To check \nbd{\sigma}weak continuity of an arbitrary semigroup it is enough to check weak continuity, that is, to check continuity with functionals $\AB{g_1,\bullet g_2}$ for \bf{all} $g_1,g_2\in G$. (This follows by a twofold application of the \it{principle of uniform boundedness}.)

We are now ready to prove a generalization of \cite[Theorem 12.1]{BhSk00}, fixing also the gap in the proof of that theorem.

\bdefi
A product system of \nbd{W^*}correspondences is \hl{strong type I} if it is generated by a \hl{strongly continuous} set $S$ of units, that is, for each $\xi^\odot,\xi'^\odot\in S$ the semigroup $\bfam{\AB{\xi_t,\bullet\xi'_t}}_{t\in\R_+}$ is strongly continuous. We use similar definitions for \hl{weak type I}. A strongly continuous product system is \hl{strong/weak type I} if the generating strongly/weakly continuous set of units can be chosen from the strongly continuous sections.
\edefi

By the discussion preceding the definition, a weakly continuous set of units is also \nbd{\sigma}weakly continuous. Therefore, we dispense with defining \it{\nbd{\sigma}weak type I}.

\bthm\label{scdilthm}
Let $E^\odot$ be a product system of \nbd{W^*}correspondences over a \nbd{W^*}algebra $\cB$ of weak type I with generating weakly continuous set $S$ of units. Furthermore, suppose that $\xi^\odot\in S$ is a unital unit. Denote by $(\ol{E^\xi}^s,\vt^\xi,\xi)$ the strong closure of $E^\xi$, the normal extension of the \nbd{E_0}semigroup on $\sB^a(E^\xi)$ to $\sB^a(\ol{E^\xi}^s)$ and the unit vector $\xi$ as constructed for Theorem \ref{indcontthm}. Then $\vt^\xi$ is strongly continuous.
\ethm

\proof
Recall that the elements of the form
\beq{\label{form}
\xi b_n\xi^n_{t_n}\ldots b_1\xi^1_{t_1}b_0
~~~~~~(n\in\N,t_i>0,{\xi^i}^\odot\in S,b_i\in\cB)
}\eeq
are total in $E^\xi$. Thus, they are strongly total in $\ol{E^\xi}^s$.

\item
Let $x=\xi b_n\xi^n_{t_n}\ldots b_1\xi^1_{t_1}b_0$. We will show that the map
\beqn{
t
~\longmapsto~
x\xi_t
}\eeqn
is strongly continuous, that is, that $t\mapsto x\xi_t\odot g$ is continuous for $g$ in the representation space $G$ of $\cB$. Put ${\xi^0}^\odot:=\xi^\odot$ and define the \nbd{\sigma}weakly continuous semigroups $T^{i,j}$ and $T^i$ $(0\le i,j\le n)$ by setting $T^{i,j}_t:=\AB{\xi^i_t,\bullet\xi^j_t}$ and $T^i:=T^{i,i}$. By Lemma \ref{s*contlem},
\beqn{
t
~\longmapsto~
\norm{x\xi_t\odot g}^2
~=~
\bAB{g,T^0_t(b_0^*T^1_{t_1}(b_1^*\ldots T^n_{t_n}(b_n^*b_n)\ldots b_1)b_0)g}
}\eeqn
is continuous. It remains to show that 
\beqn{
t
~\longmapsto~
\AB{x\xi_s\odot g,x\xi_t\odot g}
}\eeqn
depends continuously on $t$ in a neighbourhood of $s$. To see that this is so, we observe that in either case, $t\ge s$ and $t\le s$, we find $m\in\N;s_j>0;0\le i_j,k_j\le n;c_j,c'_j\in\cB$ such that the right-hand side becomes
\beqn{
\bAB{g,T^{i_1,k_1}_{s_1}(c_1^*\ldots T^{i_m,k_m}_{s_m}(c_m^*c'_m)\ldots c'_1)g}
}\eeqn
with $s_j$ depending jointly continuously on $t$ and $s$. (Simply factor in $x\xi_s$ and in $x\xi_t$ the pieces $\xi^i_{t_i}$ of the units into products of smaller pieces, so that the involved time points in both coincide and the inner product can be calculated.) By induction, Lemma \ref{s*contlem} tells us that this depends continuously on $t$ in either case.

\item
Let $x,y,z\in E^\xi$ have the form in \eqref{form}, put $a:=yz^*$, so that $ax=y\AB{z,x}$ also has the form in \eqref{form}, and choose $g\in G$. Recall that $\vt_t(a)x\xi_t=(ax)\xi_t$. Like in \cite{BhSk00}, we compute
\bmun{
\vt_t(a)x\odot g-ax\odot g
~=~
\vt_t(a)x\odot g-\vt_t(a)x\xi_t\odot g+\vt_t(a)x\xi_t\odot g-ax\odot g
\\
~=~
\vt_t(a)(x\odot g-x\xi_t\odot g)+((ax)\xi_t\odot g-(ax)\odot g).
}\emun
Since $x$ and $ax$ have the form in \eqref{form}, by (i) this converges to $0$ for $t\to0$. By boundedness of $t\mapsto\vt_t(a)$, this shows that $\vt_t(a)$ is strongly continuous at $0$ at least for all $a$ of the given form.

\item
Let $x_1,y_1,z_1,x_2,y_2,z_2\in E^\xi$ have the form in \eqref{form}, and choose $g_1,g_2\in G$. Fix an arbitrary $a\in\sB^a(\ol{E^\xi}^s)$. Observe that $y_1x_1^*ax_2y_2^*=y_1\AB{x_1,ax_2}y_2^*$ is an operator of the form dealt with in (ii). Observe also that the elements $xy^*z\odot g:=x\AB{y,z}\odot g$ still form a total subset. Also here $xy^*z\xi_t\odot g=\vt_t(xy^*)z\xi_t\odot g$. We compute
\baln{
\AB{x_1y_1^*z_1&\odot g_1,\vt_t(a)x_2y_2^*z_2\odot g_2}-\AB{x_1y_1^*z_1\odot g_1,ax_2y_2^*z_2\odot g_2}
\\
~=~
&
\bAB{x_1y_1^*z_1\odot g_1,\vt_t(a)x_2y_2^*z_2\odot g_2}-\bAB{x_1y_1^*z_1\xi_t\odot g_1,\vt_t(a)x_2y_2^*z_2\xi_t\odot g_2}
\\
&
~~~+\bAB{\vt_t(x_1y_1^*)z_1\xi_t\odot g_1,\vt_t(a)\vt_t(x_2y_2^*)z_2\xi_t\odot g_2}-\bAB{x_1y_1^*z_1\odot g_1,ax_2y_2^*z_2\odot g_2}
\\
~=~
&
\bAB{x_1y_1^*z_1v\odot g_1-x_1y_1^*z_1\xi_t\odot g_1\,,\,\vt_t(a)x_2y_2^*z_2\odot g_2}
\\
&
~~~+~~~~~~~~~~~~~~~~~\bAB{x_1y_1^*z_1\xi_t\odot g_1\,,\,\vt_t(a)(x_2y_2^*z_2\odot g_2-x_2y_2^*z_2\xi_t\odot g_2)}
\\
&
~~~+~~~~~~~~~\bAB{z_1\xi_t\odot g_1-z_1\odot g_1\,,\,(y_1x_1^*ax_2y_2^*)z_2\xi_t\odot g_2}
\\
&
~~~+~~~~~~~~~~~~~~~~~~~~~~~~~~~\bAB{z_1\odot g_1\,,\,(yx^*ax_2y_2^*)(z_2\xi_t\odot g_2-z_2\odot g_2)}.
}\ealn
This converges to $0$ for $t\to0$. By boundedness of $t\mapsto\vt_t(a)$, this shows that $\vt_t(a)$ is \nbd{\sigma}weakly continuous at $0$. This means, the pre-dual semigroup of $\vt$ is weakly continuous at $0$ and, therefore, by \cite[Proposition I.1.3]{EnNa06}, it is strongly continuous everywhere. In other words, $\vt$ is \nbd{\sigma}weakly continuous.

\item
A \nbd{\sigma}weakly continuous semigroup of endomorphisms is strongly continuous. Indeed,
\beqn{
\norm{\vt_t(a)h-\vt_s(a)h}^2
~=~
\AB{h,\vt_t(a^*a)h}-\AB{\vt_s(a)h,\vt_t(a)h}-\AB{\vt_t(a)h,\vt_s(a)h}+\AB{h,\vt_s(a^*a)h}.
}\eeqn
For fixed $s$ and $t\to s$ this converges to $0$.\qed

\bcor \label{scscor}
A weak type I product system whose generating set of units contains a unital unit admits a (unique, by Theorem \ref{vNscsecthm}) strongly continuous structure making that unit a strongly continuous section.
\ecor

We reprove the following result from Markiewicz and Shalit \cite{MaSha10}.

\bcor\label{MScor}
A weakly continuous contractive CP-semigroup on a von Neumann algebra is strongly continuous. 
\ecor

\proof
Unitalizing if necessary (see \cite[Section 8]{BhSk00}), we may assume the semigroup is Markov. Since the semigroup is bounded, weak continuity implies \nbd{\sigma}weak continuity. If we apply the theorem to the single strongly generating unit $\xi^\odot$ of the GNS-system, we find that the minimal dilation is strongly continuous. It follows that $T_t=\AB{\xi_t,\bullet\xi_t}=\AB{\xi,\vt^\xi_t(\xi\bullet\xi^*)\xi}$ is strongly continuous.\qed

\lf
We say a weak/strong type I system is \hl{contractive} weak/strong type I if the generating set of units can be chosen to contain only contractive units.

\bcor \label{wstIcor}
A contractive weak type I product system is strong type I.
\ecor

\proof
This is equivalent to the statement that for every finite subset $S'=\CB{s_1,\ldots,s_n}$ of $S$ the contractive(!) CP-semigroup $\bfam{b_{ij}}_{i,j}\mapsto\AB{\xi^{s_i}_t,b_{ij}\xi^{s_j}_t}_{i,j}$ on $M_n(\cB)$ is strongly continuous.\qed

\bob\label{0cob}
It is even sufficient to check weak continuity at $t=0$. Indeed, for whatever semigroup $T$ of \nbd{\sigma}weak maps on $\cB$ we have ($T$ weakly continuous at $0$) $\Rightarrow$ ($T$ \nbd{\sigma}weakly continuous at $0$) $\Rightarrow$ ($T_*$ strongly continuous at $0$) $\Rightarrow$ ($T_*$ strongly continuous everywhere) $\Rightarrow$ ($T$ \nbd{\sigma}weakly (\it{a fortiori} weakly) continuous everywhere).
\eob

Let us fix a unital unit from the generating set of a weak type I system. Theorem \ref{scdilthm} tells us that the \nbd{E_0}semigroup constructed from it is strongly continuous, so that the product system inherits a strongly continuous structure having $\xi^\odot$ among the strongly continuous sections (Corollary \ref{scscor}). It is important to know whether the other units in $S$ are strongly continuous sections, too.

The next result gives general criteria to check if a unit is a strongly continuous section. The if-direction of the first part is a strong version of a considerable improvement of $3\Rightarrow 1$ in \cite[Lemma 4.4.11]{BBLS04}. Also the other statements, true also in the situation of \cite[Lemma 4.4.11]{BBLS04} (for continuous units in continuous product systems), are new.

\blem\label{contunilem}
Let $\xi^\odot$ be a unital unit among the strongly continuous sections in a strongly continuous product system $E^\odot$.
\begin{enumerate}
\item\label{CU1}
Suppose $\zeta^\odot$ is another unit in $E^\odot$. Then $\zeta^\odot$ is a strongly continuous section if and only if the functions $t\mapsto\AB{\xi_t,\zeta_t}$, $t\mapsto\AB{\zeta_t,\xi_t}$, and $t\mapsto\AB{\zeta_t,\zeta_t}$ are weakly continuous at $t=0$ (and, therefore, everywhere).

\item\label{CU2}
If $\zeta^\odot$ is a unit among the strongly continuous sections of $E^\odot$, and if $\zeta'^\odot$ is another unit such that the functions $t\mapsto\AB{\zeta_t,\zeta'_t}$, $t\mapsto\AB{\zeta'_t,\zeta_t}$, and $t\mapsto\AB{\zeta'_t,\zeta'_t}$ are weakly continuous at $t=0$, then $\zeta'^\odot$ is a strongly continuous section.
\end{enumerate}
\elem

\bcor \label{Sscscor}
All units in the generating set $S$ of a weak type I product system with a unital unit $\xi^\odot\in S$ are strongly continuous sections. Therefore, by Theorem \ref{vNscsecthm}, the strongly continuous structure does not depend on the choice of $\xi^\odot\in S$.
\ecor

\bcor
If the unit $\zeta'^\odot$ is a strongly continuous section in the strongly continuous subsystem generated by $\zeta^\odot$, then $\zeta'^\odot$ is a strongly continuous section for $E^\odot$, too.
\ecor

\proof[Proof of Lemma \ref{contunilem}.~]
Define the CP-semigroup $S_t=\AB{\zeta_t,\bullet\zeta_t}$ generated by $\zeta^\odot$. If $\zeta^\odot$ is a strongly continuous section, then $S$ is weakly continuous at $t=0$ and, hence, by Observation \ref{0cob} (\nbd{\sigma})weakly continuous everywhere. In particular, $t\mapsto\AB{\zeta_t,\zeta_t}=S_t(\U)$ is weakly continuous. The same argument applies to the mixed inner products. This is the only-if-direction of Part \ref{CU1}.

If we assume strong continuity of the inner products everywhere, then the proof of the if-direction in Part \ref{CU1} is quite similar to the proof of \cite[Lemma 4.4.11]{BBLS04}. But the weak version at $t=0$ requires a refined argument. It is this refinement that allows to show the improvement that, actually, weak continuity (continuity in the case of \cite[Lemma 4.4.11]{BBLS04}) at $t=0$ is sufficient.

So, let now $\zeta^\odot$ be a unit satisfying the weak continuity condition on the inner products. We first show that $S$ is (\nbd{\sigma})weakly continuous. Indeed, since
\bmun{
S_\ve(b)
~=~
\AB{(\zeta_\ve-\xi_\ve)+\xi_\ve,b((\zeta_\ve-\xi_\ve)+\xi_\ve)}
\\
~=~
\AB{\zeta_\ve-\xi_\ve,b(\zeta_\ve-\xi_\ve)}
+
\AB{\xi_\ve,b(\zeta_\ve-\xi_\ve)}
+
\AB{\zeta_\ve-\xi_\ve,b\xi_\ve}
+
T_\ve(b),
}\emun
since, by assumption,
\beq{\label{strong}
\abs{\zeta_\ve-\xi_\ve}^2
~=~
(\AB{\zeta_\ve,\zeta_\ve}-\U)-(\AB{\zeta_\ve,\xi_\ve}-\U)-(\AB{\xi_\ve,\zeta_\ve}-\U)+(\AB{\xi_\ve,\xi_\ve}-\U)
}\eeq
goes to $0$ weakly, and since $T$ is weakly continuous, it follows by Cauchy-Schwarz inequality that $S_\ve(b)-b=(S_\ve(b)-T_\ve(b))+(T_\ve(b)-b))$ goes to zero weakly. From this, by Observation \ref{0cob}, (\nbd{\sigma})weak continuity of $S$ everywhere follows.

So far, all this works for an arbitrary unit that is a strongly continuous section. Now, recall from Theorem\ref{vNscsecthm} that for unital $\xi^\odot$ we may consider the strongly continuous structure of $E^\odot$ to be derived from $\vt^\xi$. Therefore, a section $x$ of $E^\odot$ is strongly continuous if and only if the function $t\mapsto\xi x_t\odot g$ is continuous for all $g\in G$ or, equivalently, if for all $t$ the function $\abs{\xi x_s-\xi x_t}^2$ goes weakly to $0$ for $s\to t$. From $\xi\zeta_{t+\ve}-\xi\zeta_t=\xi(\zeta_\ve-\xi_\ve)\zeta_t$ (see \cite{BBLS04}) and \eqref{strong}, it follows that
\beqn{
\abs{\xi\zeta_{t+\ve}-\xi\zeta_t}^2
~=~
S_t(\AB{\zeta_\ve,\zeta_\ve}-\U)-S_t(\AB{\zeta_\ve,\xi_\ve}-\U)-S_t(\AB{\xi_\ve,\zeta_\ve}-\U)+S_t(\AB{\xi_\ve,\xi_\ve}-\U).
}\eeqn
Since $S_t$ is normal, this implies weak right continuity of $t\mapsto\xi\zeta_t$. Substituting $t>0$ with $t-\ve$ $(\ve<t)$, an appropriate application of Lemma \ref{s*contlem} shows also weak left continuity.

Part \ref{CU2} follows by applying Part \ref{CU1} and \eqref{strong} to $\AB{\zeta'_t,\xi_t}=\AB{\zeta'_t,\xi_t-\zeta_t}+\AB{\zeta'_t,\zeta_t}$. Indeed, while, by assumption, the second term converges weakly to $\U$ at $0$, an application of Cauchy-Schwarz inequality yields (using also Part \ref{CU1} and \eqref{strong}) that the first term converges weakly to $0$. Now we are ready to apply Part \ref{CU1} to the unit $\zeta'^\odot$.\qed

\brem\label{scdilrem}
If the unit $\xi^\odot$ is only contractive, the statement of Theorem \ref{scdilthm} remains true for the \nbd{E}semigroup with \it{pre-assigned product system} constructed in Skeide \cite{Ske08}. The property stated in \cite[Theorem 1.2(2)]{Ske08} (\nbd{W^*}version) guarantees also in this case that the product system inherits a strongly continuous structure from that \nbd{E}semigroup. We do not give details, as in these notes we are only interested in Markov semigroups. But we mention that this implies that Corollary \ref{scscor} and Lemma \ref{contunilem} with its corollaries remain true if we replace `unital unit $\xi^\odot$' with `contractive unit $\xi^\odot$'. Therefore, all these properties depend only on the question if the generating weakly continuous set $S$ contains at least one contractive unit.  The positive  answer to question if this set $S$ is even strongly continuous continuous depends via Corollary \ref{MScor} on the question if the set $S$ can be chosen contractive.
\erem

We close by showing that for \nbd{W^*}algebras no spatial extension of the GNS-system is required. This makes the proof of Theorem \ref{vNHPthm} independent of the construction of the spatial extension from \cite{BLS10} (involving CPD-semigroups and their GNS-systems).

We start with a lemma that allows to determine the product systems of elementary CP-semigroups.

\blem\label{elemtrivlem}
The trivial product system has no proper strongly continuous subsystems.
\elem

\proof
Denote by $E^\odot$ a strongly continuous product subsystem of the trivial one. Observe that all $E_t$ are just strongly closed ideals. Denote by $q_t$ the unique central projection such that $E_t=q_t\cB$. It follows that
\beqn{
q_{s+t}\cB
~=~
E_{s+t}
~=~
E_s\odot E_t
~=~
q_s\cB\odot q_t\cB
~=~
q_sq_t\cB,
}\eeqn
or $q_{s+t}=q_sq_t$. The only semigroups of projections are constant for $t>0$. (Indeed, for $s\ge t$ from the semigroup property it follows $q_{s+t}=q_{s-t}q_t^2=q_{s-t}q_t=q_s$. Therefore, for arbitrary $s,t>0$ and sufficiently big $n$, we get $q_sq_t=q_{ns}q_t=q_{ns+t}=q_{ns}=q_s$. By symmetry, $q_s=q_sq_t=q_t$.) Since $E^\odot$ is assumed strongly continuous, the set $\CB{\AB{x_t,x_t}\colon t>0,x_t\in E_t}$ is strongly total in $\cB$. So, the only possibility for that constant $q_t$ is $q_t=\U$. So, $E_t=\cB$ and the product system is the trivial one.\qed

\lf
Note that this lemma does not require that the subsystem inherits its strongly continuous structure from the containing trivial subsystem. It just has to possess a strongly continuous structure on its own.

\bcor\label{elemtrivcor}
Let $c=\bfam{c_t}_{t\in\R_+}$ be a strongly continuous contractive semigroup in the \nbd{W^*}al\-ge\-bra $\cB$. Then the strongly continuous product system of the elementary CP-semigroup $c_t^*\bullet c_t$ is the trivial one with generating unit $c^\odot=\bfam{c_t}_{t\in\R_+}$.
\ecor

\proof
Effectively, the trivial product system contains the unit $c^\odot:=c$ and $T^c_t=\AB{c_t,\bullet c_t}=c_t^*\bullet c_t$. That unit is weakly continuous, so, by Theorem \ref{scdilthm} the product subsystem $E^\odot$ generated by $c^\odot$ is strongly continuous. By Lemma \ref{elemtrivlem}, $E^\odot$ must be the whole trivial product system. Since the unit $\U^\odot=\bfam{\U}_{t\in\R-+}$ is strongly continuous for both, the continuous structures coincide.\qed

\lf
\bthm\label{vNspthm}
Let $T$ be a strongly continuous normal Markov semigroup on a \nbd{W^*}algebra $\cB$. Then $T$ is spatial if and only if (the strong closure of) its GNS-system is spatial.
\ethm

\proof
The backwards direction we know already. So let us assume that $T$ dominates the elementary CP-semigroup $c_t^*\bullet c_t$ for some strongly continuous semigroup $c$ in $\cB$. In particular, $c$ is contractive. By \cite[Theorem 14.3]{BhSk00} and its proof, for every CP-semigroup $S$ dominated by $T$ there exists a unique contractive positive endomorphism $w^\odot$ of the GNS-system $E^\odot$ of $T$ such that the unit $\zeta^\odot:=\sqrt{w}\xi^\odot=\bfam{\sqrt{w_t}\xi_t}_{t\in\R_+}$ generates $S$. If $S$ is elementary, then by Corollary \ref{elemtrivcor}, the subsystem generated by that unit is the trivial one. So, the only thing that remains to be shown, is that the unit $\zeta^\odot$ is strongly continuous, because in that case, by Lemma \ref{contunilem}\eqref{CU2}, the unital central unit $\bfam{\U}_{t\in\R_+}$ of that subsystem is strongly continuous also in $E^\odot$.

The map $t\mapsto\AB{\zeta_t,\zeta_t}=c_t^*c_t$ is weakly continuous. Observe that $\U-\AB{\xi_t,\zeta_t}=\AB{\xi_t,(\U-\sqrt{w_t})\xi_t}$. From
\beqn{
0
~\le~
\AB{\xi_t,(\U-\sqrt{w_t})\xi_t}
~\le~
\AB{\xi_t,(\U-\sqrt{w_t})(\U+\sqrt{w_t})\xi_t}
~=~
\AB{\xi_t,(\U-w_t)\xi_t}
~=~
\AB{\xi_t,\xi_t}-\AB{\zeta_t,\zeta_t},
}\eeqn
it follows that also the map $t\mapsto\AB{\xi_t,\zeta_t}=\AB{\zeta_t,\xi_t}$ is weakly continuous at $t=0$, as required by Lemma \ref{contunilem}\eqref{CU1}.\qed

\brem
Apart from being crucial for the proof of existence of a Hudson-Parthasarathy dilation for a spatial Markov semigroup, this result is also important for the classification of strong type I systems. A slight modification asserts that a strongly continuous weak type I system is spatial if and only if the \it{CPD-semigroup} generated by the generating set $S$ is spatial in the sense of Skeide \cite{Ske10}.
\erem

\section*{Appendix B: $E_0$--Semigroups and representations for strongly continuous product systems}
\addcontentsline{toc}{section}{Appendix B: $E_0$--Semigroups and representations for strongly continuous product systems}
	\stepcounter{section}
        \renewcommand{\thesection}{B}

The principal scope of this appendix is to prove Theorem \ref{vNexithm}: Every strongly full strongly continuous product system is the product system associated with a strongly continuous normal \nbd{E_0}semigroup on some $\sB^a(E)$. However, since the existence of a nondegenerate faithful representation of a faithful product system is a closely related problem (dual to Theorem \ref{vNexithm} in the sense of commutant of von Neumann correspondences in Skeide \cite{Ske03c}), we include this result here. As corollaries we reprove a result by Arveson and Kishimoto \cite{ArKi92}, we show that that CP-semigroups admit so-called \it{elementary dilations}, and we show that \it{faithful} strongly continuous product systems have a strongly continuous \it{commutant system}.

According to our convention in these notes, $\cB$ is a \nbd{W^*}algebra acting a Hilbert space $G$. For the balance of this appendix, $E^\odot$ is a product system of \nbd{W^*}correspondences over $\cB$ that is strongly continuous with respect to a family of isometric embeddings $i_t\colon E_t\rightarrow\wh{E}$ in the sense of Definition \ref{scPSdefi}. We also use the other notations introduced there.

We start by proving some properties that hold for all strongly continuous product systems, strongly full or not, and faithful left action or not. First of all, we note that the embeddings $i_t\colon E_t\rightarrow\wh{E}$ give rise to embeddings $H_t:=E_t\odot G\rightarrow\wh{H}:=\wh{E}\odot G$, also denoted by $i_t$, defined by $x_t\odot g\mapsto (i_tx_t)\odot g$. We, therefore, may speak about continuous sections $h=\bfam{h_t}_{t\in\R_+}$ of $E^\odot\odot G$, in the sense that $t\mapsto i_th_t$ is continuous. We denote the set of all \hl{continuous sections} of $E^\odot\odot G$ by $CS_i(E^\odot\odot G)$. By definition, whenever $x,y\in CS_i^s(E^\odot)$ and $g\in G$, then the functions $t\mapsto i_tx_t\odot g$ and $(s,t)\mapsto i_{s+t}(x_sy_t)\odot g$ are continuous.

\bcor\label{lacontcor}
If $y\in CS_i^s(E^\odot)$ and $g\in G$, then for every $b\in\cB$ the function $t\mapsto i_t(by_t)\odot g$ is continuous.
\ecor

\proof
Choose a strongly continuous section $x$ such that $x_0=b$, and in the continuous function $(s,t)\mapsto i_{s+t}(x_sy_t)\odot g$ put $s=0$.\qed

\lf
The following lemma is just \it{Dini's theorem} for nets. Later on, we will apply its corollary to the functionals $\AB{y_t\odot g,\bullet y_t\odot g}$ on $\cB$.

\blem\label{uniapprlem}
Let $\bfam{\vp_t}_{t\in\R_+}$ be a family of normal positive linear functionals on $\cB$ such that $t\mapsto\vp_t(c)$ is continuous for all $c\in\cB$. Suppose the net $\bfam{c_\lambda}_{\lambda\in\Lambda}$ in $\cB$ increases to $c\in\cB$. Then for every $0\le a<b<\infty$ and every $\ve>0$ there exists a $\lambda_0\in\Lambda$ such that
\beqn{
\vp_\alpha(c-c_\lambda)
~<~
\ve.
}\eeqn
for all $\alpha\in\SB{a,b}$ and all $\lambda\ge\lambda_0$.
\elem

\bcor\label{uniapprcor}
The positive linear functional on $\cB$ defined by $c\mapsto\int_a^b\vp_\alpha(c)\,d\alpha$ is normal.
\ecor

\proof[Proof of Lemma \ref{uniapprlem}.~]
This is a standard application of compactness of the interval $\SB{a,b}$, like many others that follow in this appendix. For each $\beta\in\SB{a,b}$ choose $\lambda_\beta\in\Lambda$ such that $\vp_\beta(c-c_\lambda)<\ve$ for all $\lambda\ge\lambda_\beta$. Define $I_\beta$ to be the largest subinterval of $\SB{a,b}$ containing $\beta$ such that $\vp_\alpha(c-c_{\lambda_\beta})<\ve$ for all $\alpha\in I_\beta$. Since $c_\lambda$ increases to $c$, we get $\vp_\alpha(c-c_\lambda)<\ve$ for all $\alpha\in I_\beta$ and all $\lambda\ge\lambda_\beta$. Since $\alpha\mapsto\vp_\alpha(c-c_{\lambda_\beta})$ is continuous, every $I_\beta$ is open in $\SB{a,b}$. Since $I_\beta\ni\beta$, the family of all $I_\beta$ forms an open cover of the compact interval $\SB{a,b}$. So, we may choose $\beta_1,\ldots,\beta_m$ such that the union over $I_{\beta_i}$ is $\SB{a,b}$. Since every $\alpha\in\SB{a,b}$ is contained in at least one of the intervals $I_{\beta_i}$, it follows that $\lambda_0=\max_{i=1,\ldots,m}\lambda_{\beta_i}$ does the job.\qed

\lf
The following density result is analogue to \cite[Proposition 2.6]{Ske11a}. The proof applies \it{a fortiori} also to the situation in \cite[Proposition 2.6]{Ske11a}, and simplifies its proof quite a bit.

\bprop\label{Hcontprop}
Every continuous section $h\in CS_i(E^\odot\odot G)$ may be approximated locally uniformly by elements in $\ls CS_i^s(E^\odot)\odot G$. Moreover:
\begin{enumerate}
\item\label{Hcsdef}
For every $k_t\in H_t$ we can find a continuous section $h\in CS_i(E^\odot\odot G)$ such that $h_t=k_t$.

\item \label{Hcpsdef}
For every pair $x\in CS_i^s(E^\odot)$ and $h\in CS_i(E^\odot\odot G)$ of sections the function
$$
(s,t)
~\longmapsto~
i_{s+t}(x_sh_t)
$$
is continuous.
\end{enumerate}
\eprop

\proof
Once we have the density statement, \eqref{Hcsdef} is a standard result about continuous fields of Banach spaces (proved, for instance, like \cite[Proposition 7.9]{Ske03b}), and \eqref{Hcpsdef} follows by three epsilons, approximating $h\in CS_i(E^\odot\odot G)$ with an element in $\ls CS_i^s(E^\odot)\odot G$ on a suitably big interval. So, let us prove the density statement.

Let $h\in CS_i(E^\odot\odot G)$ and choose $0\le a<b<\infty$ and $\ve>0$. By Definition \ref{scPSdefi}, for every $\beta\in\SB{a,b}$ there exists a section $h^\beta$ in $\ls CS_i^s(E^\odot)\odot G$ such that $\snorm{h_\beta-h^\beta_\beta}<\ve$. For every $\beta$ define $I_\beta$ to be the largest interval containing $\beta$ such that $\snorm{h_\alpha-h^\beta_\alpha}<\ve$ for all $\alpha\in I_\beta$. Every $I_\beta$ is open in $\SB{a,b}$ and contains at least $\beta$. Therefore, the family of all $I_\beta$ forms an open cover of the compact interval $\SB{a,b}$. So, we may choose $\beta_1,\ldots,\beta_m$ such that the union over $I_{\beta_i}$ is $\SB{a,b}$. By standard theorems about \it{partitions of unity} there exist continuous functions $\vp_i$ on $\SB{a,b}$ with the following properties:
\baln{
0
~\le~
&
\vp_i,
&
\vp_i\upharpoonright I_{\beta_i}^\complement
&
~=~
0,
&
\sum_{i=1}^m\vp_i
&
~=~
1.
}\ealn
From these properties, one easily verifies that $\norm{h_\alpha-\sum_{i=1}^m\vp_i(\alpha)h^{\beta_i}_\alpha}<\ve$ for all $\alpha\in\SB{a,b}$. This shows that $\sum_{i=1}^m\vp_ih^{\beta_i}\in\ls CS_i(E^\odot)\odot G$ approximates $h$ uniformly up to $\ve$ on the interval $\SB{a,b}$.\qed

\bthm\label{scisothm}
Let $i_t\colon E_t\rightarrow E^i$ and $k_t\colon E_t\rightarrow E^k$ be two strongly continuous structures on the product system $E^\odot=\bfam{E_t}_{\in\R_+}$. If the identity morphism is a strongly continuous isomorphism from $E^\odot$ with respect to the embeddings $\bfam{i_t}$ to $E^\odot$ with respect to the embeddings $\bfam{k_t}$, then it is a strongly continuous isomorphism for the other direction, too.
\ethm

\proof
For continuous product systems and continuous isomorphisms, this is \cite[Theorem 2.2]{Ske09a}. In the proof of that theorem we also mentioned that the theorem, actually, is a statement that is valid for every continuous field of Banach spaces that takes its continuous structure from a family of isometric embeddings into a fixed Banach space. Making use of Proposition \ref{Hcontprop}, we will apply \cite[Theorem 2.2]{Ske09a} to the the continuous field of Hilbert spaces $\bfam{H_t}_{t\in\R_+}$ with respect to the two embeddings.

What we have to show is that if $x\odot g\in CS_i(E^\odot\odot G)$ for every $g\in G$ implies $x\odot g\in CS_k(E^\odot\odot G)$ for every $g\in G$, then $x\odot g\in CS_k(E^\odot\odot G)$ for every $g\in G$ implies $x\odot g\in CS_i(E^\odot\odot G)$ for every $g\in G$. What we have from \cite[Theorem 2.2]{Ske09a} is the statement that if $h\in CS_i(E^\odot\odot G)$ implies $h\in CS_k(E^\odot\odot G)$, then $h\in CS_k(E^\odot\odot G)$ implies $h\in CS_i(E^\odot\odot G)$. It is, therefore, sufficient to show that if $x\odot g\in CS_i(E^\odot\odot G)$ for every $g\in G$ implies $x\odot g\in CS_k(E^\odot\odot G)$ for every $g\in G$, then $h\in CS_i(E^\odot\odot G)$ implies $h\in CS_k(E^\odot\odot G)$. But, this last statement follows if we approximate $h$ locally uniformly in the by sections in $\ls C_i^s(E^\odot)\odot G\subset C_k^s(E^\odot)$ and take into account, like in the proof of \cite[Theorem 2.2]{Ske09a}, that the locally uniform approximation does not depend on the choice of the continuous structure $i$ or $k$.\qed

\subsection*{B.1: $E_0$--Semigroups}

After these general properties of strongly continuous product systems, we now come to the proof of Theorem \ref{vNexithm}. For this part of this appendix we shall assume that $E^\odot$ is strongly full.

We already mentioned in Section \ref{vNalgSEC} that the algebraic part of the proof in Skeide \cite{Ske06} for the Hilbert space case allows to construct an \nbd{E_0}semigroup that is is definitely not continuous and acts on a $\sB^a(E)$ whose pre-dual that cannot be separable. The \nbd{E_0}semigroup is constructed in terms of a left dilation, and the construction, using the identifications in \eqref{idea}, involves the choice of a left dilation of the discrete subsystem $\bfam{E_n}_{n\in\N_0}$ of $E^\odot$. Existence of a left dilation of $\bfam{E_n}_{n\in\N_0}$ is granted by Skeide \cite[Theorem 6.3]{Ske09}.

Like in Skeide \cite{Ske07}, where the \nbd{C^*}case is treated, also here it is convenient to adapt Arveson's construction in \cite{Arv06} of our \nbd{E_0}semigroup in \cite{Ske06}, rather than using our construction directly. (We refer to \cite{Ske07} for motivation.) This construction is based on the choice of a unit vector $\xi_1\in E_1$. It is a feature of continuous product systems that every $E_t$ has a unit vector; see \cite[Lemma 3.2]{Ske07}. Although we dare to conjecture that the same is also true for strongly continuous product systems of \nbd{W^*}correspondences that are strongly full (it is difficult to imagine a counter example), we could not yet prove it. (Note that \cite[Lemma 3.2]{Ske07} is about \bf{every} continuous product system of \nbd{C^*}correspondences over a unital \nbd{C^*}algebra $\cB$, which, therefore, are full automatically. If we drop strong fullness, then the statement for product systems of \nbd{W^*}correspondence is surely false: Consider \cite[Example 4.13]{Ske04p} for a non-unital \nbd{C^*}algebra and close it strongly.)

One basic idea in the construction of the left dilation of the discrete subsystem $\bfam{E_n}_{n\in\N_0}$ in \cite{Ske09} was that even if the strongly full \nbd{W^*}correspondence $E_1$ over $\cB$ has no unit vector, then $\ol{M_\en(E_1)}^s$ has one; see \cite[Proposition 6.2]{Ske09}. Here $\en$ is a sufficiently big cardinal number and $\ol{M_\en(E_1)}^s=E_1\sbars{\otimes}\sB(\C^\en)$ is the (spatial) external tensor product of \nbd{W^*}modules, which is a \nbd{W^*}module over $\cB\sbars{\otimes}\sB(\C^\en)=\ol{M_\en(\cB)}^s$. The elements of $\ol{M_\en(E_1)}^s$ are understood best as \nbd{E_1}valued matrices, and the \nbd{W^*}module operations are the natural matrix operations; see \cite[Section 6]{Ske09} for details. One easily verifies that $\sB^a(\ol{M_\en(E_1)}^s)=\sB^a(E_1)\sbars{\otimes}\sB(\C^\en)=\ol{M_\en(\sB^a(E_1))}^s$.

We follow the construction in Skeide \cite{Ske09a} where we proved the result of \cite{Ske07}  for not necessarily unital \nbd{C^*}algebras. We also refer to \cite{Ske09a} for a detailed motivation and an explanation why what follows is the proper generalization of Arveson's idea \cite{Arv06}. Technically, the whole proof of the \nbd{W^*}case here, is very similar to \cite{Ske09a} (involving also technical results from \cite{Ske07} and \cite{Ske11a}). For reasons of space we dispense with giving full proofs of these technical results, and often refer to either identical or at least very similar proofs in the cited papers.

So, let us fix a unit vector $\Xi_1\in\ol{M_\en(E_1)}^s$. To facilitate notation we fix a set $S$ of cardinality $\#S=\en$ and denote the elements of $\ol{E_\alpha^\en}^s$ as $X_\alpha=\bfam{X_\alpha^s}_{s\in S}$. Also, $\Xi_1$ is given by the matrix  $\bfam{(\Xi_1)_{ss'}}_{ss'}$. Like a vector $\xi_1\in E_1$ can act on $x_\alpha\in E_\alpha$ as $\xi_1x_\alpha\in E_{\alpha+1}$, the vector $\Xi_1$ can act on $X_\alpha\in\ol{E_\alpha^\en}^s$ as $\Xi_1X_\alpha:=\bfam{\sum_{s'\in S}(\Xi_1)_{ss'}X_\alpha^{s'}}_{s\in S}\in\ol{E_{\alpha+1}^\en}^s$.

Next, we define the direct integrals we need. The family of embeddings $i_t\colon E_t\rightarrow\wh{E}$ gives rise to embeddings $i_t^\en\colon\ol{E_t^\en}^s\rightarrow\ol{\wh{E}^\en}^s$. Every section $X=\bfam{X_t}_{t\in\R_+}$ with $X_t\in\ol{E_t^\en}^s$ gives rise to a function $t\mapsto X(t):=i_t^\en X_t$ with values in $\ol{\wh{E}^\en}^s$. We denote by
\beqn{
CS_i^{\en,s}(E^\odot)
~=~
\BCB{\,X\colon t\mapsto X(t)\text{~is strongly continuous~}}
}\eeqn
the set of all sections that are strongly continuous. Let $0\le a<b<\infty$. By $\int_a^b E_\alpha^\en\,d\alpha$ we understand the self-dual extension of the pre-Hilbert \nbd{\cB}module that consists of continuous sections $X\in CS_i^{\en,s}(E^\odot)$ restricted to $\RO{a,b}$ with inner product
\beqn{
\AB{X,Y}_{\SB{a,b}}
~:=~
\int_a^b\AB{X_\alpha,Y_\alpha}\,d\alpha
~=~
\int_a^b\AB{X(\alpha),Y(\alpha)}\,d\alpha.
}\eeqn
By an application of the \it{principle of uniform boundedness}, all strongly continuous sections are bounded on the compact interval $\SB{a,b}$. Therefore, the integral exists as a \it{Riemann integral} in the weak operator topology of $\sB(G)$.

\bob
Some care is in place in calculating inner products of arbitrary elements $X$ and $Y$, when thinking of them as an ``integral'' over ``sections''. Even when the bundle is trivial, these elements need no longer have an interpretation as sections with values in the fibers. \cite[Example 4.3.13]{Ske01} shows that this already fails under norm completion.
\eob

\bprop\label{srcllprop}
$\int_a^b E_\alpha^\en\,d\alpha$ contains as a pre-Hilbert submodule the space $\es\eR_{\RO{a,b}}$ of restrictions to $\RO{a,b}$ of those sections $X$ for which $t\mapsto X(t)$ has strong left limits everywhere and is strongly right continuous with a finite number of jumps in $\RO{a,b}$.
\eprop

\proof
The proof is more or less like that of \cite[Proposition 4.2]{Ske07} for right continuous (not only strongly continuous) functions with left limits (not only strong limits) in each point. Just that now the approximation must be done for each of the functions $t\mapsto X(t)\odot g$ separately. (A linear subspace $C$ of $\es\eR_{\RO{a,b}}$ is strongly dense if for each $g\in G$, each $\ve>0$, and each $X\in\es\eR_{\RO{a,b}}$ the function $t\mapsto X(t)\odot g$ can be approximated up to $\ve$ by a function $t\mapsto Y_t\odot g$ with $Y_t\in C$. By the proof \cite[Proposition 4.2]{Ske07}, this is true for $C$ being the space of strongly continuous sections restricted to $\RO{a,b}$, because the approximation of the jumps is done in a way that does not depend on $g$. Better: It is uniform in $\norm{g}\le1$.)\qed

\bprop\label{stfullprop}
$\int_a^b E_\alpha^\en\,d\alpha$ is strongly full.
\eprop

\proof
Let $b=\AB{x_a,x_a}\in\cB$ for some $x_a\in E_a$. These elements are strongly total, so it is sufficient if we strongly approximate each such $b$. By definition, there exists a section $y\in CS_i^s(E^\odot)$ such that $y_a=x_a$. It is not difficult to show that for the strongly right continuous sections $y^\lambda:=\bfam{\frac{y_t\I_{\RO{a,a+\lambda}(t)}}{\sqrt{\lambda}}}_{t\in\R_+}$ $(\lambda>0)$, the expression $\AB{y^\lambda,y^\lambda}_{\SB{a,b}}$ converges weakly to $b$ for $\lambda\to0$. So the the span of all $\AB{y^\lambda,y^\lambda}$ is weakly, hence, strongly dense in $\cB$.\qed

\lf
Let $\s\sS$ denote the right \nbd{\cB}module of all sections $X$ that are \hl{locally $\es\eR$}, that is, for every $0\le a<b<\infty$ the restriction of $X$ to $\RO{a,b}$ is in $\es\eR_{\RO{a,b}}$, and that are \hl{stable} with respect to the unit vector $\Xi_1$, that is, there exists an $\alpha_0\ge0$ such that
\beqn{
X_{\alpha+1}
~=~
\Xi_1X_\alpha
}\eeqn
for all $\alpha\ge\alpha_0$. By $\s\sN$ we denote the subspace of all sections in $\s\sS$ which are eventually $0$, that is, of all sections $X\in\s\sS$ for which there exists an $\alpha_0\ge0$ such that $X_\alpha=0$ for all $\alpha\ge\alpha_0$. A straightforward verification shows that
\beqn{
\AB{X,Y}
~:=~
\lim_{m\to\infty}\int_m^{m+1}\AB{X(\alpha),Y(\alpha)}\,d\alpha
}\eeqn
defines a semiinner product on $\s\sS$ and that $\AB{X,X}=0$ if and only if $X\in\s\sN$. Actually, we have
\beqn{
\AB{X,Y}
~=~
\int_T^{T+1}\AB{X(\alpha),Y(\alpha)}\,d\alpha
}\eeqn
for all sufficiently large $T>0$; see \cite[Lemma 2.1]{Arv06}. So, $\s\sS/\s\sN$ becomes a pre-Hilbert module with inner product $\AB{X+\s\sN,Y+\s\sN}:=\AB{X,Y}$. By $E$ we denote its self-dual extension.

\bprop\label{denseprop}
For every section $X$ and every $\alpha_0\ge0$ define the section $X^{\alpha_0}$ as
\beqn{
X^{\alpha_0}_\alpha
~:=~
\begin{cases}
0&\alpha<\alpha_0
\\
\Xi_1^nX_{\alpha-n}&\alpha\in\RO{\alpha_0+n,\alpha_0+n+1},n\in\N_0.
\end{cases}
}\eeqn

If $X$ is in $CS_i^{\en,s}(E^\odot)$, then $X^{\alpha_0}$ is in $\s\sS$. Moreover, the set $\bCB{X^{\alpha_0}+\s\sN\colon X\in CS_i^{\en,s}(E^\odot),\alpha_0\ge0}$ is a strongly dense submodule of $E$.
\eprop

\proof
The proof is like that of \cite[Proposition 4.3]{Ske07}. It only takes a moments thought to convince oneself that the approximating sequence in that proof may be replaced without problem by a strongly approximating net.\qed

\lf
After these preparations it is completely plain to see that for every $t\in\R_+$ the map $X\odot y_t\mapsto Xy_t$, where the section $Xy_t\in E$ is defined by
\beqn{
(Xy_t)_\alpha
~=~
\begin{cases}
X_{\alpha-t}y_t&\alpha\ge t,
\\
0&\text{else}
\end{cases}
}\eeqn
($X_\alpha y_t\in\ol{E_{\alpha+t}^\en}^s$ is the componentwise product $\ol{E_\alpha^\en}^s\times E_t\rightarrow\ol{E_{\alpha+t}^\en}^s$ of $X_\alpha$ and $y_t$), defines an isometry $v_t\colon E\sodots\!\! E_t\rightarrow E$, and that these isometries iterate associatively with the product system structure.

\bprop\label{surprop}
Each $v_t$ is surjective.
\eprop

\proof
By Proposition \ref{denseprop} it is sufficient to approximate every section of the form $X^{\alpha_0}$ with $X\in CS_i^{\en,s}(E^\odot),\alpha_0\ge0$ in the (semi-)inner product of $\s\sS$ by finite sums of sections of the form $Yz_t$ for $Y\in\s\sS,z_t\in E_t$. As what the section does on the finite interval $\RO{0,t}$ is not important for the inner product, we may even assume that $\alpha_0\ge t$. And as in the proof of Proposition \ref{denseprop} the approximation can be done by approximating $X$ in $\es\eR_{\RO{\alpha_0,\alpha_0+1}}$ and then extending the restriction to $\RO{\alpha_0,\alpha_0+1}$ stably to the whole axis. (This stable extension is the main reason why we worry to introduce the subspace of strongly right continuous sections.) On $\es\eR_{\RO{\alpha_0,\alpha_0+1}}$, however, the approximation may be done for each direct summand in $E$ separately. (Since for $\alpha\in\RO{\alpha_0,\alpha_0+1}$ the operator $\Xi_1$ in the definition of $X^{\alpha_0}$ does not occur, the maps $v_t$ decompose into the components of $\ol{E^\en}^s$.) In other words, it is sufficient to prove the density statement only on the interval $\RO{\alpha_0,\alpha_0+1}$ and for $\en=1$.

The continuous version of Proposition \ref{surprop} for $\en=1$ is done in \cite[Proposition 4.6]{Ske07}. Like in the proof of  Proposition \ref{Hcontprop}, for the strongly continuous version we have to modify the proof of \cite[Proposition 4.6]{Ske07}. As this modification is a bit too much for just referring to the continuous version, we give full detail.

Let $\alpha_0\ge t$ and let $x$ be a strongly continuous section and choose $g\in G$. We will approximate the continuous section $\alpha\mapsto x_\alpha\odot g$ uniformly on the compact interval $\SB{\alpha_0,\alpha_0+1}$ (and, therefore, in $L^2$) by finite sums over sections of the form $\alpha\mapsto y_{\alpha-t}z_t\odot g$. Choose $\ve>0$. For every $\beta\in\SB{\alpha_0,\alpha_0+1}$ choose $n^\beta\in\N,y_k^\beta\in E_{\beta-t},z_k^\beta\in E_t$ such that $\snorm{(x_\beta-\sum_{k=1}^{n^\beta}y_k^\beta z_k^\beta)\odot g}<\ve$. Choose continuous sections $\bar{y}_k^\beta=\bfam{(\bar{y}_k^\beta)_\alpha}_{\alpha\in\R_+}\in CS_i^s(E^\odot)$ such that $(\bar{y}_k^\beta)_{\beta-t}=y_k^\beta$. For every $\beta$ chose the maximal interval $I_\beta\subset\SB{\alpha_0,\alpha_0+1}$ containing $\beta$ such that $\bnorm{\bfam{x_\alpha-\sum_{k=1}^{n^\beta}(\bar{y}_k^\beta)_{\alpha-t}z_k^\beta}\odot g}<\ve$ for all $\alpha\in I_\beta$. Like in the other proofs, $I_\beta$ is open in $\SB{\alpha_0,\alpha_0+1}$ and contains at least $\beta$. So, we may choose finitely many  $\beta_1,\ldots,\beta_m\in\SB{\alpha_0,\alpha_0+1}$ such that the union of all $I_{\beta_i}$ is $\SB{\alpha_0,\alpha_0+1}$. Putting $I_i:=I_{\beta_i}\backslash(I_{\beta_1}\cup\ldots\cup I_{\beta_{i-1}})$, we define a finite partition $I_1,\ldots,I_m$ of $\SB{\alpha_0,\alpha_0+1}$. Taking away the point $\alpha_0+1$ and adjusting the endpoints of the $I_i$ suitably, we may assume that all $I_i$ are right open. Denote by $\I_i$ the indicator function of $I_i$. Then, restriction of the piecewise continuous section
\beqn{
\alpha
~\longmapsto~
\begin{cases}
0&\alpha<t
\\
\displaystyle\sum_{i=1}^m\sum_{k=1}^{n_{\beta_i}}(\bar{y}_k^{\beta_i})_{\alpha-t}z_k^{\beta_i}\odot g\I_i(\alpha)&\alpha\ge t
\end{cases}
}\eeqn
to $\RO{\alpha_0,\alpha_0+1}$ is in $\es\eR_{\RO{\alpha_0,\alpha_0+1}}$ and approximates $\alpha\mapsto x_\alpha\odot g$ uniformly on $\RO{\alpha_0,\alpha_0+1}$ up to $\ve$.\qed

\lf
So, the $v_t$ form a left dilation of $E^\odot$ to $E$. This left dilation is \hl{strongly continuous} in the following sense. 

\bprop\label{scdprop}
For every $X\in E$, every strongly continuous section $y\in CS_i^s(E^\odot)$, and every $g\in G$ the function $t\mapsto Xy_t\odot g$ is continuous.
\eprop

\proof
In the proof of the analogue \cite[Proposition 4.7]{Ske07} for continuous product systems, we made use of the fact that, there, $X$ could be approximated in norm by a section of the form $X^{\alpha_0}$, so that it was sufficient to show the statement of \cite[Proposition 4.7]{Ske07} only for such sections. Also here proving the statement first for sections $X^{\alpha_0}$ will be an important step in the proof. However, Proposition \ref{denseprop} guarantees only strong approximation of $X$, and the argument that this is sufficient differs considerably from the proof in \cite{Ske07}.

So, suppose for a moment we had proved that $t\mapsto X^{\alpha_0}y_t\odot g$ is continuous for every $X\in C_i^{\en,s}(E^\odot)$, every $\alpha_0\ge0$, every $y\in CS_i^s(E^\odot)$, and every $g\in G$. In order to prove that $t\mapsto Yz_t\odot g'$ is continuous for every $Y\in E$, every $z\in CS_i^s(E^\odot)$, and every $g'\in G$, it is sufficient to show that $t\mapsto\norm{Yz_t\odot g'}$ is continuous, and that $t\mapsto Yz_t\odot g'$ is weakly continuous. Since that function is bounded uniformly on finite intervals, it is even sufficient to check weak continuity on a total subset of $H$, only. Continuity of the norm follows via $\norm{Yz_t\odot g'}=\sqrt{\AB{z_t\odot g',\AB{Y,Y}z_t\odot g'}}$ from Corollary \ref{lacontcor}. To see weak continuity at a fixed point $t$, we observe that the elements $X^{\alpha_0}y_t\odot g$ form a total subset of $H$. We compute
\bmun{
\AB{X^{\alpha_0}y_t\odot g,Yz_s\odot g'-Yz_t\odot g'}
\\
~=~
\AB{X^{\alpha_0}y_t\odot g-X^{\alpha_0}y_s\odot g,Yz_s\odot g'}
+
\AB{X^{\alpha_0}y_s\odot g,Yz_s\odot g'}
-
\AB{X^{\alpha_0}y_t\odot g,Yz_t\odot g'}.
}\emun
For $s\to t$ the first summand goes to $0$ by continuity of $t\mapsto X^{\alpha_0}y_t\odot g$. Like for continuity of the norm, also the difference of the last two summands goes to $0$ by Corollary \ref{lacontcor}.

It remains, therefore, to show continuity of $t\mapsto X^{\alpha_0}y_t\odot g$ for  $X^{\alpha_0}$ for $X\in CS_i^{\en,s}(E^\odot)$ and $\alpha_0\ge0$. To calculate $\norm{(X^{\alpha_0}y_t-X^{\alpha_0}y_s)\odot g}^2$ we have to integrate over $\alpha$ the values of $\snorm{(X^{\alpha_0}_{\alpha-t}y_t-X^{\alpha_0}_{\alpha-s}y_s)\odot g}^2$ for $\alpha$ in some unit interval such that  $\alpha-t$ and $\alpha-s$ are not smaller than $\alpha_0$. So
\beq{\label{secnd}
\norm{(X^{\alpha_0}y_t-X^{\alpha_0}y_s)\odot g}^2
~=~
\int_d^{d+1}\norm{(X^{\alpha_0}_{\alpha+s}y_t-X^{\alpha_0}_{\alpha+t}y_s)\odot g}^2\,d\alpha
}\eeq
for all $d\ge\alpha_0$. The function $(\alpha,t)\mapsto X^{\alpha_0}_\alpha y_t\odot g$ is uniformly continuous on each of the intervals $\RO{\alpha_0+n,\alpha_0+n+1}\times\SB{a,b}$ and it is bounded on each $\R_+\times\SB{a,b}$. We fix a $t$. Choose $d\ge\alpha_0$ such that $d+t=n+\alpha_0$ for a suitable $n\in\N_0$. Then $\alpha+t$ in the integral in \eqref{secnd} goes over the interval $\RO{\alpha_0+n,\alpha_0+n+1}$. Choose $\ve\in\bfam{0,\frac{1}{2}}$. Then in
\beqn{
\norm{(X^{\alpha_0}y_t-X^{\alpha_0}y_s)\odot g}^2
\\
~=~
{\Biggl[\textstyle\int_d^{d+\ve}+\int_{d+\ve}^{d+1-\ve}+\int_{d+1-\ve}^{d+1}\Biggr]}~\bnorm{\smash{(X^{\alpha_0}_{\alpha+s}y_t-X^{\alpha_0}_{\alpha+t}y_s)\odot g}}^2\,d\alpha
}\eeqn
the first and the last integral are bounded by $\ve$ times a constant which can be chosen independent of $s$ as long as $s$ varies in a bounded set. If we choose $s\in\bfam{t-\ve,t+\ve}$, in the middle integral also $\alpha+s$ is in the same interval $\RO{\alpha_0+n,\alpha_0+n+1}$, so that both $X^{\alpha_0}_{\alpha+s}y_t$ and $X^{\alpha_0}_{\alpha+t}y_s$ depend uniformly continuously on $(\alpha,s)$ in $(d+\ve,d+1-\ve)\times(t-\ve,t+\ve)$. In particular, if $s$ is sufficiently close to $t$, then both $X^{\alpha_0}_{\alpha+s}y_t\odot g$ and $X^{\alpha_0}_{\alpha+t}y_s\odot g$ are close to their common limit $X^{\alpha_0}_{\alpha+t}y_t\odot g$ uniformly in $\alpha$. It follows that the middle integral goes to $0$ for $s\to t$. Sending also $\ve\to0$, the proposition is proved.\qed

\bcor
The \nbd{E_0}semigroup $\vt^v$ is strongly continuous.
\ecor

\proof
Fix $a\in\sB^a(a)$. Since $\vt^v_t(a)$ is bounded uniformly in $t$ by $\norm{a}$, we may check strong continuity at $t$ on the total subset of elements of the form $Xy_t\odot g$ ($X\in E$, $y\in CS_i^s(E^\odot)$, $g\in G$). We find
\bmun{
(\vt^v_s(a)-\vt^v_t(a))Xy_t\odot g
\\
~=~
\bfam{\vt^v_s(a)(Xy_t\odot g)-\vt^v_s(a)(Xy_s\odot g)}
+
\bfam{\vt^v_s(a)(Xy_s\odot g)-\vt^v_t(a)(Xy_t)\odot g}
\\
~=~
\vt^v_s(a)\bfam{Xy_t\odot g-Xy_s\odot g}
+
\bfam{(aX)y_s\odot g-(aX)y_t\odot g}.
}\emun
By Proposition \ref{scdprop}, both expressions are small, whenever $s$ is sufficiently close to $t$.\qed

\bcor
The continuous structure induced by the \nbd{E_0}semigroup $\vt^v$ coincides with the continuous structure of $E^\odot$.
\ecor

\proof
By amplifying if necessary, we may assume that $E$ has a unit vector $\xi$. Each section $y$ corresponds to the section $\bfam{\xi y_t}_{t\in\R_+}$ with respect to the new embedding of $E^\odot$ into $E$ induced by $\vt^v$. By Proposition \ref{scdprop} (which, clearly, remains valid also under amplification of $\vt^v$), if $y$ is strongly continuous, then so is $t\mapsto\xi y_t$. By Theorem \ref{scisothm}, the two continuous structures coincide.\qed

\lf
This ends the proof of Theorem \ref{vNexithm}.

\bob
Note that the proofs of the two preceding corollaries do not depend on the concrete form of the left dilation. We, therefore, showed the following more general statement: If $v_t$ is a left dilation of a strongly continuous product system $E^\odot$ that is strongly continuous in the sense of Proposition \ref{scdprop}, then the induced \nbd{E_0}semigroup $\vt^v$ is strongly continuous and the strongly continuous structure induced by that \nbd{E_0}semigroup upon $E^\odot$ coincides with the original one.
\eob

\brem \label{countvNrem}
We do not need countability hypotheses. But of course, like in the \nbd{C^*}case, if countability hypotheses are fulfilled, then the constructions \nbd{E_0}semigroup $\leftrightarrow$ product system in either direction preserve these. If $\sB^a(E)$ and, thus, $\cB$ have separable pre-duals, then the product system $E^\odot$ of a strongly continuous \nbd{E_0}semigroup on $\sB^a(E)$ is \hl{countably generated} in the sense that there exists a countable subset of $C_i^s(E^\odot)$ that generates $C_i^s(E^\odot)$ by locally uniform strong limits. Conversely, if $\cB$ has separable pre-dual and if $E^\odot$ is countably  generated, then $E$ and, therefore, also $\sB^a(E)$ has separable pre-dual.
\erem

\subsection*{B.2: Nondegenerate representations}

Since the introduction of Arveson systems in \cite{Arv89}, \it{representations} of product systems have been recognized as an important concept. For product systems of correspondences the definition is due to Muhly and Solel \cite{MuSo02}. A concept equivalent to \it{faithful} and \it{nondegenerate} representation is that of \it{right dilation}, introduced in Skeide \cite{Ske06,Ske07,Ske11a}. The naming, \it{left} and \it{right} dilation underlines a deep symmetry between the concepts. In fact, the \it{commutant} of von Neumann correspondences introduced in Skeide \cite{Ske03c} (and, independently, in Muhly and Solel \cite{MuSo04} under the name of \it{\nbd{\sigma}dual}) turns a left dilation of a product system into a right dilation of its commutant system, and \it{vice versa}; see Skeide \cite[Theorem 9.9]{Ske09} (preprint 2004).

While on the algebraic level, the correspondence between left and right dilations is perfect, this situation changes, when we take into consideration continuity. Under the hypothesis that all occurring \nbd{W^*}algebras and \nbd{W^*}modules have separable pre-dual, Muhly and Solel have proved in \cite{MuSo07} that the duality is perfect for \it{weakly measurable} product systems, by reducing it to a result by Effros \cite{Eff65}: Measurable fields of von Neumann algebras have measurable commutant fields.

Since we do not know a similar result for \it{strongly continuous} fields of von Neumann algebras, we cannot imitate the reduction from \cite{MuSo07}. Instead, we have to prove from scratch existence of strongly continuous right dilations under a condition, faithfulness of all left actions, which under commutant corresponds to strong fullness on the commutant side. In the case of strongly continuous product systems that are both strongly full and faithful, we get symmetry: Also the commutant system possesses a strongly continuous structure. Note, however, that this result relies heavily on the ``semigroup structure'' encoded by the product system; a product system is a monoid of correspondences. In fact, we shall use that also a strongly continuous right dilation is related with a certain \nbd{E_0}semigroup and that the product system of that \nbd{E_0}semigroup is the commutant of the one for which we construct a right dilation. If that right dilation is strongly continuous, then so is the \nbd{E_0}semigroup. Therefore, the commutant system inherits a strongly continuous structure.

In the sequel, we discuss only \nbd{W^*}versions. \nbd{C^*}Versions have been discussed in \cite{Ske09,Ske11a,Ske08}. So, $E^\odot$ is a product system of \nbd{W^*}correspondences over the \nbd{W^*}algebra $\cB$.

\bemp[Definition \cite{MuSo02}.~]
A \hl{representation} of $E^\odot$ on a Hilbert space $H$ is a family $\eta=\bfam{\eta_t}_{t\in\R_+}$ of linear maps $\eta_t\colon E_t\rightarrow\sB(H)$ that fulfill
\baln{
\eta_t(x_t)\eta_s(y_s)
&
~=~
\eta_{t+s}(x_ty_s),
&
\eta_t(x_t)^*\eta_t(y_t)
&
~=~
\eta_0(\AB{x_t,y_t}).
}\ealn
A representation is \hl{normal} (\hl{faithful}) if $\eta_0$ is normal (faithful). A representation is \hl{nondegenerate} if each $\eta_t$ acts nondegenerately on $H$ (that is, $\cls\eta_t(E_t)H=H$).
\eemp

Note that $\eta_t$ is a ternary homomorphism into the Hilbert \nbd{\sB(H)}module $\sB(H)$. Therefore, $\eta_t$ is completely contractive. $\eta_0$ is a representation. (Simply, put $t=s=0$.) In particular, it makes sense to speak of $\eta_0$ being normal (that is, order continuous). If $\eta$ is faithful, then each $\eta_t$ is faithful. If $\eta$ is normal, then each $\eta_t$ is \nbd{\sigma}weak; see the remark following \cite[Lemma 2.16]{MuSo02}. Note, too, that the pair $(\eta_t,\eta_0)$ is what is called a \it{representation} of the single correspondence $E_t$.

\brem
Nondegenerate representations have been called \it{essential} in Arveson \cite{Arv89} and Hirshberg \cite{Hir05a}. Muhly and Solel \cite{MuSo02} defined as a \it{covariant representation} a family that fulfills only $\eta_t(x_t)\eta_s(y_s)=\eta_{t+s}(x_ty_s)$ plus the requirement that $\eta_0$ is a representation (which is no longer automatic). A covariant representation also fulfilling $\eta_t(x_t)^*\eta_t(y_t)=\eta_0(\AB{x_t,y_t})$, is called \it{isometric}, while a covariant representation fulfilling the nondegeneracy condition, is called \it{fully coisometric}. Also, covariant representations are usually required to be completely contractive. In the isometric case (that is, in particular, for our representations), this is automatic.
\erem

\bemp[Definition \cite{Ske06,Ske11a}.~]
A \hl{right dilation} of $E^\odot$ to a faithful \nbd{W^*}correspondence $H$ from $\cB$ to $\C$ (that is, a Hilbert space with a nondegenerate normal faithful left action of $\cB$) is a family of bilinear unitaries $w_t\colon E_t\odot H\rightarrow H$, such that the product $(x_t,h)\mapsto x_th:=w_t(x_t\odot h)$ iterates associatively with the product system structure.
\eemp

Like left dilations to $E$ require that $E$ is strongly full, right dilations to $H$ require that $H$ is faithful. If either condition is missing, then we speak of left and right \it{quasi} dilations.

Note that existence of a right dilation $E^\odot$ implies that $E^\odot$ is \hl{faithful} in the sense that each $E_t$ has faithful left action.

Without the obvious proof, we state that faithful nondegenerate representations and right dilations are equivalent concepts.

\bprop\label{rep=rdprop}
The relation $\eta_t(x_t)h=w_t(x_t\odot h)$ establishes a one-to-one correspondence between faithful nondegenerate normal representations $\eta$ of $E^\odot$ on $H$ and right dilations $w$ of $E^\odot$ to $H$.
\eprop

\bemp[Facts.~]\label{rdfacts}
We collect some facts about the relation between right dilations of product systems of \nbd{W^*}correspondences, \nbd{E_0}semigroups and commutants of product systems.

\begin{enumerate}
\item\label{rdf1}
If $w_t$ is a right dilation of $E^\odot$ to $H$, then $\vt^w_t(a):=w_t(\id_t\odot a)w_t^*$ defines an \nbd{E_0}semigroup on $\sB^{bil}(H)$, the von Neumann algebra of bilinear operators on $H$; see \cite{MuSo02}.

\item\label{rdf2}
Conversely, if $H$ is a faithful \nbd{W^*}correspondence from $\cB$ to $\C$ (that is, a Hilbert space with a faithful, nondegenerate, and normal representation of $\cB$), and if $\vt$ is a normal \nbd{E_0}semigroup on $\sB^{bil}(H)$, then
\beqn{
E_t
~:=~
\BCB{x_t\in\sB(H)\colon\vt_t(a)x_t=x_ta~(a\in\sB^{bil}(H))}
}\eeqn
is a correspondence over $\cB$ with left and right action given by the left action of $\cB$ on $H$, and with $\AB{x_t,y_t}$ being that unique element in $\cB$ that acts on $H$ like $x_t^*y_t\in\sB^{bil}(H)'\cong\cB$. Moreover, $E^\odot=\bfam{E_t}_{t\in\R_+}$ becomes a product system with $u_{t,s}(x_t\odot y_s):=x_ty_s$, and $w_t(x_t\odot h):=x_th$ defines a right dilation such that $\vt^w=\vt$.

This product system has been constructed in \cite{Ske03c}. In the case $\cB=\C$ (so that $\sB^{bil}(H)=\sB(H)$), we recover Arveson's construction \cite{Arv89} of an Arveson system for an \nbd{E_0}semigroup on $\sB(H)$.

\item\label{rdf3}
Considering also $G$ as correspondence from $\cB$ to $\C$, the space $E':=\sB^{bil}(G,H)$ is a \nbd{W^*}module over $\cB'$, acting nondegenerately on $G$ in the sense that $E'G$ is total in $H$; see Rieffel \cite[Proposition 6.10]{Rie74a}, Muhly and Solel \cite[Lemma 2.10]{MuSo02}, and Skeide \cite{Ske03c,Ske05c}. Moreover, $\sB^a(E')=\sB^{bil}(H)$ and $E'$ is strongly full. So, the \nbd{E_0}semi\-group $\vt$ on this algebra also has a product system $E'^\odot$ of \nbd{W^*}correspondences over $\cB'$.

This product system $E'^\odot$ has been introduced in \cite{Ske03c} as the commutant system of $E^\odot$. In \cite{Ske03c} it is also explained that the two product systems used in \cite{BhSk00} and in \cite{MuSo02} to construct in two different ways the unique minimal weak dilation of a Markov semigroup on $\cB'$, actually are commutants of each other.

\end{enumerate}
\eemp

\brem
We emphasize that everything about the commutant works for an arbitrary faithful nondegenerate normal representation of $\cB$ on $G$. The commutant depends on that representation (up to Morita equivalence), but most results do not.
\erem

\bob
We mentioned in the beginning of Appendix B.2 that the commutant also takes left dilations to right dilations and \it{vice versa}. (In fact, the \nbd{\cB}\nbd{\C}correspondence $H$ under commutant goes precisely to the \nbd{\C}\nbd{\cB'}correspondence $E'$ (note: $\C'=\C$!) and \it{vice versa}. Also, strong fullness on one side corresponds to faithfulness on the other. Moreover, the commutant functor is anti-multiplicative for the tensor product. See \cite[Section 9]{Ske09}.) So, as far as existence of an algebraic right dilation is concerned, we are done, because we know about existence of left dilations.

However, the transition from left to right dilations (and back) via the commutant functor involves the construction of intertwiner spaces. But we do not know theorems asserting that bundles of intertwiner spaces have enough sufficiently continuous sections. (With this in mind, it is even more noteworthy, that Arveson, actually, proved in \cite[Lemma 2.3]{Arv89} that the intertwiner spaces for an \nbd{E_0}semigroup on $\sB(H)$ ($H$ separable!) do have enough continuous sections. We pose as an open problem, to see if that proof has a chance to be generalized.) This bad behavior is the reason, why we have to do continuous right dilations separately, instead of inferring them from continuous left dilations.
\eob

The proof of existence of a strongly continuous right dilation for a faithful strongly continuous product system of \nbd{W^*}correspondences is, a bit surprisingly, much more similar to that for the \nbd{C^*}case in \cite{Ske11a}, than the proof for left dilations in Appendix B.1 is to that in \cite{Ske07}. This is mainly so, because already in the \nbd{C^*}case we had to deal with sections in $CS_i(E^\odot\odot G)$, and the few places where substituting $C_i(E^\odot)$ with $C_i^s(E^\odot)$ causes some differences have already been dealt with in Appendix B.1.

In \cite{Ske11a}, the basic ingredient for the construction of a right dilation was granted by \cite[Theorem 1.2]{Ske11a}. In the introduction of \cite{Ske11a} we spent some time to explain why this theorem corresponds to the existence of a unit vector in $\ol{M_\en(E_1)}^s$ in the construction of a left dilation. Here, we just ask the reader to keep in mind the following consideration based on the duality via commutant: If in the construction of a left dilation of $E^\odot$ a unit vector in (a suitable multiple of) $E_1$ plays a crucial role, then in the construction of a right dilation of that product system this role should be played by a unit vector in (a suitable multiple of) the commutant of $E_1$, $E'_1$. This explains, roughly speaking, why we get as main ingredient an element in an intertwiner space and not in $E_1$ itself.

\bthm\label{Xi'thm}
There exists a cardinality $\en$ such that $\sB^{bil}(G^\en,E_1\odot G^\en)$ admits an isometry $\Xi'_1$.
\ethm

\brem
Note that $\sB^{bil}(G^\en,E_1\odot G^\en)$ is nothing but $\ol{M_\en(E'_1)}^s$.
\erem

\proof[Proof of Theorem \ref{Xi'thm}.~]
In principle, this is the \nbd{W^*}analogue of \cite[Theorem 1.2]{Ske11a}, and proved exactly in the same way from existence of a nondegenerate faithful representation of a correspondence $E_1$. Existence of that representation is granted, in the \nbd{W^*}case, by \cite[Theorem 8.2]{Ske09}, and \cite[Observation 8.5]{Ske09} tells us that we may even choose that representation in the desired form. (Actually, our proof of the \nbd{C^*}case in \cite[Theorem 8.3]{Ske09}, providing the generalization of Hirshberg's original result \cite{Hir05a} to the non-full case, is by a subtle reduction to the \nbd{W^*}case \cite[Theorem 8.2]{Ske09}.)\qed

\lf
Note that $E_t\odot G^\en$ is nothing but $H_t^\en$. In particular, $E_1\odot G^\en=H_1^\en$. By standard results on tensor products, $\Xi'_1$ gives rise to an operator $(u_{t,1}\odot\id_{G^\en})(\id_t\odot\Xi'_1)\in\sB^{bil}(E_t\odot G^\en,E_{t+1}\odot G^\en)=\sB^{bil}(H_t^\en,H_{t+1}^\en)$. By abuse of notation, we shall denote this operator by $\Xi'_1$, too. With this notation, we obtain the suggestive commutation relation $\Xi'_1x_t=x_t\Xi'_1$.

In order to avoid conflicts with the Hilbert space $G$, we shall denote elements in $G^\en$ by $g=\bfam{g_s}_{s\in S}$ (not by $G=\bfam{g_s}_{s\in S}$). We apply the same convention to other Hilbert spaces like $H^\en$ and $H_t^\en$.

It now is really important to observe that results like Proposition \ref{Hcontprop} do not depend on our choice of the representation of $\cB$ on $G$. Therefore, Proposition \ref{Hcontprop} remains valid if we replace $G$ with $G^\en$ (and the canonical representation of $\cB$ on $G^\en$), and $CS_i(E^\odot\odot G)$ with $CS_i^\en(E^\odot\odot G):=CS_i(E^\odot\odot G^\en)$, the space of continuous sections $h\colon t\mapsto h_t\in H_t^\en$. 

For a section $h$ of $E^\odot\odot G^\en$ we shall denote $h(t):=i_t^\en h_t$. Let $0\le a<b<\infty$. By $\int_a^b H_\alpha^\en\,d\alpha$ we understand the norm completion of the pre-Hilbert space that consists of continuous sections $h\in CS_i^\en(E^\odot\odot G)$ restricted to $\RO{a,b}$ with inner product
\beqn{
\AB{h,h'}_{\SB{a,b}}
~:=~
\int_a^b\AB{h_\alpha,h'_\alpha}\,d\alpha
~=~
\int_a^b\AB{h(\alpha),h'(\alpha)}\,d\alpha.
}\eeqn
As in \cite[Proposition 4.2]{Ske07} or Proposition \ref{srcllprop} we show:

\bprop\label{rcllprop}
$\int_a^b H_\alpha^\en\,d\alpha$ contains the space $\eR_{\RO{a,b}}$ of restrictions to $\RO{a,b}$ of those sections $h$ for which $t\mapsto h(t)$ is right continuous with finite jumps (by this we mean, in particular, that there exists a left limit) in finitely many points of $\RO{a,b}$, and bounded on $\RO{a,b}$, as a pre-Hilbert subspace.
\eprop

In Appendix B.1 we considered sections that were stable under multiplication with $\Xi_1$. Here we have to multiply with $\Xi'_1$. Let $\sS$ denote the subspace of all sections $h=\bfam{h_t}_{t\in\R_+}$ of $E^\odot\odot G^\en$ which are \hl{locally $\eR$}, that is, for every $0\le a<b<\infty$ the restriction of $h$ to $\RO{a,b}$ is in $\eR_{\RO{a,b}}$, and which are \hl{stable} with respect to the isometry $\Xi'_1$, that is, there exists an $\alpha_0\ge0$ such that
\beq{\label{stabsec}
\Xi'_1h_\alpha
~=~
h_{\alpha+1}
}\eeq
holds for all $\alpha\ge\alpha_0$. By $\sN$ we denote the subspace of all sections in $\sS$ which are eventually $0$, that is, of all sections $h\in\sS$ for which there exists an $\alpha_0\ge0$ such that $h_\alpha=0$ for all $\alpha\ge\alpha_0$. A straightforward verification shows that
\beqn{
\AB{h,h'}
~:=~
\lim_{m\to\infty}\int_m^{m+1}\AB{h(\alpha),h'(\alpha)}\,d\alpha
}\eeqn
defines a semiinner product on $\sS$ and that $\AB{h,h}=0$ if and only if $h\in\sN$. Actually, we have
\beqn{
\AB{h,h'}
~=~
\int_T^{T+1}\AB{h(\alpha),h'(\alpha)}\,d\alpha
}\eeqn
for all sufficiently large $T>0$; see \cite[Lemma 2.1]{Arv06}. So, $\sS/\sN$ becomes a pre-Hilbert space with inner product $\AB{h+\sN,h'+\sN}:=\AB{h,h'}$. By $H$ we denote its completion.

\bprop\label{shiftcor}
If $h$ is in $CS_i^\en(E^\odot\odot G)$, then the shifted section
\beqn{
t
~\longmapsto~
\begin{cases}
0&t<1
\\
\Xi'_1h_{t-1}&t\ge1
\end{cases}
}\eeqn
is continuous for $t\ge1$.
\eprop

\proof
(Like our proof of Proposition \ref{Hcontprop}, also this proof is a simplification compared with the proof of \cite[Corollary 2.4]{Ske11a}.) By Proposition \ref{Hcontprop} the elements in $\ls CS_i(E^\odot)\odot G^\en$ approximate $t\mapsto h_t$ locally uniformly. So, it is enough to show the statement for sections of the form $t\mapsto x_t\odot g$ $(x\in CS_i(E^\odot),g\in G^\en)$. Again by Proposition \ref{Hcontprop} there is a section $h\in CS_i^\en(E^\odot\odot G)$ such that $h_1=\Xi'_1g$. Once more, by Proposition \ref{Hcontprop}, the map
\beqn{
(t,s)
~\longmapsto~
(i_{t+s}u_{t,s}\odot\id_{G^\en})x_th_s
}\eeqn
is continuous. This holds \it{a fortiori} if we fix $s=1$.\qed

\lf
From Proposition \ref{shiftcor} we easily deduce the following analogue of \cite[Proposition 4.3]{Ske07} or Proposition \ref{denseprop}.

\bcor\label{densecor}
For every section $h$ and every $\alpha_0\ge0$ define the section $h^{\alpha_0}$ as
\beqn{
h^{\alpha_0}_\alpha
~:=~
\begin{cases}
0&\alpha<\alpha_0
\\
{\zeta'_1}^nh_{\alpha-n}&\alpha\in\RO{\alpha_0+n,\alpha_0+n+1},n\in\N_0.
\end{cases}
}\eeqn
If $h$ is in $CS_i^\en(E^\odot\odot G)$, then $h^{\alpha_0}$ is in $\sS$. Moreover, the set $\bCB{h^{\alpha_0}+\sN\colon h\in CS_i(E^\odot\odot G),\alpha_0\ge0}$ is a dense subspace of $H$.
\ecor

Observe that on $H$ we have a canonical representation of $\cB$ that acts simply pointwise on sections. By a simple application of continuity, we see that this representation is faithful, because $E^\odot$ is faithful. Of course, $\U\in\cB$ acts as identity, so that the representation is nondegenerate.

\bprop
The canonical action of $\cB$ on $H$ is normal. (Equivalently, $\int_a^b E_\alpha^\en\,d\alpha$ is a \nbd{W^*}cor\-re\-spon\-dence!)
\eprop

\proof
This follows from strong density in Proposition \ref{denseprop} by an  application of Corollary \ref{uniapprcor}. (Note how important it is that Proposition \ref{denseprop} guaranties density and not just totality.)\qed

\lf
It is now completely plain to see that for every $t\in\R_+$ the map $x_t\odot h\mapsto x_th$, where
\beqn{
(x_th)_\alpha
~=~
\begin{cases}
x_th_{\alpha-t}&\alpha\ge t,
\\
0&\text{else},
\end{cases}
}\eeqn
defines an isometry $w_t\colon E_t\odot H\rightarrow H$, and that these isometries iterate associatively as required for a right dilation.

\bprop\label{rsurprop}
Each $w_t$ is surjective.
\eprop

\proof
The proof goes presicely like that of Proposition \ref{surprop}, just that in all expressions that concern $z_t$ and $y$ the oder must be inverted.\qed

\bprop\label{cdprop}
The $w_t$ are \hl{strongly continuous} in the sense that for every section $x\in CS_i^s(E^\odot)$ and every $h\in H$ the function $t\mapsto x_th$ is continuous.
\eprop

\proof
The proof goes like that of Proposition \ref{scdprop}. (Expressions like $X^{\alpha_0}y_t$ or $X^{\alpha_0}_\alpha y_t$ must be replaced with expressions like $y_t\odot x^{\alpha_0}$ or $y_tx^{\alpha_0}_\alpha$. Actually, the proof here is quite a bit simpler, because now the preparation in the first two paragraphs of that proof is no longer necessary, and we may immediately start showing the continuity statement only for sections $x^{\alpha_0}$.)\qed

\lf
We summarize:

\bthm\label{rdilthm}
Let $E^\odot$ be a faithful strongly continuous product system of \nbd{W^*}cor\-re\-spond\-ences. Then $E^\odot$ admits a strongly continuous right dilation.
\ethm

\brem
Also here the construction preserves countability hypotheses: If $\cB$ has separable pre-dual and if $E^\odot$ is countably generated, then $H$ is separable.
\erem

We mentioned already the following consequence of Theorem \ref{rdilthm}.

\bthm\label{scrdthm}
Let $E^\odot$ be a faithful strongly continuous product system of \nbd{W^*}cor\-re\-spond\-ences over $\cB\subset\sB(G)$ and let $w_t\colon E_t\odot H\rightarrow H$ be a strongly continuous right dilation, so that $E':=\sB^{bil}(G,H)$ is a von Neumann \nbd{\cB'}module with $\sB^a(E')=\sB^{bil}(H)$ (see \cite{Ske05c}) and the product system $E'^\odot$ of the \nbd{E_0}semigroup $\vt^w$ on $\sB^a(E')$ is the commutant system of $E^\odot$ (see \cite{Ske03c}). Since $\vt^w$ is strongly continuous, the commutant $E'^\odot$ possesses a strongly continuous structure.
\ethm

\bthm\label{scpsdualthm}
The commutant is a duality between strongly full and faithful strongly continuous product systems of correspondences over $\cB$ and strongly full and faithful strongly continuous product systems of correspondences over $\cB'$. 
\ethm

\noindent
\sc{Proof.~}
(Sketch. We were not very specific about the commutant and about the continuous structure under Morita equivalence. A detailed version will be provided in \cite{Ske07p}, where we intend to present a full duality via commutant for arbitrary strongly continuous product systems. But the following indications are sufficient to produce a detailed proof.)

Applying the same procedure as in Theorem \ref{scrdthm} to $E'^\odot$ (constructing a strongly continuous right dilation of $E'^\odot$), gives back the algebraic structure of $E^\odot$ (see \cite{Ske03c}) with some strongly continuous structure. The decisive questions that separate Theorem \ref{scpsdualthm} from Theorem \ref{scrdthm}, are: Does the continuous structure of $E'^\odot$ depend on the right dilation? Does the second iteration give back the continuous structure of $E^\odot$ we started with?

\lf
(i)~ Let us assume we have two strongly continuous right dilations $w^i_t$ of $E^\odot$ to $H_i$ ($i=1,2$). We shall illustrate that the \nbd{E_0}semigroups $\vt^i:=\vt^{w^i}$ are stably unitary cocycle inner conjugate via a strongly continuous cocycle, that is, by Theorem \ref{vNalgsclthm} the two strongly continuous structures induced on $E'^\odot$ coincide. One easily checks:
\begin{itemize}
\item
After stabilizing suitably, we get right dilations $(w^i_t)^\en\colon E_t\odot H_i^\en$ where $H_1$ and $H_2$ are isomorphic as \nbd{W^*}correspondences from $\cB$ to $\C$ (and both to $G^\en$).

\item
And if there is a bilinear unitary $H^1$ to $H^2$ then this allows to `lift' the right dilation $w^2_t$ to $H_2$ to a right dilation to $H_1$ generating an \nbd{E_0}semigroup inner conjugate to $\vt^2$ on $\sB^{bil}(H_1)$.
\end{itemize}
Therefore, after stabilizing (which does not change the continuous structure induced on $E'^\odot$) and identifying the \nbd{\cB}\nbd{\C}modules $H_1^\en$ and $H_2^\en$ (which amounts to an inner conjugacy and also does not change the induced continuous structure), we get two right dilations $w^i_t$ to the same \nbd{\cB}\nbd{\C}module $H\cong G^\en$.

Now it is easy to check that:
\begin{itemize}
\item
$u_t\colon w^1_t(x_t\odot h)\mapsto w^2_t(x_t\odot h)$ defines a unitary left cocycle for $\vt^1$ such that $(\vt^1)^u=\vt^2$.

\item
$u_t$ is strongly continuous. ($u_sw^1_t(x_t\odot h)\approx u_sw^1_s(x_s\odot h)=w^2_s(x_s\odot h)\approx w^2_t(x_t\odot h)=u_tw^1_t(x_t\odot h)$.)
\end{itemize}
So, indeed, the continuous structure induced on $E'^\odot$ does not depend on the choice of the right dilation of $E^\odot$.

\lf
(ii)~ To obtain that the induced strongly continuous structure on $E^\odot$ coming from some strongly continuous right dilation of $E'^\odot$ is the same as the original one, by Part (i) it is enough to show that also the original one is induced from some right dilation of $E'^\odot$. If we assume (as, by Theorem \ref{vNexithm}, we always may do) that the original strongly continuous structure of $E^\odot$ is derived from a strongly continuous \nbd{E_0}semigroup $\vt$ on some $\sB^a(E)$, then the preceding argument means that it is enough to recognize that $\vt$ is induced by a strongly continuous representation $\eta'_t$ of $E'^\odot$.

To that goal, it is important to write down $E'^\odot$ in a form that makes it easy to understand what are its strongly continuous sections. We now make full use of the fact that we may choose the identification of $\cB\subset\sB(G)$ to our liking. Changing that identification, means changing $\cB'$ and, therefore, $E'^\odot$ by a Morita equivalence. The relation among strongly continuous structures are precisely those in Remark \ref{MecPSrem} that make Morita equivalence an equivalence of (strongly) continuous product systems. Recall also \cite[Theorem 5.12]{Ske09}, which asserts that a product system of \nbd{W^*}correspondences is the one associated with an \nbd{E_0}semigroup if and only if it is Morita equivalent to a one-dimensional product system---namely, (for instance) the one-dimensional product system of the \nbd{E_0}semigroup. So, not having chosen $G$ yet, we simply choose the $G$ of a strongly continuous right dilation $w_t\colon E_t\odot G\rightarrow G$ (which exists by Theorem \ref{rdilthm}). Then $E'=\sB^{bil}(G)=\cB'=\sB^a(E')$ and the product system $E'^\odot$ of $\vt^w_t\colon\cB'\rightarrow\cB'$ is just the one-dimensional product system $E'_t:={_{\vt^w_t}}\cB'$; see the \it{methodological introduction} in Section \ref{intro}. Clearly, the strongly continuous sections of $E'^\odot$ induced by that $\vt^w$ are just the strongly continuous functions with values in $\cB'=E'_t$.

Now, construct the (strongly continuous) left dilation $v_t\colon E\sodots E_t\rightarrow E$ of $E^\odot$  to $E$ that associates $\vt$ with $E^\odot$ via $\vt=\vt^v$. By the proof of Theorem \ref{ArKiW*thm} (below), the unitaries $u_t:=(v_t\id_G)(\id_E\odot w_t^*)$ form a strongly continuous unitary semigroup in $\sB(H)$, where $H:=E\odot G$. (Observe that $E\sodots E_t\sodots G=E\sodots E_t\odot G=E\odot E_t\odot G$. This is so for all multiple tensor products of von Neumann correspondences where the right factor is a von Neumann module over a finite-dimensional algebra.) Directly from that definition we get $u_t(a\odot\id_G)u_t^*=\vt^v_t(a)\odot\id_G$ ($a\in\sB^a(E)$) and $u_t^*(\id_E\odot b')u_t=\id_E\odot\vt^w_t(b')$ ($b'\in\cB'$).

Now recall that $E'_t=\cB'\subset\sB(G)$ as sets (well, actually as von Neumann \nbd{\cB'}module). This means, every element $x'_t\in E'_t$ gives rise to an operator $\id_E\odot x'_t\in\sB(H)$. We claim that
\beqn{
\eta'_t(x'_t)
~:=~
u_t(\id_E\odot x'_t)
}\eeqn
defines a (faithful normal nondegenerate) representation of $E'^\odot$ on $H$ (so that, by Proposition \ref{rep=rdprop}, $w'_t\colon x'_t\odot h\mapsto\eta'_t(x'_t)h$ defines a right dilation). The verification of this statement is tedious but straightforward. (It is obvious that $\eta'_t(x'_t)^*\eta'_t(y'_t)=\eta'_0(x'^*_ty'_t)=\id_E\odot\AB{x'_t,y'_t}$ is just the action of $\AB{x'_t,y'_t}\in\cB'$ on $H$. For verifying that $\eta'_t$ is fulfills the representation property, it is convenient to denote (similar to the left action $b'.x'_t:=\vt^w_t(b')x'_t$ of $b'\in\cB'$ on $x'_t\in E'_t$) the product of $E'^\odot$ by $x'_t.y'_s:=u'_{t,s}(x'_t\odot y'_s):=\vt^w_s(x'_t)y'_s$. Also the associativity properties
\beqn{
v_t(v_s\odot\id_t)
~=~
v_{t+s}(\id_E\odot u_{s,t})
\text{~~~~~~~~~and~~~~~~~~~}
(\id_t\odot w_s^*)w_t^*
~=~
(u_{s,t}^*\odot\id_G)w_{t+s}^*
}\eeqn
play a role.) Since for each strongly continuous section $\bfam{x'_t}$ of $E'^\odot$ and each $g\in G$ the function $t\mapsto x'_tg\in G$ is continuous, and since the unitary semigroup $u_t$ is strongly continuous, the representation $\eta'_t$ (and, therefore, also the right dilations $w'_t$) is strongly continuous. Like for any right dilation of any product system, the action of the induced semigroup on $\sB^{bil}(H)=\sB^a(E)\odot\id_G=\sB^a(E)$ can be expresses in terms of the representation of that product system as
\beqn{
\vt^{w'}(a\odot\id_G)\eta'_t(x'_t)h
~=~
\eta'_t(x'_t)(a\odot\id_G)h,
}\eeqn
(This expresses the fact that elements of the right dilated product system are identified with the intertwiners of the induced \nbd{E_0}semigroup; see \cite{Ske03c}.) Concretely,
\bmun{
\vt^{w'}(a\odot\id_G)\eta'_t(x'_t)(x\odot g)
~=~
\eta'_t(x'_t)(a\odot\id_G)(x\odot g)
~=~
\eta'_t(x'_t)(ax\odot g)
~=~
u_t(ax\odot x'_tg)
\\
~=~
\vt_t(a\odot\id)u_t(x\odot x'_tg)
~=~
\vt_t(a\odot\id)\eta'_t(x'_t)(x\odot g),
}\emun
where in the step from the first line to to the second line we made use of $u_t(a\odot\id)=\vt_t(a\odot\id)u_t$. In other words, $\vt=\vt^{w'}$, so that the continuous structure induced on $E^\odot$ by $w'_t$ (and, therefore, by Part (i), by every strongly continuous right dilation of $E'^\odot$) coincides with the with original continuous structure of $E^\odot$ (induced by the \nbd{E_0}semigroup $\vt$).~\qedsymbol

\brem
The theorem holds, in particular, for spatial strongly continuous product systems, which are strongly full and faithful, automatically, because they contain the trivial product system, which is strongly full and faithful.
\erem

Precisely as the \nbd{C^*}case in \cite[Theorem 3.1]{Ske11a}, we provide a completely different proof of the following result due to Arveson and Kishimoto \cite{ArKi92}.

\bthm\label{ArKiW*thm}
Let $E$ be a strongly full \nbd{W^*}module over a \nbd{W^*}algebra $\cB$ (for instance, let $E=\cB$ itself!) and let $\vt$ be a faithful strongly continuous normal \nbd{E_0}semigroup on $\sB^a(E)$ (for instance $\cB$ itself if $E=\cB$). Then there exists a faithful \nbd{W^*}correspondence $K$ from $\sB^a(E)$ to $\C$ with strict left action (that is, a Hilbert space with a faithful nondegenerate normal representation of \,$\sB^a(E)$) and a strongly continuous unitary group $u$ on $K$ such that $\vt_t(a)k=u_tau_t^*k$ for all $a\in\sB^a(E)$, $t\in\R_+$, $k\in K$.
\ethm

\proof
The idea is already described in the proof of Theorem \ref{scpsdualthm}: If $v_t\colon E\sodots E_t\rightarrow E$ it a left dilation of a product system $E^\odot$ (so that $E^\odot$ is necessarily strongly full) and if $w_t\colon E_t\odot H\rightarrow H$ it a right dilation of a product system $E^\odot$ (so that $E^\odot$ is necessarily faithful), then one easily verifies that $u_:=(v_t\odot\id_K)(\id_E\odot w_t^*)$ defines a unitary semigroup in $\sB(K)$ with $K:=E\odot H$ (which may be extended to a group) and that $u_t(a\odot\id_H)u_t^*=(v_t(a\odot\id_t)v_t^*)\odot\id_H$. The proof that $u_t$ is continuous provided $v_t$ and $w_t$ are strongly continuous, goes exactly like that of the \nbd{C^*}case in \cite[Theorem 3.1]{Ske11a}.

Now if $\vt$ is a strongly continuous normal \nbd{E_0}semigroup on $\sB^a(E)$, then construct the associated product system $E^\odot$ and left dilation $v_t$. If $\vt$ is faithful, so is $E^\odot$. Therefore, we also get a right dilation $w_t$, and the unitary group $u$ does the job.\qed

\lf
An \hl{elementary dilation} of a CP-semigroup $T$ on $\cB$ is a \nbd{C^*}algebra $\cA$ with an embedding $\vp\colon\cB\rightarrow\cA$ and a semigroup $c=\bfam{c_t}_{t\in\R_+}$ of elements in $\cA$ such that
\beqn{
\vp\circ T(b)
~=~
c^*_t\vp(b)c_t
}\eeqn
for all $b\in\cB$ and $t\in\R_+$. A CP-semigroup $T$ is \hl{semifaithful} if its GNS-system is faithful. (From this, it follows that an arbitrary weak dilation $\vt$ of $T$ is faithful.) Putting together Theorems \ref{scrdthm} and \ref{scdilthm}, we get the following \nbd{W^*}analogue of \cite[Theorem 3.4]{Ske11a}. (Note that the strong continuity condition here is referring to the strong operator topology, and is much weaker than the one in the \nbd{C^*}case.)

\bthm\label{W*elemdilthm}
Every semifaithful strongly continuous normal CP-semigroup on a \nbd{W^*}al\-ge\-bra admits a (strongly continuous) elementary dilation to some $\sB(H)$ with normal embedding $\cB\rightarrow\sB(H)$.
\ethm


\addcontentsline{toc}{section}{References}

\setlength{\baselineskip}{2.5ex}


\newcommand{\Swap}[2]{#2#1}\newcommand{\Sort}[1]{}
\providecommand{\bysame}{\leavevmode\hbox to3em{\hrulefill}\thinspace}
\providecommand{\MR}{\relax\ifhmode\unskip\space\fi MR }
\providecommand{\MRhref}[2]{%
  \href{http://www.ams.org/mathscinet-getitem?mr=#1}{#2}
}
\providecommand{\href}[2]{#2}

\lf\noindent
Michael Skeide:
{\small\itshape Dipartimento E.G.S.I., Universit\`a\ degli Studi del Molise, Via de Sanctis, 86100 Campobasso, Italy, E-mail: \href{mailto:skeide@unimol.it}{\tt{skeide@unimol.it}}}\\
{\small{\itshape Homepage: \url{http://web.unimol.it/skeide/}}}


\end{document}